\documentclass[a4paper,11pt,final]{amsart}

\usepackage[utf8]{inputenc}
\usepackage[T1]{fontenc}
\usepackage[english]{babel}
\usepackage[margin=2cm]{geometry}
\usepackage{amsmath}
\usepackage{amsthm}
\usepackage{amssymb}
\usepackage{dsfont}
\usepackage[all]{xy}
\usepackage{wasysym}
\usepackage{mathtools}
\usepackage{bbold}
\usepackage{empheq}
\usepackage{tensor}
\usepackage{bbm}
\usepackage{enumitem}
\usepackage{rotating}
\usepackage{multirow}
\usepackage{setspace}
\usepackage{fancybox}
\usepackage{xargs}
\usepackage{indentfirst}
\usepackage{xcolor}
\definecolor{InternalLinks}{rgb}{0.33, 0.29, 0.31}
\usepackage{hyperref}
\usepackage[final]{showkeys}
\usepackage[colorinlistoftodos,prependcaption,textsize=tiny]{todonotes}
\usepackage{draftwatermark}
\SetWatermarkScale{5}
\SetWatermarkColor[gray]{0.9}
\usepackage[final]{fixme}
\usepackage{comment} 
\newif\ifshow 
\showtrue 
\ifshow
  \includecomment{wrap}
\else
  \excludecomment{wrap} 
\fi

\theoremstyle{plain}
\newtheorem{theo}{Theorem}[section]
\newtheorem*{theo*}{Theorem}
\newtheorem{lemm}[theo]{Lemma}
\newtheorem{prop}[theo]{Proposition}
\newtheorem{coro}[theo]{Corollary}

\theoremstyle{definition}
\newtheorem{defi}[theo]{Definition}
\newtheorem{exam}[theo]{Example}

\theoremstyle{remark}
\newtheorem{rema}[theo]{Remark}
\newtheorem*{rema*}{Remark}
\newtheorem{nota}[theo]{Notation}
\newtheorem*{nota*}{Notation}

\numberwithin{equation}{section}

\newcommand{\pr}{\noindent{\it Proof}\quad}
\newcommand{\fin}{\qed\smallskip}

\newcommand{\C}{\mathbb{C}}

\newcommand{\ku}{\Bbbk}
\newcommand{\id}{\mathrm{id}}

\newcommand{\W}{\mathbb{W}}
\newcommand{\V}{\mathbb{V}}
\newcommand{\U}{\mathbb{U}}
\newcommand{\hb}[1]{\hat{\bar{#1}}}

\newcommand{\G}{\mathbb{G}}

\newcommand{\dG}{\du{\mathbb{G}}}
\newcommand{\He}{\mathcal{H}}

\newcommand{\aYD}{\mathfrak{YD}}
\newcommand{\alg}{\mathfrak{Alg}}
\newcommand{\D}{\mathcal{D}}
\newcommand{\T}{\mathcal{T}}
\newcommand{\p}{\mathbf{p}}
\newcommand{\pol}{\mathcal{O}}

\newcommand{\m}{\mathfrak{m}}
\newcommand{\M}{\mathrm{M}}
\newcommand{\du}[1]{\widehat{#1}}
\newcommand{\com}{\Delta}
\newcommand{\dcom}{\du{\Delta}}

\newcommand{\cou}{\varepsilon}

\newcommand{\var}{\varphi}
\newcommand{\dvar}{\du{\varphi}}

\newcommand{\op}{\mathrm{op}}
\newcommand{\co}{\mathrm{co}}
\newcommand{\ops}{\circ}
\newcommand{\con}{\mathrm{c}}

\newcommand{\alp}{\alpha}

\newcommand{\odo}{\otimes}
\newcommand{\oda}{\otimes}

\newcommand{\bicroi}{\bowtie}

\newcommand{\blhd}{\blacktriangleleft}
\newcommand{\brhd}{\blacktriangleright}

\newcommand{\ome}{\omega}
\newcommand{\ad}{\mathrm{Ad}}
\newcommand{\as}{$*$}

\newcommand{\tr}{\text{trv}}

\newcommand{\ida}{\Rightarrow}
\newcommand{\vuelta}{\Leftarrow}

\newcommand{\fle}[1]{\xymatrix@C=3pc@R=3pc{#1}}

\newcommand{\iso}{\cong}

\newcommand\restr[2]{{
  \left.\kern-\nulldelimiterspace 
  #1 
  \vphantom{\big|} 
  \right|_{#2} 
}}

\newcommand{\gpd}[3]{\xymatrix{\mathcal{{#1}} \ar@<0.5ex>[r]^-{{#2}}\ar@<-0.5ex>[r]_-{{#3}} & \mathcal{{#1}}^{(0)}}}

\definecolor{ToDo}{RGB}{30,144,255}
\definecolor{Question}{RGB}{220,20,60}
\definecolor{Attention}{RGB}{255,215,0}

\newcommandx{\td}[2][1=]{\hfill\todo[inline,size=\footnotesize,linecolor=ToDo,backgroundcolor=ToDo!40,bordercolor=black,#1]{{\bf ToDo:} #2}}
\newcommandx{\qs}[2][1=]{\hfill\todo[inline,size=\normalsize,linecolor=Question,backgroundcolor=Question!40,bordercolor=black,#1]{{\bf Question:} #2}}
\newcommandx{\rmq}[2][1=]{\hfill\todo[inline,size=\normalsize,linecolor=Attention,backgroundcolor=Attention!40,bordercolor=black,#1]{{\bf Remark:} #2}}

\makeatletter
\newcommand{\addresseshere}{\enddoc@text\let\enddoc@text\relax}
\makeatother
\makeatletter
\newcommand{\pushright}[1]{\ifmeasuring@#1\else\omit\hfill$\displaystyle#1$\fi\ignorespaces}
\newcommand{\pushleft}[1]{\ifmeasuring@#1\else\omit$\displaystyle#1$\hfill\fi\ignorespaces}
\makeatother

\title[Yetter--Drinfeld algebras over a pairing of multiplier Hopf algebras]{Yetter--Drinfeld algebras over a pairing \protect\\ of multiplier Hopf algebras}
\author{Frank Taipe}
\address{Université Paris-Saclay, CNRS, Laboratoire de Mathématiques d'Orsay, 91405 Orsay, France}
\email{frank.taipe@universite-paris-saclay.fr}


\hypersetup{
pdfauthor={Frank Taipe},
pdftitle={Yetter-Drinfeld algebras over a pairing of multiplier Hopf algebras},
pdfsubject={Yetter-Drinfeld algebras in the framework of multiplier Hopf algebras},
pdfkeywords={Yetter-Drinfeld algebras, Multiplier Hopf algebras}
}

\begin{document}

\maketitle

\begin{abstract}
In the present work, we study Yetter--Drinfeld algebras over a pairing of multiplier Hopf algebras. Our main motivation is the construction of a self-dual theory of (C*-)algebraic quantum transformation groupoids. Instead of the standard characterization of Yetter--Drinfeld algebras given in the case of Hopf algebras, we develop an equivalent ``only coaction'' characterization in the framework of multiplier Hopf algebras. Finally, as a special case, we focus on Yetter--Drinfeld structures over Van Daele's algebraic quantum groups.
\end{abstract}

\section{Introduction}\label{sec:introduction}

In order to study the interplay between knot invariants, Hopf algebras and monoidal categories, Yetter--Drinfeld modules over bialgebras were introduced by D. Yetter in \cite{Y90}, originally under the name of crossed bimodules over bialgebras. In \cite{M91}, S. Majid has observed that Yetter--Drinfeld modules over finite-dimensional Hopf algebras are in one-to-one correspondence with modules over Drinfeld double of finite-dimensional Hopf algebras. Thereafter, Yetter--Drinfeld modules were used by D. Radford \cite{R93}, under the name of quantum Yang--Baxter modules to found new solutions to the quantum Yang--Baxter equation. In \cite{LR92}, L. Lambe and D. Radford also proved that those kind of modules over Hopf algebras are in one-to-one correspondence with modules over the Drinfeld double of Hopf algebras, regardless the finite-dimension of the Hopf algebra. The same result was found by M. Takeuchi in \cite{T92} in the operator algebra setting. From a modern categorical point of view, Yetter--Drinfeld modules over a Hopf algebra are in one-to-one correspondence with modules belonging to the Drinfeld center of the category of modules of the Hopf algebra.

A special class of Yetter--Drinfeld modules over a Hopf algebra is that one of braided commutative Yetter--Drinfeld algebras. This kind of objects were used by J.-H. Lu as main ingredient to construct non-trivial examples in her theory of Hopf algebroids \cite{L96}. These non-trivial Hopf algebroids can be considered as quantized versions of transformation groupoids. To illustrate this idea, consider an action of a finite group $G$ on a finite set $S$, hence two canonical examples arise:
\begin{enumerate}[label=(E\arabic*)]
\item The commutative algebra $\ku^{S}$ of $\ku$-valued functions on the set $S$ yields a braided commutative Yetter--Drinfeld algebra structure over the finite-dimensional commutative Hopf algebra $(\ku^{G},\com)$ of $\ku$-valued functions on the group $G$. By Lu's construction, we obtain the Hopf algebroid of $\ku$-valued functions on the transformation groupoid $G \ltimes S$ over the base algebra of $\ku$-valued functions on the set $S$.

\item The commutative algebra $\ku^{S}$ yields a braided commutative Yetter--Drinfeld algebra structure over the finite-dimensional group Hopf algebra $(\ku[G],\com)$. By Lu's construction, we obtain the groupoid Hopf algebroid of the transformation groupoid $G \ltimes S$ over the base algebra of $\ku$-valued functions on the set $S$.
\end{enumerate}

A generalization of the link between Yetter--Drinfeld structures and Hopf algebroids was done by T. Berzezi{\'n}ski and G. Milataru in \cite{BM02}, where they gave the conditions to characterize the Hopf algebroids arising from braided commutative Yetter--Drinfeld algebras over Hopf algebras. A piquant example using this nice construction of Hopf algebroids was obtain by A. Semikhatov in \cite{Se11}, where he showed the existence of a Yetter--Drinfeld structure on the Heisenberg algebra over the Drinfeld double of a Hopf algebra. This exotic Yetter--Drinfeld structure is related to the Kazhdan--Lusztig duality between logarithmic conformal field theories and Drinfeld--Jimbo's quantum groups \cite{Se10}.

Multiplier Hopf algebras were introduced by A. Van Daele in \cite{VD94} as a natural generalization of Hopf algebras. Subsequently, the self-dual theory of multiplier Hopf algebras with integrals \cite{VD98}, the theory of actions and coactions of multiplier Hopf algebras \cite{DVDZ99,VDZ99}, the corepresentation theory of multiplier Hopf algebras \cite{VDZ99_I,KVDZ00}, the construction of twisted objects arising from multiplier Hopf algebras \cite{D03,D04} and the construction of the Drinfeld double and the Heisenberg algebra of a multiplier Hopf algebras \cite{DVD01,DVD04_1,DVD04_2}, haven shown the fact that if a result is true in the framework of Hopf algebras then it can also be true in the case of multiplier Hopf algebras. Following this philosophy, in the direction of Yetter--Drinfeld structures, some important results were satisfactorily extended to the setting of multiplier Hopf algebras by L. Delvaux \cite{D05,D07}. A recent generalization of the Semikhatov's construction to the multiplier Hopf algebra setting can be found in \cite{YZC17}.

\medskip

The aim of this work is two-fold, first we want to continue the study of Yetter--Drinfeld structures in the framework of multiplier Hopf algebras and second we want to develop a modern characterization of Yetter--Drinfeld structures which is more suitable to work with in the self-dual theory of algebraic quantum transformation groupoids \cite{Ta22_3}. This modern approach of Yetter--Drinfeld structures is an algebraic adaptation of the characterization of Yetter--Drinfeld structures used in the setting of C*-algebraic quantum groups \cite{NV10}. An algebraic quantum transformation groupoid will be consider as the purely algebraic counterpart of a measured quantum transformation groupoid \cite{ET16}.

\medskip

This paper is organized as follows. In Section \ref{sec:preli}, we fix notations and we recall some important facts of the theory of actions and coactions of regular multiplier Hopf algebras. We also recall the standard characterization of Yetter--Drinfeld algebras in the framework of regular multiplier Hopf algebras. Then, in Section \ref{sec:oc_yd}, we introduce the {\em ``only coaction'' characterization} of Yetter--Drinfeld algebras over a pairing of multiplier Hopf algebras and we establish the equivalence with the standard characterization. We give the relationship between the ``only coaction'' characterization over a pairing and coactions of the Drinfeld codouble associated with the pairing. In Section \ref{sec:exa}, we present some examples of Yetter--Drinfeld algebras over a paring of multiplier Hopf algebras. And to illustrate the advantage of the ``only coaction'' characterization in the setting of multiplier Hopf algebras, we give another proof of the Yetter--Drinfeld structure on the Heisenberg algebra over the Drinfeld double without the use of the {\em Van Daele's covering technique}. Finally, in Section \ref{sec:yd-quantum}, we focus on the special case of Yetter--Drinfeld \as-algebras over an algebraic quantum group. 

\subsection*{Notations and conventions} An algebra in this work means a not necessary unital $\ku$-algebra where $\ku$ is a field. We only work with non-degenerate and idempotent algebras. Keeping that in mind, for a given non-degenerate idempotent algebra $A$, we can consider its multiplier algebra $\M(A)$ (see below for the definition), which is a unital idempotent algebra. We denote by ${}_{A}M$ a left $A$-module $M$ and by $\rhd : A \odo M \to M$ its module map. The multiplication map on an algebra $A$ is always denoted by $\m_{A}: A \odo A \to A$. We write ${}_{A}A$ (or $A_A$) when we regard $A$ as a left (or right) module over itself with respect to the left (or right) module map $\m_{A}$, respectively. Let $A, B$ be two algebras. For any invertible element $W$ in $A \odo B$, we have an invertible element $\ad(W) \in \textrm{End}(A \odo B)$ defined by
$\ad(W)(a \odo b) = W(a \odo b)W^{-1}$ for all $a \in A$, $b \in B$. Its inverse map is given by $\ad(W^{-1})$. We use the notation $\Sigma_{A,B}$ for the flip map $a \odo b \in A \odo B \mapsto b \odo a \in B \odo A$. Let $A, B, C$ be algebras. Through this work, given an element $E \in A \odo B$, we use the following legs notations
\[
E_{12} = E \odo 1_{C}, \quad E_{23} = 1_{C} \odo E, \quad E_{13} = (\id_{A} \odo \Sigma_{B,C})(E_{12}) = (\Sigma_{C,A} \odo \id_{B})(E_{23})
\]
for the elements in $A \odo B \odo C$, $C \odo A \odo B$ and $A \odo C \odo B$, respectively. For any algebra $A$, its opposite algebra, denoted by $A^{\op}$, is the $\ku$-module $A$ endowed with the multiplication map given by $\m_{A^{\op}}: A^{\op} \times A^{\op} \to A^{\op}$, $(a^{\op},a'^{\op}) \mapsto (a'a)^{\op}$. We use the superscript ${}^{\op}$ on an algebra element to indicate that we are seeing the element as an element of the opposite algebra. Sometimes, ${}^{\op}$ denotes the identity $\ku$-module map $\id_{A}: A^{\op} \to A$ and $\id_{A}: A \to A^{\op}$. 

\subsection*{Acknowledgments} An early version of the present work can be found in the author's PhD thesis. The author would like to express his thanks to Leonid Vainerman for the very fruitful discussions and Alfons Van Daele for his remarks and suggestions.

\section{Preliminaries}\label{sec:preli}

\subsection{Basics on non-unital algebras and its modules}

Let $A$ be a algebra.  A left module ${}_{A}M$ is called {\em faithful}, if for each non-zero $a \in A$ there exists $m \in M$ such that $a \rhd m$ is non-zero; {\em non-degenerate}, if for each non-zero $m \in M$ there exits $a \in A$ such that $a \rhd m$ is non-zero; {\em unital}, if $A \rhd M := \text{span}\{a \rhd m \,:\, a \in A, m \in M\} = M$. 
We say that ${}_{A}M$ {\em admits local units in $A$} if for every finite subset $F \subset M$ there exists $a \in A$ with $a \rhd m = m$ for all $m \in F$.

\begin{rema}
If ${}_{A}M$ admits local units, then it is non-degenerate and unital. The module ${}_{A}A$ (or $A{}_{A}$) is non-degenerate if and only if $A{}_{A}$ (or ${}_{A}A$) is faithful, respectively. For right modules, we use similar notations and terminology by identifying right $A$-modules with left $A^{\op}$-modules.
\end{rema}

An algebra $A$ is called {\em involutive or \as-algebra} if $A$ is a $\C$-algebra and there exists an anti-linear map $\phantom{\cdot}^{*} : A \to A$ called {\em the involution of $A$} such that $(a^{*})^{*} = a$ for all $a \in A$ and $(ab)^{*} = b^{*}a^{*}$ for all $a,b \in A$; {\em non-degenerate} if the modules ${}_{A}A$, $A_A$ both are non-degenerate; {\em idempotent}, if one of the modules ${}_{A}A$, $A_{A}$ is unital.
We say that $A$ {\em admits local units} if the modules ${}_{A}A$, $A_A$ both admit local units in $A$.

\begin{rema}\label{rem:sits}
The right module $A{}_{A}$ is non-degenerate if and only if the natural map $A \to \textrm{End}(A_A)$, $a \mapsto [L_{a} : b \mapsto ab]$ is injective. Similarly, the left module ${}_{A}A$ is non-degenerate if and only if the natural map $A \to \textrm{End}(\tensor[_{A}]A{})^{\op}$, $a \mapsto [R_{a} : b \mapsto ba]$ is injective.
\end{rema}


Let $A$ be a non-degenerate algebra, not necessarily unital. The {\em multiplier algebra of $A$} is the algebra
\[
\M(A) := \{ (L,R) \in \textrm{End}(A_A) \times \textrm{End}(\tensor[_{A}]A{})^{\op} \; : \; R(a)b = aL(b)\,\text{ for all }a,b\in A \}.
\]
with multiplication given by $(L,R)\cdot(L',R') := (L L',R' R)$, for all $(L,R),(L',R')\in \M(A)$. If $A$ is an \as-algebra, $\M(A)$ is also an \as-algebra with involution given by $(L,R)^{*} := (R^{\dagger},L^{\dagger})$ for all $(L,R) \in \M(A)$. Here, $R^{\dagger}(a) := R(a^{*})^{*}$ and $L^{\dagger}(a) := L(a^{*})^{*}$ for all $a\in A$.

\begin{rema}
It follows from Remark~\ref{rem:sits} that $A$ sits in $\M(A)$ as a sub-algebra with the identification $a \mapsto (L_{a},R_{a})$ for all $a \in A$. In particular, if $A$ admits a unity, then $\M(A) = A$.
\end{rema}

Given an element $T=(L,R) \in \M(A)$ and $a \in A$, we can identify $Ta$ with $L(a)$ and $aT$ with $R(a)$ as elements in $\M(A)$. An homomorphism of algebras $f: A \to M(B)$ is called {\em non-degenerate} if $f(A)B = B$ and $B f(A) = B$. In that case, the homomorphism $f: A \to \M(B)$ can be extended to $\M(A)$ and this extension, denoted by $f: \M(A) \to \M(B)$, is defined by $f(T)f(a) := f(Ta)$ and $f(a)f(T) := f(aT)$ for every $T \in \M(A)$ and $a \in A$.

\subsubsection{The opposite \as-algebra} Let $A$ be an \as-algebra and $\gamma$ be an algebra automorphism on $A$ such that $\gamma \circ * \circ \gamma \circ * = \id_{A}$. The anti-linear map $A^{\op} \to A^{\op}$, $a^{\op} \mapsto \gamma(a^{*})^{\op}$, defines an involution on $A^{\op}$. We will write $A^{\op}_{\gamma}$, to indicate the choice of the algebra $A^{\op}$ with the involution defined above. In the case $\gamma = \id_{A}$, the \as-algebra $A^{\op}_{\gamma}$ is the opposite algebra $A^{\op}$ with the trivial involution $(a^{\op})^{*} = a^{* \op}$ for all $a \in A$.

If $\hat\gamma$ denotes the inverse of $\gamma$, the map $\hat\gamma^{\ops}: A^{\op}_{\gamma} \to A^{\op}_{\gamma}$, defined by $\hat\gamma^{\ops}(a^{\op}) = \hat\gamma(a)^{\op}$, yields an algebra automorphism on $A^{\op}_{\gamma}$, with inverse $\gamma^{\ops}:a^{\op} \mapsto \gamma(a)^{\op}$, such that $\hat\gamma^{\ops} \circ * \circ \hat\gamma^{\ops} \circ * = \id_{A^{\op}_{\gamma}}$. The last implies that $(A^{\op}_{\gamma})^{\op}_{\hat\gamma^{\ops}} = A$ as \as-algebras.


\subsection{Basics on regular multiplier Hopf algebras} For further reading, we send the interested reader to the original paper \cite{VD94}. Let $A$ be a non-degenerate algebra. A {\em comultiplication} on $A$ is a homomorphism $\com : A \to \M(A \odo A)$ such that $\com(a)(1 \oda b)$ and $(a \oda 1)\com(b)$ belong to $A \odo A$ for all $a,b \in A$, and $\com$ is coassociative in the following sense
\[
(a \oda 1 \oda 1)(\com \odo \id)(\com(b)(1 \oda c)) = (\id \odo \com)((a \oda 1)\com(b))(1 \oda 1 \oda c)
\]
for all $a,b,c \in A$. We say that $(A,\com)$ is a {\em multiplier Hopf algebra} where $A$ is a non-degenerate algebra and $\com$ is a comultiplication on $A$ such that the linear maps from $A \odo A$ to itself, defined  by
\[
T_{1} : a \oda b \mapsto \com(a)(1 \oda b), \quad\quad T_{2} : a \oda b \mapsto (a \oda 1)\com(b)
\]
are bijective. A multiplier Hopf algebra $(A,\com)$ is called {\em regular} if $\com^{\co} := \Sigma_{A,A}  \com$ is again a comultiplication on $A$ such that $(A,\com^{\co})$ is also a multiplier Hopf algebra, called {\em the co-opposite multiplier Hopf algebra of $(A,\com)$}. Any unital Hopf algebra is a particular example of a regular multiplier Hopf algebra. The converse is also true, if $(A,\com)$ be a regular multiplier Hopf algebra, with $A$ being a unital algebra, then $(A,\com)$ is a Hopf algebra. In the case when $A$ is a \as-algebra, we call $\com$ a {\em comultiplication} on $A$ if it is also a \as-homomorphism. A multiplier Hopf \as-algebra is a \as-algebra with a comultiplication, making it into a multiplier Hopf algebra. For a multiplier Hopf \as-algebra, the regularity is automatic.

It was shown in \cite{VD94} that given a regular multiplier Hopf algebra $(A,\com)$, there exists a unique homomorphism $\cou : A \to \ku$, named {\it the counit}, satisfying
\[
(\cou \odo \id)(\com(a)(1 \oda b)) = ab \quad \text{ and } \quad (\id \odo \cou)((a \oda 1)\com(b)) = ab
\]
for all $a,b \in A$, and a unique anti-isomorphism $S :A \to A$, named {\it the antipode}, satisfying
\[
\m_{A}(S \odo \id)(\com(a)(1 \oda b)) = \cou(a)b \quad \text{ and } \quad \m_{A}(\id \odo S)((a \oda 1)\com(b))= \cou(b)a
\]
for all $a,b \in A$. When $A$ is a \as-algebra, we also have that $\cou$ is a \as-homomorphism and $S$ is an anti-isomorphism such that $S\circ * \circ S \circ * = \id_{A}$.

The following examples are the canonical multiplier Hopf \as-algebras arising from a not necessary finite group.

\begin{exam}\label{example_group}
Let $G$ be a group and let $K(G)$ be the (\as-)algebra of $\ku$-valued finitely supported functions on $G$. In this case, $\M(K(G))$ is the (\as-)algebra of all $\ku$-valued functions on $G$, $K(G) \odo K(G)$ can be naturally identified with the (\as-)algebra of $\ku$-valued finitely supported functions on $G \times G$, and $\M(K(G) \odo K(G))$ can be identified with the space of all $\ku$-valued functions on $G \times G$. We define the (\as-)homomorphism $\com: K(G) \to \M(K(G) \odo K(G))$, $p \mapsto [\com(p):(g,h) \mapsto p(gh)]$. For all $p, q \in K(G)$, the functions
\[
\begin{array}{rccc}
\com(p)(1 \oda q)  : & G \times G & \to & \ku \\
& (g,h) & \mapsto & p(gh)q(h)
\end{array}
,\qquad
\begin{array}{rccc}
(p \oda 1)\com(q)  : & G \times G & \to & \ku \\
& (g,h) & \mapsto & p(g)q(gh)
\end{array}
\]
have finite support, so belong to $K(G) \odo K(G)$; the coassociativity of $\com$ is an immediate consequence of the associativity of the multiplication on $G$. So $\com$ is a comultiplication on $K(G)$. The map $T_1 $, given by $(T_1 p)(g,h) = f(gh,h)$, is bijective and its inverse $R_1$ is given by $(R_1 p)(g,h) = p(gh^{-1},h)$. Similarly, the map $T_2$, given by $(T_2 p)(g,h) = p(g,gh)$, is bijective and its inverse $R_2$ is given by $(R_2 p)(g,h) = p(g,g^{-1}h)$. Then $(K(G),\com)$ is a commutative multiplier Hopf (\as-)algebra. In this case, it is easy to show that its counit and its antipode are given by
\[
\begin{array}{rccc}
\cou : & K(G) & \to & \ku \\
& p & \mapsto & p(e)
\end{array}
,\qquad
\begin{array}{rccc}
S  : & K(G) & \to & K(G) \\
& p & \mapsto & [S(p): g \mapsto p(g^{-1})]
\end{array},
\]
respectively. Let $\ku[G]$ be the group (\as-)algebra of $G$, i.e. the $\ku$-module generated by formal elements $\{\lambda_{g}\}_{g \in G}$ with algebra structure given by $\lambda_{g}\lambda_{h} = \lambda_{gh}$ for all $g,h \in G$; in the \as-algebra case, the involution is given by $(\lambda_{g})^* = \lambda_{g^{-1}}$. For each element $x \in \ku[G]$, there exists a unique $\ku$-valued finitely supported function $p_{x} : G \to \ku$, such that $x = \sum_{g \in G}p_{x}(g)\lambda_{g}$. We define the (\as-)homomorphism $\com': \ku[G] \to \ku[G] \odo \ku[G]$, $\lambda_{g} \mapsto \lambda_{g} \oda \lambda_{g}$, thus the pair $(\ku[G],\com')$ is a Hopf (\as-)algebra, its counit and its antipode are given by
\[
\begin{array}{rccc}
\cou' : & \ku[G] & \to & \ku \\
& \lambda_{g} & \mapsto & 1
\end{array}
, \qquad
\begin{array}{rccc}
S'  : & \ku[G] & \to & \ku[G] \\
& \lambda_{g} & \mapsto & \lambda_{g^{-1}}
\end{array}
,
\]
respectively.
\end{exam}




Let $(A,\com)$ be a multiplier Hopf algebra. For $a \in A$ and $\ome$ a $\ku$-linear functional on $A$, we will denote by $(\id \odo \ome)(\com(a))$ the element in $\M(A)$ such that
\[
[(\id \odo \ome)(\com(a))]a' := (\id \odo \ome)(\com(a)(a' \oda 1)) \quad \text{ and } \quad a'[(\id \odo \ome)(\com(a))] := (\id \odo \ome)((a' \oda 1)\com(a))
\]
for every $a' \in A$. We will use similar formulas for define the element $(\ome \odo \id)(\com(a))$.

\begin{rema}
Sometimes, for define a multiplier or made computations with a multiplier, we use brackets to indicate that we are working with a multiplier.
\end{rema} 

Let $(A,\com)$ be a multiplier Hopf algebra and $\var, \psi$ be non-zero $\ku$-linear functionals on $A$. We say that $\var$ is {\em left invariant} if $(\id \odo \var)(\com(a)) = \var(a)1_{\M(A)}$ for every $a \in A$. Similarly, we say that $\psi$ is {\em right invariant} if $(\psi \odo \id)(\com(a)) = \psi(a)1_{\M(A)}$ for every $a, b \in A$. %
%
%
%
%
If $(A,\com)$ admits a non-zero left invariant functional, then it is unique (up to a scalar) and there exists also a unique right invariant functional. It can be proof that these functionals are faithful.

\begin{rema}
Usually, non-zero invariant functionals are called integrals in the Hopf algebra theory and Haar measures in the quantum group setting. We will, throughout this work, speak about integrals.
\end{rema}

%

\subsection*{The opposite multiplier Hopf \as-algebra}

Given a multiplier Hopf \as-algebra $(A,\com)$, consider the vector space $A^{\op} := A$ with multiplication, comultiplication and involution given by
\[
a^{\op}a'^{\op} := (a'a)^{\op}, \qquad \com^{\op}(a^{\op})(a'^{\op} \odo 1^{\op}) := (a'a_{(1)})^{\op} \odo (a_{(2)})^{\op},
\]
\[
(a^{\op})^{*} := S^{-2}(a^{*})^{\op} = (S^{2}(a)^{*})^{\op},
\]
for all $a,a' \in A$, respectively. Then, $(A^{\op},\com^{\op})$ yields a multiplier Hopf \as-algebra with antipode and counit given by $S^{\ops}(a^{\op}) := S^{-1}(a)^{\op}$ and $\cou^{\ops}(a^{\op}) := \cou(a)$ for all $a \in A$, respectively. This multiplier Hopf \as-algebra will be called the {\em opposite multiplier Hopf \as-algebra of $(A,\com)$}.

\subsection*{The Sweedler type leg notation}

It is known that in the framework of Hopf algebras, the Sweedler notation for the comultiplication is a fundamental tool for give definitions and make computations. In the general setting of multiplier Hopf algebras, it is also possible use the same notation, but we need to be care because we are dealing with multiplier elements. Indeed, given a multiplier Hopf algebra $(A,\com)$, for any $a \in A$, the element $\com(a)$ belongs to $\M(A \odo A)$ instead of $A \odo A$. In the next line, we want to explain the {\em Van Daele's covering technique} in order to use the Sweedler notation in our work. 

Let $(A,\com)$ be a regular multiplier Hopf algebra. For an element $a \in A$, the Sweedler type leg notation $\com(a) =: a_{(1)} \oda a_{(2)}$ is used in the following sense: Given an element $b \in A$, because $\com(a)(b \oda 1)$ is an element of $A \odo A$, we can write
\[
a_{(1)}b \oda a_{(2)} = \com(a)(b \oda 1),
\]
where we have $a_{(1)}b, a_{(2)} \in A$. In the last equality, we say that the element $b$ is covering by the right the first leg of $\com(a)$ in order to have $a_{(1)}b \in A$. Similarly, we can also write $ba_{(1)} \oda a_{(2)} = (b \oda 1)\com(a)$, $a_{(1)} \oda a_{(2)}b = \com(a)(1 \oda b)$ and $a_{(1)} \oda ba_{(2)} = (1 \oda b)\com(a)$.

\medskip

The following propositions are useful in order to do computations with the Sweedler type leg notation.

\begin{prop}
Let $(A,\com)$ be a regular multiplier Hopf algebra. Given a finite set of elements of $A$, $\{a_{i}\}_{1\leq i\leq n}$, then there is an element $e \in A$ such that $ea_{i} = a_{i}e = a_{i}$ for every $1 \leq i \leq n$. This means that the algebra $A$ admits local units. 
\end{prop}

\begin{prop}
Let $(A,\com)$ be a regular multiplier Hopf algebra and $M$ be a left unital $A$-module through the module map $\rhd: A \odo M \to M$. We have
\begin{enumerate}[label=\textup{(\roman*)}]
\item For each element $m \in M$, there is $a \in A$ such that $a \rhd m = m$.
\item There is a unique left module map $\rhd: \M(A) \odo M \to M$ turning $M$ into a left $\M(A)$-module.
\end{enumerate}
\end{prop}

\subsection{Actions of multiplier Hopf algebras}

Let $(A,\com)$ be a regular multiplier Hopf algebra and $X$ be a non-degenerate algebra. Consider a left unital $A$-module structure on $X$, $\rhd: A \odo X \to X$. Fix $a \in A$ and $x, y \in X$, then there are $b,c \in A$ such that $b \rhd x = x$ and $c \rhd y = y$. Because $a_{(1)}b \oda a_{(2)}c = \com(a)(b \oda c) \in A \odo A$, we can define $(a_{(1)} \rhd x)(a_{(2)} \rhd y) := (a_{(1)}b \rhd x)(a_{(2)}c \rhd y) \in X$. To see that this is well defined, suppose there are two other elements $b',c' \in A$ such that $b' \rhd x = x$ and $c' \rhd y = y$. Because, $A$ admits local units, there are $b'', c'' \in A$ such that $b''b = b$, $b''b' = b'$, $c''c =c$ and $c''c' =c'$, thence
\begin{align*}
(a_{(1)}b' \rhd x)(a_{(2)}c' \rhd y) & = (a_{(1)}b''b' \rhd x)(a_{(2)}c''c' \rhd y) = (a_{(1)}b'' \rhd (b' \rhd x))(a_{(2)}c'' \rhd (c' \rhd y)) \\
& = (a_{(1)}b'' \rhd (b \rhd x))(a_{(2)}c'' \rhd (c \rhd y)) = (a_{(1)}b''b \rhd x)(a_{(2)}c''c \rhd y) \\
& = (a_{(1)}b \rhd x)(a_{(2)}c \rhd y).
\end{align*}

A {\em left action of $(A,\com)$} on $X$ is a left unital module map $\rhd: A \odo X \to X$ satisfying the condition
\[
a \rhd xy = (a_{(1)} \rhd x)(a_{(2)} \rhd y) 
\]
for all $a\in A$ and $x,y \in X$. Sometimes, we will say simply that $(X,\rhd)$ is a left $(A,\com)$-module algebra. When $A$ and $X$ are \as-algebras, we will say that $\rhd: A \odo X \to X$ is a left action if additionally, we have
\[
(a \rhd x)^{*} = S(a)^{*} \rhd x^{*}
\]
for all $a\in A$ and $x\in X$. In this case, we will say that $(X,\rhd)$ is a left $(A,\com)$-module \as-algebra. Similarly, right actions of $(A,\com)$ on non-degenerate algebras can be defined.

\begin{exam}
The {\em left and right adjoint action} of $(A,\com)$ on $A$ are given by
\[
a \brhd_{\mathrm{ad}} a' = a_{(1)}a' S(a_{(2)}) \quad \text{ and } \quad a' \blhd_{\mathrm{ad}} a = S(a_{(1)})a'a_{(2)},
\]
respectively.
\end{exam}

\begin{rema}
The Van Daele's covering technique not only use the multiplication in order to cover the legs in the Sweedler type leg notation of the comultiplication, it can also be extended to use the actions of multiplier Hopf algebras. For example, given a left action of a multiplier Hopf algebra $(A,\com)$ on a no-degenerate algebra $X$, for two elements $a \in A$ and $x \in X$, the element $(a_{(1)} \rhd x) \odo a_{(2)}$ is well defined and it belongs to $X \odo A$. Indeed, for $x \in X$ there is an element $b \in A$ such that $b \rhd x = x$, then we can write
\[
(a_{(1)} \rhd x) \odo a_{(2)} = (a_{(1)}b \rhd x) \odo a_{(2)} \in X \odo A,
\]
where we observe that the element $b$ covers by the right the first leg of $\com(a)=a_{(1)} \odo a_{(2)}$. Because $b$ arise from $x$ through the action, we will say that $x$ is covering $a_{(1)}$ through the action.
\end{rema}

It is known that we can pass from left module algebras to right module algebras over multiplier Hopf algebras using the opposite algebra construction. In the case of \as-algebras, we need to pay attention to the compatibility condition of the action with the antipode as it is show in the following remark.

\begin{rema}
If $(X,\rhd)$ is a left $(A,\com)$-module \as-algebra, then the map $\dot\lhd: X \odo A^{\op} \to X$, $x \odo a^{\op} \mapsto a \rhd x$ yields a right action of the opposite multiplier Hopf \as-algebra $(A^{\op},\com^{\op})$ on $X$. Indeed, observe that
\begin{align*}
(x \dot\lhd a^{\op})^{*} & = (a \rhd x)^{*} = S(a)^{*} \rhd x^{*} = x^{*} \dot\lhd (S^{2}(S^{-1}(a))^{*})^{\op} = x^{*} \dot\lhd (S^{-1}(a)^{\op})^{*} = x^{*} \dot\lhd S^{\ops}(a^{\op})^{*},
\end{align*}
for each $x\in X$ and $h\in A$. Similarly, we have $x^{*} \dot\lhd (a^{\op})^{*} = (x \dot\lhd S(a)^{\op})^{*}$. Moreover, the maps $S: (A^{\op},\com^{\op}) \to (A,\com^{\co})$, $a^{\op} \mapsto S(a)$ and $S: (A^{\op},\com^{\op,\co}) \to (A,\com)$, $a^{\op} \mapsto S(a)$ yield isomorphisms of multiplier Hopf \as-algebras.
\end{rema}

\subsubsection*{Smash product algebra}

Consider a left $(A,\com)$-module algebra $(X,\rhd)$. The $\ku$-module $X \# A := X \odo A$ endowed with the $\ku$-linear map $\m_{X \# A}: X \# A \odo X \# A \to X \# A$, $x \# a \odo y \# a' \mapsto  x(a_{(1)} \rhd y) \# a_{(2)}a'$ defines a non-degenerate algebra called {\em the smash product algebra associated to $\rhd$}. Sometimes, we will also use the notation $X \#_{\rhd} A$ to remark the use of the action $\rhd$. In case $\rhd$ is an action of a multiplier Hopf \as-algebra $(A,\com)$ on a non-degenerate \as-algebra $X$, the algebra $X \# A$ is endowed with a \as-operation defined by $(x \# a)^{*} = (a^{*}_{(1)} \rhd x^{*}) \# a^{*}_{(2)}$ for all $a\in A$, $x\in X$.

\subsection{Coactions of multiplier Hopf algebras}

A {\em left coaction of $(A,\com)$} on $X$ is a non-degenerate homomorphism $\Gamma: X \to \M(A \odo X)$ satisfying the conditions
\begin{equation}\label{eq:c1_coaction}\tag{C1}
(A \odo 1_{\M(X)})\Gamma(X) = A \odo X = \Gamma(X)(A \odo 1_{\M(X)})
\end{equation}
\begin{equation}\label{eq:c2_coaction}\tag{C2}
(\com \odo \id_{\M(X)})\Gamma = (\id_{\M(A)} \odo \Gamma)\Gamma, \qquad (\cou \odo \id_{\M(X)})\Gamma = \id_{\M(X)}.
\end{equation}
Sometimes, we will use the terminology $(X,\Gamma)$ is a left $(A,\com)$-comodule algebra. When $A$ and $X$ are \as-algebras, we will say that $\Gamma: X \to \M(A \odo X)$ is a left coaction if additionally $\Gamma$ is a \as-homomorphism. In this case, we will say that $(X,\Gamma)$ is a left $(A,\com)$-comodule \as-algebra. Similarly, right coactions of $(A,\com)$ on non-degenerate algebras can be defined.

\begin{rema}
It can be proof that the condition $(\cou \odo \id_{\M(X)})\circ\Gamma = \id_{\M(X)}$ for a left coaction is equivalent to the injectivity of the map $\Gamma$.
\end{rema}

\begin{rema}
Similar to the covering technique for the comultiplication of a multiplier Hopf algebras, when we work with a left coaction $\Gamma: X \to \M(A \odo X)$, we can use the Sweedler type leg notation $\Gamma(x) =: x_{[-1]} \odo x_{[0]}$ in the following sense: Given $a \in A$ and $x \in X$, we can write
\[
ax_{[-1]} \oda x_{[0]} = (a \oda 1)\Gamma(x), \qquad  x_{[-1]}a \oda x_{[0]} = \Gamma(x)(a \oda 1).
\]
In case we work with a right coaction $\Gamma: X \to \M(X \odo A)$, we can use the Sweedler type leg notation $\Gamma(x) =: x_{[0]} \odo x_{[1]}$ in a similar way. Moreover, similar to the covering by actions in the case of the Sweedler type leg notation for the comultiplication of a multiplier Hopf algebras, we can use the covering technique by actions in case of work with the Sweedler type leg notation for coactions.
\end{rema}

\begin{rema}
Given an invertible element $W \in \M(A \odo X)$ such that $(\com \odo \id_{\M(X)})(W) = W_{23}W_{13}$ and $(a \odo 1)W(1 \odo x), (1 \odo x)W^{-1}(a \odo 1) \in A \odo X$ for all $a\in A$, $x\in X$, the linear map $\alp_{W}: X \to \M(A \odo X)$, defined by $\alp_{W}(x) = \ad(W)(1_{\M(A)} \odo x) = W(1_{\M(A)} \odo x)W^{-1}$ for all $x \in X$, yields a left coaction of $(A,\com)$ on $X$. 
In the involutive case, in order to obtain a coaction, the element $W$ need to be an unitary element in $\M(A \odo X)$, i.e. we need to have $W^{*} = W^{-1}$.
\end{rema}

At some point, we will deal with coactions of the opposite and the co-opposite multiplier Hopf algebra of a given multiplier Hopf algebra. The following propositions will be useful in that case.

\begin{prop}\label{prop:op_co_coaction}
Let $(A,\com)$ be a multiplier Hopf algebra, $X$ be a non-degenerate algebra and $\Gamma: X \to \M(X \odo A)$ be a map. Consider the map $\Gamma^{\op}: X^{\op} \to \M(A^{\op} \odo X^{\op})$, defined by $\Gamma^{\op}(x^{\op}) = ({}^{\op} \odo {}^{\op})(\Gamma(x))$. Consider the map $\Gamma^{\co}: X^{\op} \to \M(A^{\co} \odo X^{\op})$, defined by $\Gamma^{\co}(a^{\op}) = (S \odo {}^{\op})(\Gamma(x))$.

\begin{enumerate}[label=\textup{(\roman*)}]
\item $(X,\Gamma)$ is a left $(A,\com)$-comodule algebra;
\item $(X^{\op},\Gamma^{\op})$ is a left $(A^{\op},\com^{\op})$-comodule algebra;
\item $(X^{\op},\Gamma^{\co})$ is a left $(A,\com^{\co})$-comodule algebra. 
\end{enumerate}
Moreover, using the isomorphism of multiplier Hopf algebras $S\circ{}^{\op}: (A^{\op},\com^{\op}) \to (A,\com^{\co})$, we have
\[
(S\circ{}^{\op} \odo \id_{\M(X)})\circ\Gamma^{\op} = \Gamma^{\co}.
\]
\end{prop}
\pr
Given $\Gamma: X \to \M(A \odo X)$ a left coaction of the multiplier Hopf \as-algebra $(\pol(\G),\com_{\G})$ on the non-degenerate algebra $X$. By definition, $\Gamma$ is an injective map and we have the conditions
\[
(\com \odo \id)\circ\Gamma = (\id \odo \Gamma)\circ\Gamma, \qquad (A \odo 1_{\M(X)})\Gamma(X) = A \odo X = \Gamma(X)(A \odo 1_{\M(X)}).
\]
\begin{enumerate}[label=\textup{(\roman*)}]
\item For any non-degenerate idempotent algebra $Y$, we have $\M(Y^{\op}) = \M(Y)^{\op}$. With that in mind, observe that the map $\Gamma^{\op}$ satisfy $\Gamma^{\op} = ({}^{\op} \odo {}^{\op})\circ\Gamma\circ {}^{\op}$. It follows then that $\Gamma^{\op}$ is injective and
\[
\Gamma^{\op}(X^{\op})(A^{\op} \odo 1^{\op}_{\M(X)}) = ((A \odo 1_{\M(X)})\Gamma(X))^{(\op \odo \op)} = A^{\op} \odo X^{\op}.
\]
Similarly, we can obtain $(A^{\op} \odo 1^{\op}_{\M(X)})\Gamma^{\op}(X^{\op}) = A^{\op} \odo X^{\op}$. Moreover, by direct computation we have
\begin{align*}
(\id \odo \Gamma^{\op})\circ\Gamma^{\op}\circ{}^{\op} & = ({}^{\op} \odo {}^{\op} \odo {}^{\op})\circ(\id \odo \Gamma)\circ\Gamma \\
& = ({}^{\op} \odo {}^{\op} \odo {}^{\op})\circ(\com \odo \Gamma)\circ\Gamma \\
& = (\com^{\op} \odo \id)\circ({}^{\op} \odo {}^{\op})\circ\Gamma \\
& = (\com^{\op} \odo \id)\circ\Gamma^{\op}\circ{}^{\op}.
\end{align*}
The last statements implies that $\Gamma^{\op}$ yields a left coaction of $(A^{\op},\com^{\op})$ on the algebra $X^{\op}$.

\item Because the map $\Gamma^{\co}$ satisfy $\Gamma^{\co} = (S \odo {}^{\op})\circ\Gamma\circ{}^{\op}$, then the last equality implies that $\Gamma^{\co}$ is injective and
\[
\qquad\Gamma^{\co}(X^{\op})(A \odo 1^{\op}_{\M(X)}) = (S \odo {}^{\op})((S^{-1}(A) \odo 1_{\M(X)})\Gamma(X)) = (S \odo {}^{\op})(A \odo X) = A \odo X^{\op}.
\]
Similarly, we can obtain $(A \odo 1^{\op}_{\M(X)})\Gamma^{\co}(X^{\op}) = (A \odo X^{\op})$. Moreover, we have
\begin{align*}
(\id \odo \Gamma^{\co})\circ\Gamma^{\co}\circ{}^{\op} & = (S \odo S \odo {}^{\op})\circ(\id \odo \Gamma)\circ\Gamma \\
& = (S \odo S \odo {}^{\op})\circ(\com \odo \id)\circ\Gamma \\
& = (\com^{\co} \odo \id)\circ(S \odo {}^{\op})\circ\Gamma \\
& = (\com^{\co} \odo \id)\circ\Gamma^{\co}\circ{}^{\op}.
\end{align*}
The last statements implies that $\Gamma^{\co}$ yields a left coaction of $(A,\com^{\co})$ on the algebra $X^{\op}$.
\end{enumerate}

Finally, note that the map $S\circ{}^{\op}: A^{\op} \to A$, $a^{\op} \mapsto S(a)$ is an isomorphism between the multiplier Hopf algebras $(A^{\op},\com^{\op})$ and $(A,\com^{\co})$. Hence, it follows that
\begin{align*}
(S\circ{}^{\op} \odo \id)\circ\Gamma^{\op}\circ{}^{\op} & = (S \odo {}^{\op})\circ\Gamma = \Gamma^{\co}\circ{}^{\op}.
\end{align*}
\fin

\begin{defi}
Given a left coaction $\Gamma$ of a multiplier Hopf algebra on a non-degenerate algebra. The coactions $\Gamma^{\op}$, $\Gamma^{\co}$ defined in the proposition above are called {\em the opposite coaction and the co-opposite coaction of $\Gamma$}, respectively.
\end{defi}

In case, we are working with involutive algebras, in order to define that the opposite and co-opposite coaction are coactions of multiplier Hopf \as-algebras, we need to have an non-trivial involution on the opposite algebra compatible with the coaction as it is prove in the following proposition.

\begin{prop}\label{prop:op_co_coaction_involutive}
Let $(A,\com)$ be a multiplier Hopf \as-algebra, $X$ be a non-degenerate \as-algebra endowed with an algebra automorphism $\gamma$ such that $\gamma\circ * \circ \gamma \circ *= \id_{\M(X)}$ and $\Gamma: X \to \M(A \odo X)$ be a map. If the condition
\[
\Gamma\circ\gamma = (S^{-2} \odo \gamma)\circ\Gamma
\]
holds, then the following conditions are equivalents
\begin{enumerate}[label=\textup{(\roman*)}]
\item $(X,\Gamma)$ is a left $(A,\com)$-comodule \as-algebra
\item $(X^{\op}_{\gamma},\Gamma^{\op})$ is a left $(A^{\op},\com^{\op})$-comodule \as-algebra.
\item $(X^{\op}_{\gamma},\Gamma^{\co})$ is a left $(A,\com^{\co})$-comodule \as-algebra
\end{enumerate}
Moreover, using the isomorphism of multiplier Hopf \as-algebras $S\circ{}^{\op}: (A^{\op},\com^{\op}) \to (A,\com^{\co})$, we have
\[
(S\circ{}^{\op} \odo \id_{\M(X)})\circ\Gamma^{\op} = \Gamma^{\co}.
\]
\end{prop}
\pr
By Proposition~\ref{prop:op_co_coaction}, we only need to prove the compatibility with the involution for the opposite and the co-opposite coaction. First note that by definition the involution in $X^{\op}_{\gamma}$ is given by $(x^{\op})^{*} = \gamma(x^{*})^{\op} = \hat\gamma(x)^{* \op}$, where $\hat\gamma$ denotes the inverse isomorphism of $\gamma$. In other words,
\[
* \circ {}^{\op} = {}^{\op}\circ\gamma\circ * = {}^{\op}\circ * \circ \hat\gamma.
\]
as maps from $X$ to $X^{\op}_{\gamma}$. Assume that the condition
\[
\Gamma\circ\gamma = (S^{-2} \odo \gamma)\circ\Gamma
\]
holds. By direct computation:
\begin{enumerate}[label=\textup{(\roman*)}]
\item On one hand we have
\begin{align*}
\Gamma^{\op}\circ * \circ{}^{\op} & = \Gamma^{\op}\circ{}^{\op}\circ \gamma \circ *  = ({}^{\op} \odo {}^{\op})\circ\Gamma\circ\gamma\circ * = ({}^{\op} \odo {}^{\op})\circ(S^{-2} \odo \gamma)\circ\Gamma\circ *
\end{align*}
and on the other
\begin{align*}
(* \odo *)\circ\Gamma^{\op}\circ{}^{\op} & = (* \odo *)\circ({}^{\op} \odo {}^{\op})\circ\Gamma = ({}^{\op} \odo {}^{\op})\circ(S^{-2} \odo \gamma)\circ(* \odo *)\circ\Gamma.
\end{align*}

\item On one hand we have
\begin{align*}
\Gamma^{\co}\circ *\circ{}^{\op} & = \Gamma^{\co}\circ{}^{\op}\circ \gamma \circ *  = (S \odo {}^{\op})\circ\Gamma\circ\gamma\circ * = (S \odo {}^{\op})\circ(S^{-2} \odo \gamma)\circ\Gamma\circ * 
\end{align*}
and on the other
\begin{align*}
(* \odo *)\circ\Gamma^{\co}\circ{}^{\op} & = (* \odo *)\circ(S \odo {}^{\op})\circ\Gamma = (S \odo {}^{\op})\circ(S^{-2} \odo \gamma)\circ(* \circ *)\circ\Gamma.
\end{align*}
\end{enumerate}
The proposition follows directly from the last equalities.
\fin

\begin{nota}
Given a left action $\Gamma: X \to \M(A^{\op} \odo X)$ of a multiplier Hopf algebra $(A^{\op},\com^{\op})$ on a non-degenerate algebra $X$, sometimes it will be used the underlying module map of $\Gamma$ without the consideration of the opposite multiplication structure on the multiplier Hopf \as-algebra $(A^{\op},\com^{\op})$, i.e. we will use the map $\dot{\Gamma}: X \to \M(A \odo X)$, defined by $\dot{\Gamma} = ({}^{\op} \odo \id)\circ\Gamma$.
\end{nota}

\subsection{Pairings of multiplier Hopf algebras} Let $(A,\com_{A})$, $(B,\com_{B})$ be two regular multiplier Hopf algebras and $\p: A \times B \to \ku$ be a non-degenerate $\ku$-bilinear map. For any $a\in A$ and $b\in B$, consider the $\ku$-linear maps ${}_{a}\p := \p(a,\,\cdot\,) \in B^{\vee}$ and $\p_{b} := \p(\,\cdot\,,b) \in A^{\vee}$. With those linear maps, given $a\in A$, $b\in B$, define the multipliers $({}_{a}\p \odo \id)(\com_{B}(b))$, $(\id \odo {}_{a}\p)(\com_{B}(b)) \in \M(B)$ and $(\p_{b} \odo \id)(\com_{A}(a))$, $(\id \odo \p_{b})(\com_{A}(a)) \in \M(A)$ by
\small
\[
[(\id \odo {}_{a}\p)(\com_{B}(b))]b' = (\id \odo {}_{a}\p)(\com_{B}(b)(b' \odo 1)), \quad b'[(\id \odo {}_{a}\p)(\com_{B}(b))] = (\id \odo {}_{a}\p)((b' \odo 1)\com_{B}(b))
\]
\[
[({}_{a}\p \odo \id)(\com_{B}(b))]b' = ({}_{a}\p \odo \id)(\com_{B}(b)(1 \odo b')), \quad b'[({}_{a}\p \odo \id)(\com_{B}(b))] = (\id \odo {}_{a}\p)((1 \odo b')\com_{B}(b))
\]
\[
[(\id \odo \p_{b})(\com_{A}(a))]a' = (\id \odo \p_{b})(\com_{A}(a)(a' \odo 1)), \quad a'[(\id \odo \p_{b})(\com_{A}(a))] = (\id \odo \p_{b})((a' \odo 1)\com_{A}(a))
\]
\[
[(\p_{b} \odo \id)(\com_{A}(a))]a' = (\p_{b} \odo \id)(\com_{A}(a)(1 \odo a')), \quad a'[(\p_{b} \odo \id)(\com_{A}(a))] = (\p_{b} \odo \id)((1 \odo a')\com_{B}(b))
\]
\normalsize
for all $a'\in A$, $b'\in B$. The bilinear map $\p$ is called {\em a pre-pairing between $(A,\com_{A})$ and $(B,\com_{B})$} if for all $a,a' \in A$ and $b,b' \in B$, the following conditions hold
\begin{itemize}
\item[(i)] $({}_{a}\p \odo \id)(\com_{B}(b)), (\id \odo {}_{a}\p)(\com_{B}(b)) \in B$ and $(\p_{b} \odo \id)(\com_{A}(a)), (\id \odo \p_{b})(\com_{A}(a)) \in  A$,
\item[(ii)] ${}_{a}\p(\id \odo {}_{a'}\p)(\com_{B}(b)) = {}_{a'}\p({}_{a}\p \odo \id)(\com_{B}(b)) = {}_{aa'}\p(b)$,
\item[(iii)] $\p_{b}(\id \odo \p_{b'})(\com_{A}(a)) = \p_{b'}(\p_{b} \odo \id)(\com_{A}(a)) = \p_{bb'}(a)$.
\end{itemize}

Follows from the definition that associated with a pre-pairing $\p$, there are four canonical module maps defined by
\[
\begin{array}{lccc}
\brhd: & A \odo B & \to & B \\
& a \odo b & \mapsto & b_{(1)}\p(a,b_{(2)})
\end{array},
\qquad
\begin{array}{lccc}
\blhd: & B \odo A & \to & B \\
& b \odo a & \mapsto & \p(a,b_{(1)})b_{(2)}
\end{array},
\]
\[
\begin{array}{lccc}
\brhd: & B \odo A & \to & A \\
& b \odo a & \mapsto & a_{(1)}\p(a_{(2)},b)
\end{array},
\qquad
\begin{array}{lccc}
\blhd: & A \odo B & \to & A \\
& a \odo b & \mapsto & \p(a_{(1)},b)a_{(2)}
\end{array}.
\]
If those actions are unital, we say that $\p$ is a {\em pairing between $(A,\com_{A})$ and $(B,\com_{B})$}. In fact, it was shown in \cite{DVD01} that it is enough only have a unital action among those actions. The above unital actions are called {\em the regular actions associated to the pairing $\p$}. Given a pairing $\p$, it hold the following relations
\[
\p(aa',b) = \p(a, a' \brhd b) = \p(a',b \blhd a), \quad \p(a,bb') = \p(a \blhd b,b') = \p(b' \brhd a,b),
\]
\[
\p(a,S_{B}(b)) = \p(S_{A}(a),b)
\]
for all $a,a' \in A$ and $b,b' \in B$. In case $A$ and $B$ are multiplier Hopf \as-algebras, a pairing between between $(A,\com_{A})$ and $(B,\com_{B})$ is a pairing $\p$ satisfying additionally the conditions
\[
\p(a^*,b) = \overline{\p(a,S_{B}(b)^*)}, \quad\quad \p(a,b^*) = \overline{\p(S_{A}(a)^*,b)}
\]
for all $a \in A$ and $b \in B$.

\begin{nota}
If $\p$ is a pairing between $(A,\com_{A})$ and $(B,\com_{B})$, we have
\begin{enumerate}
\item the map $\bar{\p}: B \times A \to \ku$, $(b,a) \mapsto \p(a,b)$ yields a pairing between $(B,\com_{B})$ and $(A,\com_{A})$ called {\em the flip pairing of $\p$},

\item the map $\hat{\p}: B \times A^{\op} \to \ku$, $(b,a^{\op}) \mapsto \p(a,b)$ yields a pairing between $(B,\com^{\co}_{B})$ and $(A^{\op},\com^{\op}_{A})$ called {\em the flip co-op pairing of $\p$},

\item the map $\tilde{\p}: B^{\op} \times A^{\op} \to \ku$, $(b^{\op},a^{\op}) \mapsto \p(a,b)$ yields a pairing between $(B^{\op},\com^{\co,\op}_{B})$ and $(A^{\op},\com^{\co,\op}_{A})$ called {\em the flip co-opposite opposite pairing of $\p$},
\end{enumerate}
\end{nota}

\subsection{The canonical multiplier of a pairing}

Let $\p: A \times B \to \ku$ be a pairing between two regular multiplier Hopf algebras. Consider the non-degenerate bilinear form $\p^{2} : A \odo B \times B \odo A \to \ku$ defined by $\p^{2}(a \odo b,b' \odo a') = \p(a,b')\bar{\p}(b,a')$. We will say that $\p$ admits {\em a canonical multiplier $U$} if $U$ is an invertible element in $\M(A \odo B)$ satisfying the condition $\p^{2}(U,b \odo a) = \p(a,b)$ for all $a \in A$ and $b \in B$. By non-degeneracy of $\p^{2}$, it follows the unicity of $U$. Sometimes, we will use the notation $U(\p)$ to emphasize the relation with the pairing $\p$. Using the pairing $\p$, one can show the equalities
\[
(\com_{A} \odo \id)(U) = U_{13}U_{23} \quad \text{ and } \quad (\id \odo \com_{B})(U) = U_{12}U_{13}
\]
in $\M(A \odo A \odo B)$ and $\M(A \odo B \odo B)$, respectively. Here, we are using the {\em leg notations}
\[
U_{13}:= (\id \odo \Sigma)(U \odo 1_{\M(A)}), \quad  U_{23} := 1_{\M(A)} \odo U
\]
and
\[
U_{12} := U \odo 1_{\M(B)}, \quad U_{13}:= (\id \odo \Sigma)(U \odo 1_{\M(B)}),
\]
respectively. Moreover, we also have
\[
U = (S^{-1}_{A} \odo S_{B})(U) = (S_{A} \odo S^{-1}_{B})(U) \quad {and} \quad U^{-1} = (\id_{A} \odo S_{B})(U) = (S_{A} \odo \id_{B})(U).
\]

A pairing $\p: A \times B \to \ku$ will be called {\em admissible pairing} if it admits a canonical multiplier $U$ such that for each $a\in A$, $b \in B$, the elements of the form $(1 \odo b)U(a \odo 1), (a \odo 1)U(1 \odo b)$ belong to $A \odo B$. 

\medskip

The following proposition gives important results in computations related to the canonical multiplier of an admissible pairing. Such results will be useful in Section~\ref{sec:oc_yd}.

\begin{prop}
Let $\p: A \times B \to \ku$ be a pairing of regular multiplier Hopf algebras which admits a canonical multiplier $U$. Consider the isomorphism $T_{U} = (\id \odo S^{-1}_{B})\ad(U)(\id \odo S_{B}) \in \mathrm{End}(A \odo B)$. Then, for any $a,a' \in A$ and $b \in B$, it holds
\begin{enumerate}[label=\textup{(p\arabic*)}]
\item\label{eq:ad_p_1} $(\id \odo {}_{a'}\p)(U(a \odo b)U^{-1}) = (b \brhd a'_{(1)}) a S_{A}(a'_{(2)}) = a'_{(1)} a S_{A}(a'_{(2)} \blhd b)$;
\item\label{eq:ad_p_2} $(\id \odo {}_{a'}\p)(U^{-1}(a \odo b)U) = S_{A}(a'_{(1)})a(a'_{(2)} \blhd b) = S_{A}(b \brhd a'_{(1)})aa'_{(2)}$;
\item\label{eq:t_p_1} $(\id \odo {}_{a'}\p)(T_{U}(a \odo b)) = S^{-1}_{A}(a'_{(2)} \blhd b)aa'_{(1)} = S^{-1}_{A}(a'_{(2)})a(b \brhd a'_{(1)})$;
\item\label{eq:t_p_2} $(\id \odo {}_{a'}\p)(T^{-1}_{U}(a \odo b)) = a'_{(2)} a S^{-1}_{A}(b \brhd a'_{(1)}) = (a'_{(2)} \blhd b) a S^{-1}_{A}(a'_{(1)})$.
\end{enumerate}

\end{prop}
\pr
Fix $a,a' \in A$ and $b \in B$.
\begin{enumerate}[label=\textup{(p\arabic*)}]
\item Because
\begin{align*}
\p((\id \odo {}_{a'}\p)(U(a \odo b)U^{-1}),b') & = \p^{2}(U(a \odo b)U^{-1},b' \odo a') \\
& = \p^{2}(U,b'_{(1)} \odo a'_{(1)})\p^{2}(a \odo b,b'_{(2)} \odo a'_{(2)})\p^{2}(U^{-1},b'_{(3)} \odo a'_{(3)}) \\
& = \p^{2}(U,b'_{(1)} \odo a'_{(1)})\p(a,b'_{(2)})\p(a'_{(2)},b)\p^{2}((\id \odo S_{B})(U),b'_{(3)} \odo a'_{(3)}) \\
& = \p^{2}(U,b'_{(1)} \odo a'_{(1)})\p(a,b'_{(2)})\p(a'_{(2)},b)\p^{2}(U,b'_{(3)} \odo S_{A}(a'_{(3)})) \\
& = \p(a'_{(1)},b'_{(1)})\p(a,b'_{(2)})\p(a'_{(2)},b)\p(S_{A}(a'_{(3)}),b'_{(3)}) \\
& = \p(a'_{(1)}aS(a'_{(3)}),b')\p(a'_{(2)},b') \\
& = \p(a'_{(1)}\p(a'_{(2)},b)aS_{A}(a'_{(3)}),b') \\
& = \p((b \brhd a'_{(1)}) a S_{A}(a'_{(2)}),b')
\end{align*}
for all $b' \in B$, we have the result.

\item Similar to the above item.

\item  We have
\begin{align*}
\p(S_{B}(b) \brhd S^{-1}_{A}(a), b') & = \p(S^{-1}_{A}(a),b'S_{B}(b)) \\
& = \p(\com_{A}(S^{-1}_{A}(a)), b' \odo S_{B}(b)) \\
& = \p(\com_{A}(S^{-1}_{A}(a)), (S_{B} \odo S_{B})(S^{-1}_{B}(b') \odo b)) \\
& = \p((S_{A} \odo S_{A})(\com_{A}(S^{-1}_{A}(a)), S^{-1}_{B}(b') \odo b) \\
& = \p(\com_{A}(a), b \odo S^{-1}_{B}(b')) \\
& = \p(a,bS^{-1}_{B}(b')) \\
& = \p(a \blhd b, S^{-1}_{B}(b')) \\
& = \p(S^{-1}_{A}(a \blhd b),b')
\end{align*}
for all $a\in A$, $b,b' \in B$. Then, because
\begin{align*}
\p((\id \odo {}_{a'}\p)(T_{U}(a \odo b)),b') & = \p^{2}((\id \odo S^{-1}_{B})\ad(U)(\id \odo S_{B})(a \odo b),b' \odo a') \\
& = \p^{2}(U(a \odo S_{B}(b))U^{-1},b' \odo S^{-1}_{A}(a')) \\
& = \p((\id \odo {}_{S^{-1}_{A}(a')}\p)(U(a \odo S_{B}(b))U^{-1}),b') \\
& = \p((S_{B}(b) \brhd S^{-1}_{A}(a')_{(1)})aS_{A}(S^{-1}_{A}(a')_{(2)}),b') \\
& = \p((S_{B}(b) \brhd S^{-1}_{A}(a'_{(2)}))aa'_{(1)},b') \\
& = \p(S^{-1}_{A}(a'_{(2)} \blhd b)aa'_{(1)},b')
\end{align*}
for all $b' \in B$, we have the result.

\item Note that $T^{-1}_{U} = T_{U^{-1}} = (\id \odo S^{-1}_{B})\ad(U^{-1})(\id \odo S_{B}) \in \textrm{End}(A \odo B)$. The proof follows in a similar way that the item above.
\end{enumerate}
\fin

The following remark shows that there are other admissible pairings associated with a given admissible pairing.

\begin{rema}
If $\p$ is an admissible pairing between $(A,\com_{A})$ and $(B,\com_{B})$ with canonical multiplier $U \in \M(A \odo B)$, thus we have
\begin{enumerate}
\item The flip pairing $\bar{\p}: B \times A \to \ku$, $(b,a) \mapsto \p(a,b)$ between $(B,\com_{B})$ and $(A,\com_{A})$ is also an admissible pairing with canonical multiplier given by $\bar{U} := \Sigma(U)$.

\item The flip co-op pairing $\hat{\p}: B \times A^{\op} \to \ku$, $(b,a^{\op}) \mapsto \p(a,b)$ between $(B,\com^{\co}_{B})$ and $(A^{\op},\com^{\op}_{A})$ is also an admissible pairing with canonical multiplier given by $\hat{U} := (\id \odo {}^{\op})(\Sigma(U))$.

\item The flip co-opposite opposite pairing $\tilde{\p}: B^{\op} \times A^{\op} \to \ku$, $(b^{\op},a^{\op}) \mapsto \p(a,b)$ between $(B^{\op},\com^{\co,\op}_{B})$ and $(A^{\op},\com^{\co,\op}_{A})$ is also an admissible pairing with canonical multiplier given by $\tilde{U} := ({}^{\op} \odo {}^{\op})(\Sigma(U))$.
\end{enumerate}
\end{rema}

\subsection{The Heisenberg algebra associated with a pairing} Let $\p: A \times B \to \ku$ be a pairing between two regular multiplier Hopf algebras. The {\em Heisenberg algebra associated with $\p$} is defined as the smash product algebra $A\#_{\brhd}B$ arising from the canonical action $\brhd: B \odo A \to A$, $b \brhd a = a_{(1)}\p(a_{(2)},b)$. This algebra, denoted by $\He(\p)$, is non-degenerate and its multiplication is defined by
\[
(a \# b)(a' \# b') = a(b_{(1)} \brhd a') \# b_{(2)}b'
\]
for all $a,a'\in A$, $b,b'\in B$. In case, $\p$ is a pairing of multiplier Hopf \as-algebras $(A,\com_{A})$ and $(B,\com_{B})$, the Heisenberg algebra associated with $\p$ admits an \as-operation given by
\[
(a \# b)^{*} = (b^{*}_{(1)} \brhd a^{*}) \# b^{*}_{(2)}
\]
for all $a\in A$, $b\in B$. In this case, we say that $\He(\p)$ is the {\em Heisenberg \as-algebra associated with $\p$}.

\begin{rema}
Sometimes, when we do computations with elements of $\He(\p) = A \# B$, we will use the following notations: $ab := a \# b$ and $ba := (b_{(1)} \brhd a) \# b_{(2)}$ for $a\in A$ and $b\in B$. We have the relations
\[
ba = (b_{(1)} \brhd a)b_{(2)} = a_{(1)}(b \blhd a_{(2)}), \quad ab = b_{(2)}(S^{-1}_{B}(b_{(1)}) \brhd a) = (b \blhd S^{-1}_{A}(a_{(2)}))a_{(1)}
\]
for all $a\in A$ and $b\in B$. Those relations are called the {\em Heisenberg's relations related to the pairing $\p$}.
\end{rema}

Consider the flip pairing of $\p$, i.e. the pairing $\bar{\p}$, the next proposition shows the relation between the Heisenberg algebras $\He(\p)$ and $\He(\bar{\p})$.

\begin{prop}\label{prop:iso_heisenberg}
The map $\mathcal{L}: \He(\bar{\p}) \to \He(\p)$, $b\,\#\,a \mapsto S^{-1}_{A}(a)\;\#\,S_{B}(b)$ defines an anti-isomorphism. Moreover, if we are dealing with multiplier Hopf \as-algebras, the map $\mathcal{L}$ is an anti-\as-isomorphism if $S_{A}^{2} = \id_{A}$ and $S_{B}^{2} = \id_{B}$. 
\end{prop}
\pr
Follows from the bijectivity of the antipode of a regular multiplier Hopf algebra that the inverse map of $\mathcal{L}$ is given by $\mathcal{L}^{-1}: \He(\p) \to \He(\bar{\p})$, $a\,\#\,b \mapsto S^{-1}_{B}(b)\;\#\,S_{A}(a)$. Then, it is enough prove that $\mathcal{L}$ is an anti-homomorphism. We have
\begin{align*}
\mathcal{L}((b\,\#\,a)(b'\,\#\,a')) & = \mathcal{L}(bb'_{(1)}\,\#\,(a \blhd b'_{(2)})a') = S^{-1}_{A}((a \blhd b'_{(2)})a')\,\#\,S_{B}(bb'_{(1)}) \\
& = S^{-1}_{A}(a')S^{-1}_{A}(a \blhd b'_{(2)})\,\#\,S_{B}(b'_{(1)})S_{B}(b) \\
& = S^{-1}_{A}(a')(S_{B}(b'_{(2)}) \brhd S^{-1}_{A}(a))\,\#\,S_{B}(b'_{(1)})S_{B}(b) \\
& = S^{-1}_{A}(a')(S_{B}(b')_{(1)} \brhd S^{-1}_{A}(a))\,\#\,S_{B}(b')_{(2)}S_{B}(b) \\
& = (S^{-1}_{A}(a')\;\#\,S_{B}(b'))(S^{-1}_{A}(a)\;\#\,S_{B}(b)) \\
& = \mathcal{L}(b'\,\#\,a')\mathcal{L}(b\,\#\,a)
\end{align*}
for all $a,a'\in A$ and $b, b'\in B$. Then $\mathcal{L}$ is an anti-isomorphism. In the involutive case, if $S_{A}^{2} = \id_{A}$ and $S_{B}^{2} = \id_{B}$, we have
\begin{align*}
\mathcal{L}((b\,\#\,a)^{*}) & = \mathcal{L}(b^{*}_{(1)}\,\#\,(a^{*} \blhd b^{*}_{(2)})) = S^{-1}_{A}((a^{*} \blhd b^{*}_{(2)}))\,\#\,S_{B}(b^{*}_{(1)}) \\
& = (S_{B}(b^{*}_{(2)}) \brhd S^{-1}_{A}(a^{*}))\,\#\,S_{B}(b^{*}_{(1)}) = (S^{-1}_{A}(a^{*})^{*}\,\#\,S_{B}(b^{*})^{*})^{*} \\
& = (S_{A}(a)\,\#\,S^{-1}_{B}(b))^{*} = (S^{-1}_{A}(a)\,\#\,S_{B}(b))^{*} = (\mathcal{L}(\ome\,\#\,a))^{*}.
\end{align*}
for all $a\in A$ and $b\in B$. Then $\mathcal{L}$ is an anti-\as-isomorphism.
\fin

The next proposition shows that there is another anti-isomorphism between $\He(\p)$ and $\He(\bar{\p})$. This anti-isomorphism is the so-called Lu's isomorphism in the case of Hopf algebras.

\begin{prop}
The map $\mathcal{L'}: \He(\bar{\p}) \to \He(\p)$, $b\,\#\,a \mapsto S_{A}(a)\;\#\,S^{-1}_{B}(b)$ defines an anti-isomorphism. Moreover, if we are dealing with multiplier Hopf \as-algebras, the map $\mathcal{L'}$ is an anti-\as-isomorphism if $S_{A}^{2} = \id_{A}$ and $S_{B}^{2} = \id_{B}$. In this case $\mathcal{L}' = \mathcal{L}$.
\end{prop}
\pr
Similar to the proposition above.

\subsection{Duality between actions and coactions associated with a pairing}\label{subsec:duality_act_coact} Let $\p: A \times B \to \ku$ be an admissible pairing between two regular multiplier Hopf algebras. It was shown in \cite{VDZ99}, in the case of canonical pairings, and in \cite{D05}, in a more general case, that there is an equivalence of categories between the category of coactions of $(A,\com_{A})$ and the category of actions of $(B,\com_{B})$.

Given a left coaction $\Gamma$ of $(A,\com_{A})$ on $X$, the map $\lhd_{\Gamma,\p}: X \odo B \to X$, $x \lhd_{\Gamma,\p} b = (\p_{b} \odo \id)\Gamma(x)$ defines a right action of $(B,\com_{B})$ on $X$. In the other hand, given a right action $\lhd: X \odo B \to B$, the map $\Gamma_{\lhd,\p}: X \to \M(A \odo X)$, $x \mapsto (1_{\M(A)} \odo x) \bar{\lhd} (a \odo 1)U$ defines a left coaction of $(A,\com_{A})$ on $X$. Here, we are using the extension to multiplier algebras of the right action $\bar{\lhd}: A \odo X \odo A \odo B \mapsto A \odo X$, $(a' \odo x) \bar{\lhd} (a \odo b) = a'a \odo (x \lhd b)$. Thus, the following functor
\[
\begin{array}{rccc}
\mathfrak{D} : & {}^{A}\alg & \leadsto & \alg_{B} \\
& (X,\Gamma) & \mapsto & (X,\lhd_{\Gamma,\p}) \\
& f: (X,\Gamma) \to (X',\Gamma') & \mapsto & f: (X,\lhd_{\Gamma,\p}) \to (X',\lhd_{\Gamma',\p})
\end{array}
\]
defines an equivalence of categories with inverse functor given by
\[
\begin{array}{rccc}
\mathfrak{D}^{-1} : & \alg_{B} & \leadsto & {}^{A}\alg \\
& (X,\lhd) & \mapsto & (X,\Gamma_{\lhd,\p}) \\
& f: (X,\lhd) \to (X',\lhd') & \mapsto & f: (X,\Gamma_{\lhd,\p}) \to (X,\Gamma_{\lhd',\p})
\end{array}.
\]

\subsection{Crossed product algebra associated to a coaction}

Let $\p: A \times B \to \ku$ be an admissible pairing between two regular multiplier Hopf algebras. Consider a left $(A,\com_{A})$-comodule algebra $(X,\Gamma)$. The {\em crossed product algebra associated to $\Gamma$ and relative to $\p$} is the $\ku$-module
\[
(A,\com) \ltimes_{\Gamma,\p} X : = \text{span}_{\ku}\{(b \odo 1_{\M(X)})\Gamma(x) : b\in B, x \in X\} \subseteq \M(\He(\bar{\p}) \odo X)
\]
with algebra structure given by
\[
(b \odo 1)\Gamma(x)\cdot(b' \odo 1)\Gamma(x') := (b(x_{[-1]} \brhd b') \odo 1)\Gamma(x_{[0]}x'),
\]
for all $b,b' \in B$ and $x,x' \in X$. Here, we are using the Sweedler type leg notation $\Gamma(x) = x_{[-1]} \odo x_{[0]}$ for each $x \in X$. When, $A$, $B$ and $X$ are \as-algebras, the algebra $(A,\com) \ltimes_{\Gamma,\p} X$ is endowed with an involution defined by
\[
((b \odo 1)\Gamma(x))^{*} := \Gamma(x^{*})(b^{*} \odo 1) = ((x^{*}_{[-1]} \brhd b^{*}) \odo 1)\Gamma(x^{*}_{[0]}),
\]
for each $b \in B$, $x \in X$. In the \as-algebra $\M(\He(\bar{\p}) \odo X)$, the following relations also hold
\[
\Gamma(x)(b \odo 1) = ((x_{[-1]} \brhd b) \odo 1)\Gamma(x_{[0]}),\quad (b \odo 1)\Gamma(x) = \Gamma(x_{[0]})((S^{-1}(x_{[-1]}) \brhd b) \odo 1),
\]
and
\[
(b^{*} \odo 1)\Gamma(x^{*}) = (\Gamma(x)(b \odo 1))^{*}
\]
for all $b \in B$, $x \in X$.


\begin{prop}
Let $(X,\Gamma)$ be a left $(A,\com_{A})$-comodule algebra. Consider the right action of $(B,\com_{B})$ on $X$ arising from $\Gamma$ related to $\p$, i.e. the module map $\lhd_{\Gamma,\p}: X \odo B \to X$, $x \lhd_{\Gamma,\p} b = (\p_{b} \odo \id)(\Gamma(x))$. 
Then, the $\ku$-linear map
\[
\begin{array}{lccc}
\mathcal{I}_{\Gamma} : & (A,\com_{A}) \ltimes_{\Gamma,\p} X & \to & B\,\#_{\lhd_{\Gamma,\p}}\,X \\
& (b \odo 1)\Gamma(x) & \mapsto &b\,\#\,x
\end{array}
\]
yields an isomorphism of algebras. Moreover, when we are dealing with involutive $*$-algebras, the isomorphism $\mathcal{I}_{\Gamma}$ preserve the involution, i.e. it is an isomorphism of \as-algebras.
\end{prop}
\pr
By direct computation, we have
\begin{align*}
\mathcal{I}_{\Gamma}((b \odo 1)\Gamma(x)\cdot(b' \odo 1)\Gamma(x')) & = \mathcal{I}_{\Gamma}((b(x_{[-1]} \brhd b') \odo 1)\Gamma(x_{[0]}x')) \\
& = b(x_{[-1]} \brhd b') \# x_{[0]}x' = bb'_{(1)} \# \p(x_{[-1]},b'_{(2)})x_{[0]}x' \\
& = bb'_{(1)} \# (x \lhd_{\Gamma,\p} b'_{(2)})x' = (b \# x)(b' \# x') \\
& = \mathcal{I}_{\Gamma}((b \odo 1)\Gamma(x))\mathcal{I}_{\Gamma}((b' \odo 1)\Gamma(x'))
\end{align*}
for all $b,b'\in B$ and $x,x'\in X$. Additionally in the involutive case, given $b\in B$ and $x\in X$, consider the elements $e\in A$ such that $e \brhd b^* = b^*$, hence we have
\begin{align*}
\mathcal{I}_{\Gamma}(((b \odo 1)\Gamma(x))^{*}) & = \mathcal{I}_{\Gamma}((x^{*}_{[-1]}e \brhd b^{*}) \odo 1)\Gamma(x^{*}_{[0]})) = (x^{*}_{[-1]}e \brhd b^{*}) \# x^{*}_{[0]} \\
& = b^{*}_{(1)} \# \p(x^{*}_{[1]}e,b^{*}_{(2)})x^{*}_{[0]} = b^{*}_{(1)} \# (x^{*} \lhd_{\Gamma,\p} b^{*}_{(2)}) \\
& = (b \# x)^{*} = \mathcal{I}_{\Gamma}((b \odo 1)\Gamma(x))^{*}.
\end{align*}
This end the proof.
\fin


\subsection{Twisted tensor products and coproducts of multiplier Hopf algebras}
Given two regular multiplier Hopf algebras $(A,\com_{A})$, $(B,\com_{B})$. An $\ku$-linear map $T : B \odo A \to A \odo B$ is called a {\em twisting map on $B \odo A$} if
\begin{itemize}
\item[(1)] $T(\id_{B} \odo \m_{A}) = (\m_{A} \odo \id_{B})(\id_{A} \odo T)(T \odo \id_{A})$,
\item[(2)] $T(\m_{B} \odo \id_{A}) = (\id_{A} \odo \m_{B})(T \odo \id_{B})(\id_{B} \odo T)$.
\end{itemize}
and an algebra isomorphism $T' : A \odo B \to B \odo A$ is called {\em cotwisting map on $A \odo B$} if
\begin{itemize}
\item[(3)] $(\id_{B} \odo \com_{A})T' = (T' \odo \id_{A})(\id_{A} \odo T')(\com_{A} \odo \id_{B})$,
\item[(4)] $(\com_{B} \odo \id_{A})T' = (id_{B} \odo T')(T' \odo \id_{B})(\id_{A} \odo \com_{B})$.
\end{itemize}

\begin{rema}
It is evident that the flip map $\Sigma_{A,B}: A \odo B \to B \odo A$ is a cotwisting map on $A \odo B$ and the flip map $\Sigma^{-1}_{A,B}=\Sigma_{B,A}: B \odo A \to A \odo B$ is a twisting map on $B \odo A$.
\end{rema}

\begin{rema}
The right hand side of the equations $(3), (4)$ above are well defined because isomorphisms of non-degenerate algebras can be uniquely extended to isomorphisms of its multiplier algebras.
The requirements of a cotwisting map $T$ on $A \odo B$ with respect to $\com_{A}$ and $\com_{B}$ are dual to the requirements of a twisting map on $B \odo A$ with respect to the multiplications $\m_{A}$ and $\m_{B}$.
\end{rema}


Given a twisting map $T: B \odo A \to A \odo B$, we define the map $\m_{T}:= (\m_{A} \odo \m_{B})(\id_{A} \odo T \odo \id_{B}): A \odo B \odo A \odo B \to A \odo B$. Similarly, given a cotwisting map $T' : A \odo B \to B \odo A$, we define the map $\com_{T} := (\id_{\M(A)} \odo T \odo \id_{\M(B)})(\com_{A} \odo \com_{B}): A \odo B \to \M(A \odo B \odo A \odo B)$. We have the next theorem due to L. Delvaux:

\begin{theo}[\cite{D03}]
Let $(A,\com_{A})$ and $(B,\com_{B})$ be two regular multiplier Hopf algebras.
\begin{enumerate}[label=\textup{(\arabic*)}]
\item If $T:B \odo A \to A \odo B$ is a invertible twisting map satisfying the conditions
\begin{enumerate}[label=\textup{(\roman*)}]
\item $(\id_{A} \odo \Sigma \odo \id_{B})(\com_{A} \odo \com_{B})T = (T \odo T)(\id_{B} \odo \Sigma \odo \id_{A})(\com_{B} \odo \com_{A})$
\item $\com_{\Sigma}T =  \com(\id \odo T \odo \id)(\com_{B} \odo \com_{A})$
\end{enumerate}
then $(A \times_{T} B, \com_{\Sigma})$ yields a regular multiplier Hopf algebra, called the {\em twisted tensor product of $(A,\com_{A})$ and $(B,\com_{B})$}. Here, $A \times_{T} B$ denotes the $\ku$-module $A \odo B$ endowed with the non-degenerate multiplication $\m_{T}$. In the involutive case, given a twisting map $T: B \odo A \to A \odo B$ such that $T(* \odo *)\Sigma_{A,B}T(* \odo *)\Sigma_{A,B} = \id_{A} \odo \id_{B}$, or equivalently $T(* \odo *)\Sigma_{A,B} = \Sigma_{B,A}(* \odo *)T^{-1}$, then $(A \times_{T} B,\com_{\Sigma})$ is a multiplier Hopf \as-algebra.

\item If $T:A \odo B \to B \odo A$ is a cotwisting map then $(A \times_{\Sigma} B, \com_{T})$ yields a regular multiplier Hopf algebra, called the {\em twisted tensor coproduct of $(A,\com_{A})$ and $(B,\com_{B})$}. Here, $A \times_{\Sigma} B$ denotes the $\ku$-module $A \odo B$ endowed with the canonical non-degenerate multiplication $\m_{\Sigma}$, i.e. the tensor product algebra $A \odo B$. In the involutive case, given a cotwisting map $T: B \odo A \to A \odo B$ such that $T$ is an \as-isomorphism, then $(A \odo B,\com_{T})$ is a multiplier Hopf \as-algebra.
\end{enumerate}
\end{theo}

\subsubsection*{Skew-pairings and twisted tensor coproducts}
Given two regular multiplier Hopf algebras $(A,\com_{A})$ and $(B,\com_{B})$. An invertible multiplier $W$ in $\M(B \odo A)$ is called a {\em skew-copairing} if
\begin{itemize}
\item[(i)] $(\id_{B} \odo \com_{A})(W) = W_{12}W_{13}$ in $\M(B \odo A \odo A)$,
\item[(ii)] $(\com_{B} \odo \id_{A})(W) = W_{23}W_{13}$ in $\M(B \odo B \odo A)$.
\end{itemize}

Consider now the linear map $\sigma_{W}: A \odo B \to B \odo A, a \odo b \mapsto W(b \odo a)W^{-1}$. Following the Delvaux's work~\cite{D04}, the map $\sigma_{W}$ is a cotwisting map. Then, it is possible construct the twisted tensor coproduct $(A \odo B,\com_{W})$, which is the algebra $A \odo B$ endowed with the comultiplication given by
\[
\begin{array}{lccc}
\com_{W} : & A \odo B & \to & M(A \odo B \odo A \odo B) \\
& a \odo b & \mapsto & (\id \odo \sigma_{W} \odo \id)(\com_{A}(a) \odo \com_{B}(b))
\end{array}
\]

In the involutive case if $W$ is an unitary element, i.e. $W^{*} = W^{-1}$, then the map $\sigma_{W}$ is a cotwisting map between the multiplier Hopf \as-algebras $(A,\com_{A})$ and $(B,\com_{B})$, thus $(A \odo B,\com_{W})$ is a multiplier Hopf \as-algebra.

\subsection{The Drinfeld double and codouble associated with a pairing}

Let $\p: A \times B \to \ku$ be an admissible pairing between two regular multiplier Hopf algebras. There is a invertible twisting map $T_{\p}: B^{\op} \odo A \to A \odo B^{\op}$ defining an associative non-degenerate multiplication on $A \odo B^{\op}$, this map is given by
\[
T_{\p}(b^{\op} \odo a) = (S_{B}(b_{(3)}) \brhd a \blhd b_{(1)}) \odo b^{\op}_{(2)}
\]
for $a \in A$ and $b \in B$. The Drinfeld double of the pairing $\p$ is the twisted tensor product on $A \odo B^{\op}$ associate to the twisting map $T_{\p}$, endowed with the comultiplication given by
\[
\com_{D(\p)}(a \bicroi b^{\op}) = i_{A}(\com_{A}(a))i_{B^{\op}}(\com_{B^{\op}}(b^{\op}))
\]
for $a \in A$ and $b \in B$. Here, we are using the extension of the canonical non-degenerate homomorphisms $i_{A} : a \odo a' \in A \odo A \mapsto (a \bicroi 1) \odo (a' \bicroi 1) \in \M(A \bicroi B^{\op} \odo A \bicroi B^{\op})$ and $i_{B^{\op}} : b^{\op} \odo b'^{\op} \in B^{\op} \odo B^{\op} \mapsto (1 \bicroi b^{\op}) \odo (1 \bicroi b'^{\op}) \in \M(A \bicroi B^{\op} \odo A \bicroi B^{\op})$. The Drinfeld double is a regular multiplier Hopf algebra and it will be denoted by $\D(\p) = A \bicroi B^{\op}$. 

\begin{rema}
In \cite{DVD04_1}, the authors define the Drinfeld double of a pairing as the bicrossed product $\D(A,B) := A \bicroi B^{\co}$. Observe that this regular multiplier Hopf algebra and $\D(\p)$ are isomorphic using the map $\id \bicroi S_{B}: \D(\p) \to \D(A,B)$. The choice of our characterization is motivated by the definition of the Drinfeld double for quantum groups in the operator theory setting.
\end{rema}

Consider the canonical multiplier $U$ associated to the pairing $\p:A \times B \to \ku$. Consider the multiplier $U^{\ops} := ({}^{\op} \odo \id)(U) \in \M(A^{\op} \odo B)$. The map $\sigma_{U^{^{\ops}}}: A^{\op} \odo B \to B \odo A^{\op}$, defined by $\sigma_{U^{^{\ops}}} = \Sigma\ad(U^{\ops})$, is an invertible cotwisting map, defining a comultiplication on $A^{\op} \odo B$ given by
\[
\begin{array}{lccc}
\com_{U^{^{\ops}}} : & A^{\op} \odo B & \to & M(A^{\op} \odo B \odo A^{\op} \odo B) \\
& a^{\op} \odo b & \mapsto & (\id \odo \sigma_{U^{^{\ops}}} \odo \id)(\com_{A}(a) \odo \com_{B}(b))
\end{array}
.
\]
The Drinfeld codouble of the pairing $\p$ is the twisted tensor coproduct $\T(\p) := A^{\op} \odo B$ associated to the cotwisting map $\sigma_{U^{^{\ops}}}$.

\subsubsection*{Duality} Let $\p: A \times B \to \ku$ be an admissible pairing between two regular multiplier Hopf algebras with canonical multiplier $U \in \M(A \odo B)$. The $\ku$-bilinear map $\hat{\p}: B^{\co} \times A^{\op} \to \ku$, defined by $\hat{\p}(b^{\co},a^{\op}) = \p(a,b)$ for all $a\in A$, $b\in B$, yields an admissible pairing with canonical multiplier $V = (\id \odo {}^{\op})(\Sigma(U)) \in \M(B^{\co} \odo A^{\op})$. 
Denote by $i_{\D(\hat{\p})}: B^{\co} \to \M(\D(\hat{\p}))$ and $i_{\D(\hat{\p})}: A \to \M(\D(\hat{\p}))$ the canonical inclusions in $\M(\D(\hat{\p}))$, i.e. we have $i_{\D(\hat{\p})}(b^{\co})i_{\D(\hat{\p})}(a) = b^{\co} \bicroi a$ for all $a\in A$ and $b\in B$.

\begin{prop}\label{prop:pairing_drinfeld}
The $\ku$-bilinear map $\mathds{P}: \D(\hat{\p}) \times \T(\p) \to \ku$ given by
\[
\mathds{P}(b^{\co} \bicroi a,a'^{\op} \odo b') = \hat{\p}(b^{\co},a'^{\op})\p(a,b'),
\]
for all $a,a' \in A$ and $b,b' \in B$, yields an admissible pairing between the Drinfeld double $\D(\hat{\p})$ and the Drinfeld codouble $\T(\p)$. Moreover, the canonical multiplier associated to the pairing $\mathds{P}$ is given by $\W = \V_{12}\U_{13}$, where $\V = (i_{\D(\hat{\p})} \odo \id_{A^{\op}})(V)$ and $\U = (i_{\D(\hat{\p})} \odo \id_{B})(U)$.

\end{prop}
\pr
It is straightforward show that $\mathds{P}$ is a pairing. Now, fix $a,a' \in A$ and $b,b' \in B$. Take $e \in B$ and $f^{\op} \in A^{\op}$ such that $e \brhd a' = a'$ and $f^{\op} \brhd b'^{\co} = b'^{\co}$. Denote $U(1 \odo e) = d \odo d'$ and $V(1^{\co} \odo f^{\op}) = c^{\co} \odo c'^{\op}$. Then

\begin{align*}
c^{\co} \bicroi d \odo c'^{\op} \odo d' & = (c^{\co} \bicroi 1 \odo c'^{\op} \odo 1)(1 \bicroi d \odo 1 \odo d') \\
& = (i \odo \id)(c^{\co} \odo c'^{\op})_{12}(i \odo \id)(d \odo d')_{13} \\
& = (i \odo \id)(V(1^{\co} \odo f^{\op}))_{12}(i \odo \id)(U(1 \odo e))_{13} \\
& = (i \odo \id)(V)_{12}(i \odo \id)(U)_{13}(1^{\co} \bicroi 1 \odo f^{\op} \odo e) \\
& = \V_{12}\U_{13}(1^{\co} \bicroi 1 \odo f^{\op} \odo e) = \W(1^{\co} \bicroi 1 \odo f^{\op} \odo e)
\end{align*}
and
\begin{align*}
\mathds{P}^{2}(\W,a'^{\op} \odo b' \odo b^{\co} \bicroi a) & = \mathds{P}^{2}(\W(1^{\co} \bicroi 1 \odo f^{\op} \odo e),a'^{\op} \odo b' \odo b^{\co} \bicroi a) \\
& = \mathds{P}^{2}(c^{\co} \bicroi d \odo c'^{\op} \odo d',a'^{\op} \odo b' \odo b^{\co} \bicroi a) \\
& = \mathds{P}(c^{\co} \bicroi d,a'^{\op} \odo b')\mathds{P}(b^{\co} \bicroi a,c'^{\op} \odo d') \\
& = \hat{\p}(c^{\co},a'^{\op})\p(d,b')\hat{\p}(b^{\co},c'^{\op})\p(a,d') \\
& = \hat{\p}^{2}(c^{\co} \odo c'^{\op},a'^{\op} \odo b^{\co})\p^{2}(d \odo d',b' \odo a) \\
& = \hat{\p}^{2}(V(1^{\co} \odo f^{\op}),a'^{\op} \odo b^{\co})\p^{2}(U(1 \odo e),b' \odo a) \\
& = \hat{\p}^{2}(V,a'^{\op} \odo b^{\co})\p^{2}(U,b' \odo a) \\
& = \hat{\p}(b^{\co},a'^{\op})\p(a,b') \\
& = \mathds{P}(b^{\co} \bicroi a,a'^{\op} \odo b').
\end{align*}
The last implies that $\W$ is the canonical multiplier of the pairing $\mathds{P}$. The other conditions to be an admissible pairing follows from the fact that $\W$ is constructed using the canonical multiplier $U$ of the admissible pairing $\p$.
\fin

The admissible pairing $\mathds{P}$ and its canonical multiplier $\W$ play an important role in the equivalence between actions and coactions between the Drinfeld double and the Drinfeld codouble.

\begin{prop}\label{prop:taipe_duality_right_drinfeld}
There is a correspondence one-to-one between right $(\D(\hat{\p}),\com_{\D})$-module algebras and left $(\T(\p),\com_{\T})$-comodule algebras. More explicitly, given a left coaction $\Gamma: X \to \M(\T(\p) \odo X)$ of $(\T(\p),\com_{\T})$ on an algebra $X$, the map $\lhd_{\Gamma,\mathds{P}}: X \odo \D(\hat{\p}) \to X$, defined by
\[
x \lhd_{\Gamma,\mathds{P}} (b^{\co}\bicroi a) = (\hat{\p}(b^{\co},\,\cdot\,) \odo \p(a,\,\cdot\,) \odo \id)(\Gamma(x)),
\]
yields a right action of $(\D(\hat{\p}),\com_{\D})$ on $X$. On the other hand, given a right action $\lhd: X \odo \D(\hat{\p}) \to X$ of $(\D(\hat{\p}),\com_{\D})$ on an algebra $X$, the map $\Gamma_{\lhd,\mathds{P}}: X \to \M(\T(\p) \odo X)$, defined by
\[
((a^{\op} \odo b) \odo 1)\Gamma_{\lhd,\mathds{P}}(x) = ((a^{\op} \odo b) \odo x) \lhd \Sigma(\W),
\]
yields a left coaction of $(\T(\p),\com_{\T})$ on $X$.
\end{prop}
\pr
By the proposition above the pairing $\mathds{P}$ is admissible. Hence, the result follows directly from the duality actions and coactions associated with  an admissible pairing, see section \ref{subsec:duality_act_coact}.
\fin

\subsection{Yetter--Drinfeld algebras over multiplier Hopf algebras}

We recall some important results about the standard characterization of Yetter--Drinfeld algebras in the framework of multiplier Hopf algebras \cite{D07}.

Let $(A,\com)$ be a regular multiplier Hopf algebra, and $X$ be a non-degenerate $(A,\com)$-module algebra for the right unital action $\lhd: X \odo A \to X$ and a $(A^{\op},\com^{\op})$-comodule algebra for the left coaction $\Gamma: X \to \M(A^{\op} \odo X)$. We use the Sweedler type leg notation $x^{\op}_{[-1]} \odo x_{[0]} := \Gamma(x)$ for $x \in X$. 

\begin{defi}
We say that $(X,\lhd,\Gamma)$ is a {\em (right-left) Yetter--Drinfeld $(A,\com)$-algebra} if we have the following condition
\begin{equation}\label{eq:rl-YD-s-1}\tag{s-YD}
a'x_{[-1]}a_{(1)} \odo (x_{[0]} \lhd a_{(2)}) = a'a_{(2)}(x \lhd a_{(1)})_{[-1]} \odo (x \lhd a_{(1)})_{[0]}
\end{equation}
for all $x \in X$, $a,a' \in A$. The above equality is called {\em the Yetter--Drinfeld compatibility condition} for $(X,\lhd,\Gamma)$. A (right-left) Yetter--Drinfeld $(A,\com)$-algebra $(X,\lhd,\Gamma)$ is called {\em braided commutative} or {\em $A$-commutative}, if
\begin{equation}\label{eq:rc-bc-s-1}\tag{s-BC}
xy = (y \lhd x_{[-1]})x_{[0]} \quad \text{or equivalently} \quad xy = y_{[0]}(x \lhd S(y_{[-1]}))
\end{equation}
for all $x, y \in X$.
\end{defi}

\begin{rema}
Similar to the Hopf algebra case, sometimes we denoted $(X,\lhd,\Gamma) \in {}^{A}\aYD_{A}$ for say that $X$ is a right-left Yetter--Drinfeld $(A,\com)$-algebra. We will give another equivalent condition in the next proposition.
\end{rema}

\begin{prop}
The tuple $(X,\lhd,\Gamma)$ is a (right-left) Yetter--Drinfeld $(A,\com)$-algebra if and only if
\begin{equation}\label{eq:rl-YD-s-2}\tag{s-YD'}
\Gamma(x \lhd a)(a'^{\op} \odo 1) = (a'S^{-1}(a_{(3)})x_{[-1]}a_{(1)})^{\op} \odo (x_{[0]} \lhd a_{(2)})
\end{equation}
for all $x \in X$, $a, a' \in A$.
\end{prop}
\pr
\begin{itemize}
\item[($\ida$)] Fix $a,a' \in A$ and $x \in X$. Take $e \in A$ such that $S^{-1}(a_{(3)})e = S^{-1}(a_{(3)})$, then
\small
\begin{align*}
(a'S^{-1}(a_{(3)})x_{[-1]}a_{(1)})^{\op} \odo (x_{[0]} \lhd a_{(2)}) & = (a'S^{-1}(a_{(3)})ex_{[-1]}a_{(1)})^{\op} \odo (x_{[0]} \lhd a_{(2)}) \\
& = ({}^{\op} \odo \id)(\m \odo \id)(a'S^{-1}(a_{(3)}) \odo ex_{[-1]}a_{(1)} \odo (x_{[0]} \lhd a_{(2)})) \\
& = ({}^{\op} \odo \id)(\m \odo \id)(a'S^{-1}(a_{(3)}) \odo ea_{(2)}(x \lhd a_{(1)})_{[-1]} \odo (x \lhd a_{(1)})_{[0]}) \\
& = (a'S^{-1}(a_{(3)})a_{(2)}(x \lhd a_{(1)})_{[-1]})^{\op} \odo (x \lhd a_{(1)})_{[0]} \\
& = \cou(a_{(2)})\Gamma(x \lhd a_{(1)})(a'^{\op} \odo 1) \\
& = \Gamma(x \lhd a)(a'^{\op} \odo 1)
\end{align*}
\normalsize

\item[($\vuelta$)] Fix $a,a' \in A$ and $x \in X$. Take $e \in A$ such that $a_{(2)}e = a_{(2)}$, then
\small
\begin{align*}
a'a_{(2)}(x \lhd a_{(1)})_{[-1]} \odo (x \lhd a_{(1)})_{[0]} & = a'a_{(2)}e(x \lhd a_{(1)})_{[-1]} \odo (x \lhd a_{(1)})_{[0]} \\
& = (\m \odo \id)(\id \odo {}^{\op} \odo \id)(a'a_{(2)} \odo \Gamma(x \lhd a_{(1)})(e^{\op} \odo 1)) \\
& = (\m \odo \id)(a'a_{(2)} \odo e S^{-1}(a_{(1)_{(3)}})x_{[-1]}a_{(1)_{(1)}} \odo (x_{[0]} \lhd a_{(1)_{(2)}})) \\
& = a'a_{(2)}S^{-1}(a_{(1)_{(3)}})x_{[-1]}a_{(1)_{(1)}} \odo (x_{[0]} \lhd a_{(1)_{(2)}}) \\
& = a'a_{(3)_{(2)}}S^{-1}(a_{(3)_{(1)}})x_{[-1]}a_{(1)} \odo (x_{[0]} \lhd a_{(2)}) \\
& = a'x_{[-1]}a_{(1)} \odo (x_{[0]} \lhd \cou(a_{(3)})a_{(2)}) \\
& = a'x_{[-1]}a_{(1)} \odo (x_{[0]} \lhd a_{(2)}).
\end{align*}
\normalsize
\end{itemize}
\fin

\subsection{Yetter--Drinfeld algebras and module algebras over the Drinfeld double}

We recall an adapted version of Theorem~{2.6} found in \cite{D05}.

\begin{theo}\label{th:delvaux_duality_right}
Let $\p: A \odo B \to \ku$ be an admissible pairing between two regular multiplier Hopf algebras and $X$ be a non-degenerate algebra. Consider the admissible pairing $\hat{\p}: B^{\co} \times A^{\op} \to \ku$, $(b^{\co},a^{\op}) \mapsto \p(a,b)$. It hold
\begin{enumerate}[label=\textup{(\arabic*)}]
\item $X$ is a right $\D(\hat{\p})$-module algebra with unital action $\lhd_{\D(\hat{\p})}$ if and only if $X$ is a right $(A,\com)$-module algebra with action
\[
\begin{array}{lccc}
\lhd : & X \odo A & \to & X \\
& x \odo a & \mapsto & x \lhd_{\D(\hat{\p})} (1^{\co}_{\M(B)} \bicroi a)
\end{array}
\]
and $X$ is a left $(A^{\op},\com^{\op})$-comodule algebra with coaction $\Gamma : X \to \M(A^{\op} \odo X)$ such that for all $x \in X$ and $b \in B$, we have
\begin{equation}\label{eq:YD_arising_DD}
(\id_{A^{\op}} \odo (\;\cdot\; \lhd a))(\Gamma(x)) = (S^{-1}(a_{(1)})^{\op} \odo 1_{\M(X)})\Gamma(x \lhd a_{(2)})(a^{\op}_{(3)} \odo 1_{\M(X)}).
\end{equation}
\item If all the algebras involved are involutives and $\p$ is an admissible pairing of multiplier Hopf \as-algebras, then the right action of $\D(\hat{\p})$ on $X$ is compatible with the \as-operation on $X$ in the sense that for all $a \in A$, $b \in B$ and $x \in X$
\[
(x \lhd_{\D(\hat{\p})} (b^{\co} \bicroi a))^{*} = x^{*} \lhd_{\D(\hat{\p})} S_{\D(\hat{\p})}(b^{\co} \bicroi a)^{*} 
\]
if and only if
\[
(x \lhd a)^{*} = x^{*} \lhd S(a)^{*}, \quad \Gamma(x^{*}) = \Gamma(x^{*}).
\]
\end{enumerate}
\end{theo}

\begin{rema}
One can prove that the condition~\eqref{eq:YD_arising_DD} is equivalent to the condition~\eqref{eq:rl-YD-s-2}. 
\end{rema}

With the last remark, the above theorem can be expressed in form of the following equivalence.

\begin{coro}\label{co:taipe_duality_right}
The functor
\[
\begin{array}{rccc}
\mathfrak{E}_{T} : & {}^{A}\aYD_{A} & \leadsto & \alg_{\D(\hat{\p})} \\
& (X,\lhd,\Gamma) & \mapsto & (X,\lhd_{\D(\hat{\p})}) \\
\end{array}
\]
yields an equivalence of categories. Where, given $(X,\lhd,\Gamma) \in {}^{A}\aYD_{A}$, the right action of $\D(\hat{\p})$ on $X$ is given by
\[
x \lhd_{\D(\hat{\p})} (b^{\co} \bicroi a) = (\hat{\p}(b^{\co},\cdot) \odo \id)(\Gamma(x)) \lhd a
\]
for all $x \in X$, $a \in A$ and $b \in B$.

\end{coro}

\section{The ``only coaction'' characterization of Yetter--Drinfeld algebras}\label{sec:oc_yd}

Let $\p: A \times B \to \ku$ be an admissible pairing between two regular multiplier Hopf algebras $(A,\com_{A})$ and $(B,\com_{B})$. The canonical multiplier of $\p$ will be denoted by $U$. Recall that $U$ is the unique invertible element in $\M(A \odo B)$ such that $\p^{2}(U,b \odo a) = \p(a,b)$ for all $a \in A$, $b \in B$. In all this section, we fix the notation $U^{\ops} := ({}^{\op} \odo \id)(U) \in \M(A^{\op} \odo B)$ for a given canonical multiplier $U \in \M(A \odo B)$.

\subsection{Left-left Yetter--Drinfeld algebras}


Let $\alp : X \to \M(A^{\op} \odo X)$ be a left coaction of $(A^{\op},\com^{\op}_{A})$ on $X$ and $\beta : X \to \M(B \odo X)$ be a left coaction of $(B,\com_{B})$ on $X$.

\begin{defi}
We say that $(X,\alp,\beta)$ is a {\em left-left Yetter--Drinfeld algebra over $\p$} if
\begin{equation}\tag{YD}
(\id_{\M(B)} \odo \alp)\beta = \Sigma_{12}\ad(U^{\ops}_{12})(\id_{\M(A^{\op})} \odo \beta)\alp.
\end{equation}
In other words, if the diagram
\[
\xymatrix@C=6pc@R=2pc{X\ar[r]^-{\beta}\ar[d]_-{\alp} & \M(B \odo X)\ar[r]^-{\id_{\M(B)} \odo \alp} & \M(B \odo A^{\op} \odo X) \\ \M(A^{\op} \odo X)\ar[r]_-{\id_{\M(A^{\op})} \odo \beta} & \M(A^{\op} \odo B \odo X)\ar[r]_-{\ad(U^{\ops}) \odo \id_{\M(X)}} & \M(A^{\op} \odo B \odo X)\ar[u]_-{\Sigma \odo \id_{\M(X)}}}
\]
is commutative. Sometimes, we will write $(X,\alp,\beta) \in \aYD^{\,ll}(\p)$.
\end{defi}

\begin{prop}\label{prop:equivalence_st_oc_ll}
The tuple $(X,\alp,\beta)$ is a left-left Yetter--Drinfeld algebra over $\p$ if and only if the tuple $(X,\lhd_{\beta,\p},\alp)$ is a right-left Yetter--Drinfeld $(A,\com_{A})$-algebra. Here, the map $\lhd_{\beta,\p} : X \odo A \to X$, $x \odo a \mapsto ({}_{a}\p \odo \id)(\beta(x))$ %
is the right action of $(A,\com_{A})$ on $X$ induced by the left coaction $\beta$ and the pairing $\p$.
\end{prop}
\pr
Fix $a,a' \in A$ and $x \in X$. For every $a'' \in A$ and $x' \in X$, it is obvious that
\[
(a'' \odo x')\alp(x \lhd_{\beta,\p} a)(a' \odo 1) = (a'' \odo x')({}_{a}\p \odo \id \odo \id)((\id \odo \alp)(\beta(x)))(a' \odo 1).
\]
Moreover, if we use the notation $b \odo y = (1 \odo x')\beta(x_{[0]}) \in B \odo X$, then
\footnotesize
\begin{align*}
(a'' \odo x')(S^{-1}_{A}(a_{(3)})x_{[-1]}a_{(1)}a' \odo (x_{[0]} \lhd_{\beta,\p} a_{(2)})) & =
a''S^{-1}_{A}(a_{(3)})x_{[-1]}a_{(1)}a' \odo x'(x_{[0]} \lhd_{\beta,\p} a_{(2)}) \\
& = a''S^{-1}_{A}(a_{(3)})x_{[-1]}a_{(1)}a' \odo (\p(a_{(2)},\cdot) \odo \id)((1 \odo x')\beta(x_{[0]})) \\
& = a''S^{-1}_{A}(a_{(3)})x_{[-1]}a_{(1)}a' \odo (\p(a_{(2)},\cdot) \odo \id)(b \odo y) \\
& = a''S^{-1}_{A}(a_{(3)})x_{[-1]}a_{(1)}\p(a_{(2)},b)a' \odo y \\
& = a''S^{-1}_{A}(a_{(2)})x_{[-1]}(b \brhd a_{(1)})a' \odo y \\
& = (a'' \odo 1)((\id \odo {}_{a}\p)T_{U}(x_{[-1]} \odo b)a' \odo y) \qquad\qquad\qquad\qquad\qquad \text{by } \ref{eq:t_p_1} \\
& = (a'' \odo 1)({}_{a}\p \odo \id \odo \id)(\Sigma T_{U} \odo \id)(x_{[-1]} \odo b \odo y)(a' \odo 1) \\
& = (a'' \odo 1)({}_{a}\p \odo \id \odo \id)(\Sigma T_{U} \odo \id)(x_{[-1]} \odo (1 \odo x')\beta(x_{[0]}))(a' \odo 1) \\
& = (a'' \odo x')({}_{a}\p \odo \id \odo \id)(\Sigma T_{U}({}^{\op} \odo \id) \odo \id)(\id \odo \beta)\alp(x))(a' \odo 1)
\end{align*}
\normalsize
This two equalities and the non-degeneracy of the product on $A \odo X$ implies that the equality
\[
\alp(x \lhd_{\beta,\p} a)(a' \odo 1) = ((S^{-1}_{A}(a_{(3)})x_{[-1]}a_{(1)}a')^{\op} \odo (x_{[0]} \lhd_{\beta,\p} a_{(2)}))
\]
holds for all $a,a' \in A$, $x \in X$ if and only if
\begin{align*}
(\id_{\M(B)} \odo \alp)\beta & = (\Sigma({}^{\op} \odo \id_{\M(B)})T_{U}({}^{\op} \odo \id_{\M(B)}) \odo \id_{\M(X)})(\id_{\M(A^{\op})} \odo \beta)\alp \\
& = (\Sigma\ad(U^{\ops}) \odo \id_{\M(X)})(\id_{\M(A^{\op})} \odo \beta)\alp.
\end{align*}
\fin

\subsection{Braided commutativity}\label{sec:ll-mod-bc}

By Proposition~\ref{prop:equivalence_st_oc_ll}, have a left-left Yetter--Drinfeld algebra $(X,\alp,\beta)$ over a pairing $\p:A \times B \to \ku$ is equivalent to have a (right-left) Yetter--Drinfeld $(A,\com_{A})$-algebra, namely $(X,\lhd_{\beta,\p},\alp)$. Here, we study the braided commutativity condition in the setting of a left-left Yetter--Drinfeld algebra $(X,\alp,\beta)$ over $\p$.

Recall that a (right-left) Yetter--Drinfeld $(A,\com_{A})$-algebra $(X,\lhd_{\beta,\p},\alp)$ is called {\em braided commutative} if $\m_{X}\tau_{X} = \m_{X}$, where $\tau_{X}: X \odo X \to X \odo X$ is the map defined by
\[
\tau_{X}(x \odo y) = (y \lhd_{\beta,\p} x_{[-1]}) \odo x_{[0]}
\]
for all $x, y \in X$, or equivalently if $\m_{X}\rho_{X} = \m_{X}$, where $\rho_{X}: X \odo X \to X \odo X$ is the map defined by
\[
\rho_{X}(x \odo y) = y_{[0]} \odo (x \lhd_{\beta,\p} S_{A}(y_{[-1]}))
\]
for all $x, y \in X$. Observe that $\rho_{X}\tau_{X} = \tau_{X}\rho_{X} = \id_{X \odo X}$.

Using the $\ku$-linear maps $\p: A \odo B \to \ku$ and $\p_{+}: A \odo B \to \ku$, defined by $\p_{+}(a \odo b) := \p(S_{A}(a),b) = \p(a,S_{B}(b))$, it is straightforward to show that
\[
\tau_{X} = (\p \odo \id \odo \id)\Sigma_{34}\Sigma_{23}(\dot\alp \odo \beta) \quad \text{ and } \quad \rho_{X} = (\p_{+} \odo \id \odo \id)\Sigma_{12}\Sigma_{23}\Sigma_{24}(\beta \odo \dot\alp)
\]
as maps from $X \odo X$ to $X \odo X$.

Before to give our definition of braided commutative in the ``only coaction'' characterization of Yetter--Drinfeld algebras, we recall some constructions and we fix some notations. The algebra $\He(\p) \odo X^{\op}$ is endowed with the non-degenerate multiplication given by
\[
\m_{\He(\p) \odo X^{\op}} = (\m_{\He(\p)} \odo \m_{X^{\op}})(\id_{\He(\p)} \odo \Sigma_{\He(\p),X^{\op}} \odo \id_{\M(X)^{\op}}),
\]
here $\He(\p) = A \# B$ denotes the Heisenberg algebra of the pairing $\p$ and
\[
\m_{\He(\p)}((a \# b) \odo (a' \# b')) = a(b_{(1)} \brhd a') \# b_{(2)}b' = aa'_{(1)} \# (b \blhd a'_{(2)})b'
\]
for all $a,a'\in A$, $b,b'\in B$. Denote by $\iota_{A,X^{\op}} : A \odo X^{\op} \to \M(\He(\p) \odo X^{\op})$, $a \odo x \mapsto a \# 1_{\M(B)} \odo x$ and by $\iota_{B,X^{\op}}: B \odo X^{\op} \to \M(\He(\p) \odo X^{\op})$, $b \odo x \mapsto 1_{\M(A)} \# b \odo x$, the canonical inclusions induced by the inclusions of $A$ and $B$ in the multiplier Heisenberg algebra $\M(\He(\p))$. Note that the canonical pairing can be regarded also as a $\ku$-linear map on $\He(\p)$, i.e. $\p: A \# B \to \ku$, $a \# b \mapsto \p(a,b)$.

\begin{defi}\label{df:ll-bc-mod}
A left-left Yetter--Drinfeld algebra $(X,\alp,\beta)$ over $\p$ will be called {\em braided commutative} if
\begin{equation}\label{eq:ll-bc-mod}\tag{BC}
\m_{\He(\p) \odo X^{\op}}(\iota_{A,X^{\op}} \odo \iota_{B,X^{\op}})(\alp^{\op}(x^{\op}) \odo \beta^{\co}(y^{\op})) = \m_{\He(\p) \odo X^{\op}}(\iota_{B,X^{\op}} \odo \iota_{A,X^{\op}})(\beta^{\co}(y^{\op}) \odo \alp^{\op}(x^{\op}))
\end{equation}
for each $x,y \in X$. In this case, we will denote $(X,\alp,\beta) \in \aYD^{\,ll}_{\textrm{bc}}(\p)$.
\end{defi}

In the following remark will give a more simple notation to express the braided commutativity condition. This simply notation will be used mainly in computations.

\begin{rema}
Recall the commutator notation of two elements in an algebra: Given an algebra $X$ and two elements $x,y \in X$, the commutator of $x$ and $y$ is the element $[x,y] := xy -yx$. With this notation in mind, a left-left Yetter--Drinfeld algebra $(X,\alp,\beta)$ over $\p$ is braided commutative if and only if for each $x,y \in X$, we have
\begin{equation}\label{eq:ll-bc-mod-2}\tag{BC'}
[\alp^{\op}(x^{\op}),\beta^{\co}(y^{\op})] = 0
\end{equation}
inside the algebra $\M(\He(\p) \odo X^{\op})$.
\end{rema}

As we expect this definition of a braided commutative Yetter--Drinfeld algebra over a pairing is equivalent to the braided commutativity condition in the standard characterization of Yetter--Drinfeld algebras over multiplier Hopf algebras as we prove in the next proposition. 

\begin{prop}\label{prop:equiv_bc}
Consider a pairing $\p:A \times B \to \ku$ between two regular multiplier Hopf algebras. The left-left Yetter--Drinfeld algebra $(X,\alp,\beta)$ over $\p$ is braided commutative if and only if the (right-left) Yetter--Drinfeld $(A,\com_{A})$-algebra $(X,\lhd_{\beta,\p},\alp)$ is braided commutative.
\end{prop}
\pr
In order to simplify the proof, we will suppose that $A$ and $B$ are unital Hopf algebras. Given $x,y \in X$, consider the Sweedler type leg notations $x^{\op}_{[-1]} \odo x_{[0]} := \alpha(x)$ and $y^{[-1]} \odo y^{[0]} := \beta(y)$. On one hand, it holds
\footnotesize
\begin{align*}
\m_{\He(\p) \odo X^{\op}}(\iota_{A,X^{\op}} \odo \iota_{B,X^{\op}})(\alp^{\op}(x) \odo \beta^{\co}(y)) & = \m_{\He(\p) \odo X^{\op}}(x_{[-1]} \# 1 \odo x^{\op}_{[0]} \odo 1 \# S_{B}(y^{[-1]}) \odo y^{[0]^{\op}}) \\
& = (\id_{\He(\p)} \odo {}^{\op}\m_{X})(x_{[-1]} \# S_{B}(y^{[-1]}) \odo y^{[0]} \odo x_{[0]})
\intertext{\normalsize and on the other}
\m_{\He(\p) \odo X^{\op}}(\iota_{B,X^{\op}} \odo \iota_{A,X^{\op}})(\beta^{\co}(y) \odo \alp^{\op}(x)) & = \m_{\He(\p) \odo X^{\op}}(1 \# S_{B}(y^{[-1]}) \odo y^{[0]^{\op}} \odo x_{[-1]} \# 1 \odo x^{\op}_{[0]}) \\
& = (1 \# S_{B}(y^{[-1]}))(x_{[-1]} \# 1) \odo (x_{[0]}y^{[0]})^{\op} \\
& = (S_{B}(y^{[-1]})_{(1)} \brhd x_{[-1]}) \# S_{B}(y^{[-1]})_{(2)} \odo (x_{[0]}y^{[0]})^{\op} \\
& = x_{[-1]_{(1)}} \# S_{B}(y^{[-1]}_{(1)}) \odo \p_{+}(x_{[-1]_{(2)}} \odo y^{[-1]}_{(2)})(x_{[0]}y^{[0]})^{\op} \\
& = x_{[-1]} \# S_{B}(y^{[-1]}) \odo \p_{+}(x_{[0]_{[-1]}} \odo y^{[0]^{[-1]}})(x_{[0]_{[0]}}y^{[0]^{[0]}})^{\op} \\
& = x_{[-1]} \# S_{B}(y^{[-1]}) \odo (\p_{+} \odo {}^{\op}\m_{X})(x_{[0]_{[-1]}} \odo y^{[0]^{[-1]}} \odo x_{[0]_{[0]}} \odo y^{[0]^{[0]}}) \\
& = x_{[-1]} \# S_{B}(y^{[-1]}) \odo (\p_{+} \odo {}^{\op}\m_{X})\Sigma_{12}\Sigma_{34}\Sigma_{23}(\beta \odo \dot\alp)(y^{[0]} \odo x_{[0]}) \\
& = x_{[-1]} \# S_{B}(y^{[-1]}) \odo {}^{\op}\m_{X}(\p_{+} \odo \id \odo \id)\Sigma_{34}\Sigma_{23}(\beta \odo \dot\alp)(y^{[0]} \odo x_{[0]}) \\
& = x_{[-1]} \# S_{B}(y^{[-1]}) \odo {}^{\op}\m_{X}\rho_{X}(y^{[0]} \odo x_{[0]}) \\
& = (\id_{\He(\p)} \odo {}^{\op}\m_{X}\rho_{X})(x_{[-1]} \# S_{B}(y^{[-1]}) \odo y^{[0]} \odo x_{[0]}).
\end{align*}
\normalsize

Now, we have
\begin{itemize}
\item[($\vuelta$)] If $(X,\lhd_{\beta,\p},\alp)$ is braided commutative, we have $\m_{X} = \m_{X}\tau_{X}$, follows immediately from the equalities above that the condition~(\ref{eq:ll-bc-mod}) is true for every $x,y \in X$, i.e the tuple $(X,\alp,\beta)$ is braided commutative.

\item[($\ida$)] Let $x',y' \in X$. Take $y \odo x = \tau_{X}\Sigma(y' \odo x')$, thus we have $y' \odo x' = \Sigma\rho_{X}(y \odo x)$. If $(X,\alp,\beta)$ satisfies the condition~(\ref{eq:ll-bc-mod}), then the equalities
\footnotesize
\begin{align*}
\qquad\qquad(\p \odo \id_{\M(X)^{\op}})\m_{\He(\p) \odo X^{\op}}(\iota_{A,X^{\op}} \odo \iota_{B,X^{\op}})(\alp^{\op}(x) \odo \beta^{\co}(y)) & = (\p \odo {}^{\op}\m_{X})(x_{[-1]} \# S_{B}(y^{[-1]}) \odo x_{[0]} \odo y^{[0]}) \\
& = {}^{\op}\m_{X}(\p(x_{[-1]} \# S_{B}(y^{[-1]}))y^{[0]} \odo x_{[0]}) \\
& = {}^{\op} \m_{X}(\p(S_{A}(x_{[-1]}) \# y^{[-1]})y^{[0]} \odo x_{[0]}) \\
& = {}^{\op}\m_{X}((y \lhd_{\beta,\p} S_{A}(x_{[-1]})) \odo x_{[0]})) \\
& = {}^{\op}\m_{X}\Sigma\rho_{X}(x \odo y) \\
& = {}^{\op}\m_{X}(y' \odo x')
\end{align*}
and
\begin{align*}
\qquad(\p \odo \id_{\M(X)})\m_{\He(\p) \odo X^{\op}}(\iota_{B,X^{\op}} \odo \iota_{A,X^{\op}})(\beta^{\co}(y) \odo \alp^{\op}(x)) & = (\p \odo {}^{\op}\m_{X}\rho_{X})(x_{[-1]} \# S_{B}(y^{[-1]}) \odo y^{[0]} \odo x_{[0]}) \\
& = \m_{X}\rho_{X}\Sigma\rho_{X}(x \odo y) \\
& = \m_{X}\rho_{X}(y' \odo x'),
\end{align*}
\normalsize
implies $\m_{X}(y' \odo x') = \m_{X}\rho_{X}(y' \odo x')$. Thus $(X,\lhd_{\beta,\p},\alp)$ is braided commutative.
\end{itemize}
\fin

We have an equivalent condition for the braided commutativity using the Heisenberg algebra of the flip pairing $\p$. Denote by $\du{\iota}_{A,X} : A \odo X \to \M(\He(\bar{\p}) \odo X)$, $a \odo x \mapsto 1_{\M(B)} \# a \odo x$ and by $\du{\iota}_{B,X}: B \odo X \to \M(\He(\bar{\p}) \odo X)$, $b \odo x \mapsto b \# 1_{\M(A)} \odo x$, the canonical inclusions induced by the inclusions of $A$ and $B$ in the multiplier Heisenberg algebra $\M(\He(\bar{\p})) = \M(B \# A)$.

\begin{prop}\label{prop:dual-ll-bc-mod}
A left-left Yetter--Drinfeld algebra $(X,\alp,\beta)$ over $\p$ is braided commutative if and only if
\[
\m_{\He(\bar{\p}) \odo X}(\du{\iota}_{A,X} \odo \du{\iota}_{B,X})((S_{A} \odo \id)\dot\alp(x) \odo \beta(y)) = \m_{\He(\bar{\p}) \odo X}(\du{\iota}_{B,X} \odo \du{\iota}_{A,X})(\beta(y) \odo (S_{A} \odo \id)\dot\alp(x))
\]
for every $x,y \in X$.
\end{prop}
\pr
Consider the anti-isomorphism $\mathcal{L}: \He(\bar{\p}) \to \He(\p)$, $b \# a \mapsto S_{A}^{-1}(a) \# S_{B}(b)$ introduced in Proposition~\ref{prop:iso_heisenberg}. Observe that
\[
(\mathcal{L} \odo {}^{\op})\hat\iota_{A,X}(S_{A} \odo \id)\dot\alp(x) = \iota_{A,X^{\op}}\alp^{\op}(x^{\op}) \quad \text{ and } \quad (\mathcal{L} \odo {}^{\op})\hat\iota_{B,X}\beta(x) = \iota_{B,X^{\op}}\beta^{\co}(x^{\op})
\]
for all $x \in X$. Because $\mathcal{L} \odo {}^{\op}: \He(\bar{\p}) \odo X \to \He(\p) \odo X^{\op}$ is an anti-isomorphism, it hold
\small
\begin{align*}
(\mathcal{L} \odo {}^{\op})\m_{\He(\bar{\p}) \odo X}(\du{\iota}_{A,X} \odo \du{\iota}_{B,X})((S_{A} \odo \id)\dot\alp(x) \odo \beta(y)) = 
\m_{\He(\p) \odo X^{\op}}(\iota_{A,X^{\op}} \odo \iota_{B,X^{\op}})(\alp^{\op}(x^{\op}) \odo \beta^{\co}(y^{\op}))
\end{align*}
\normalsize
and
\small
\begin{align*}
(\mathcal{L} \odo {}^{\op})\m_{\He(\bar{\p}) \odo X}(\du{\iota}_{B,X} \odo \du{\iota}_{A,X})(\beta(y) \odo (S_{A} \odo \id)\dot\alp(x)) = \m_{\He(\p) \odo X^{\op}}(\iota_{B,X^{\op}} \odo \iota_{A,X^{\op}})(\beta^{\co}(y^{\co}) \odo \alp^{\op}(x^{\co}))
\end{align*}
\normalsize
for every $x, y \in X$. The result follows from the above equalities.
\fin

\subsection{Duality}

Let $\p: A \times B \to \ku$ be an admissible pairing between two regular multiplier Hopf algebras. By Proposition~\ref{prop:op_co_coaction}, $\alp: X \to \M(A^{\op} \odo X)$ is a left coaction of $(A^{\op},\com^{\op}_{A})$ on $X$ if and only if $\alp^{\ops}: X^{\op} \to \M(A \odo X^{\op})$ is a left coaction of $(A,\com_{A})$ on $X^{\op}$. Similarly, $\beta: X \to \M(B \odo X)$ is a left coaction of $(B,\com_{B})$ on $X$ if and only if  $\beta^{\ops}: X^{\op} \to \M(B^{\op} \odo X^{\op})$ is a left coaction of $(B^{\op},\com^{\op}_{B})$ on $X^{\op}$.

\begin{prop}\label{prop:duality_ydbc}
The following statements are equivalent
\begin{enumerate}[label=\textup{(\roman*)}]
\item $(X,\alp,\beta) \in \aYD^{ll}(\p)$;
\item $(X^{\op},\beta^{\op},\alp^{\op}) \in \aYD^{ll}(\bar{\p})$.
\end{enumerate}
Moreover, $(X,\alp,\beta)$ is braided commutative if and only if $(X^{\op},\beta^{\op},\alp^{\op})$ is braided commutative.
\end{prop}
\pr
By symmetry, it suffices to prove that (i) implies (ii). First observe that $U$ is the canonical multiplier of the pairing $\p$ if and only if $\bar{U} := \Sigma(U)$ is the canonical multiplier of the flip pairing $\bar{\p}$. Because
\[
\bar{U}^{\ops} = ({}^{\op} \odo \id)(\Sigma(U)) = ({}^{\op} \odo {}^{\op})(\Sigma(U^{\ops})),
\]
we have
\[
\Sigma\circ\ad(\bar{U}^{\ops})\circ({}^{\op} \odo {}^{\op}) = ({}^{\op} \odo {}^{\op})\circ\ad(U^{\ops})^{-1}\circ\Sigma.
\]
Now, if $(X,\alpha,\beta)$ is a Yetter--Drinfeld algebra over $\p$, thus
\begin{align*}
(\id \odo \beta^{\op})\circ\alpha^{\op} & = ({}^{\op} \odo {}^{\op} \odo {}^{\op})\circ(\id \odo \beta)\circ\alpha\circ{}^{\op} \\
& = ({}^{\op} \odo {}^{\op} \odo {}^{\op})\circ\ad(U^{\ops}_{12})^{-1}\circ\Sigma_{12}\circ(\id \odo \alpha)\circ\beta\circ{}^{\op} \\
& = \Sigma_{12}\circ\ad(\bar{U}^{\ops}_{12})\circ({}^{\op} \odo {}^{\op} \odo {}^{\op})\circ(\id \odo \alpha)\circ\beta\circ{}^{\op} \\
& = \Sigma_{12}\circ\ad(\bar{U}^{\ops}_{12})\circ(\id \odo \alpha^{\op})\circ\beta^{\op}
\end{align*}
which says that $(X^{\op},\beta^{\op},\alpha^{\op})$ is a Yetter--Drinfeld algebra over $\bar{\p}$.

Finally, take $\tilde\alpha := \beta^{\op}: X^{\op} \to \M(B^{\op} \odo X^{\op})$ and $\tilde\beta = \alpha^{\op}: X^{\op} \to \M(A \odo X^{\op})$. Observe that with the last notations, it follows $\tilde\alpha^{\op} = \beta$, $\tilde\beta^{\co} = (S_{A} \odo {}^{\op})\circ\tilde\beta\circ {}^{\op} = (S_{A}\circ{}^{\op} \odo \id)\circ\alpha = (S_{A} \odo \id)\dot\alpha$ and then
\begin{equation}\label{eq:yd_dual_1}
\m_{\He(\bar{\p}) \odo X}(\hat\iota_{B,X} \odo \hat\iota_{A,X})(\tilde\alpha^{\op}(x) \odo \tilde\beta^{\co}(y)) = \m_{\He(\bar{\p}) \odo X}(\hat\iota_{B,X} \odo \hat\iota_{A,X})(\beta(x) \odo (S_{A} \odo \id)\dot\alpha(y)),
\end{equation}
\begin{equation}\label{eq:yd_dual_2}
\m_{\He(\bar{\p}) \odo X}(\hat\iota_{A,X} \odo \hat\iota_{B,X})(\tilde\beta^{\co}(y) \odo \tilde\alpha^{\op}(x)) = \m_{\He(\bar{\p}) \odo X}(\hat\iota_{A,X} \odo \hat\iota_{B,X})((S_{A} \odo \id)\dot\alpha(y) \odo \beta(x))
\end{equation}
for all $x,y \in X$. By definition, $(X^{\op},\tilde\alpha,\tilde\beta) = (X^{\op},\beta^{\op},\alpha^{\op}) \in \aYD^{ll}(\bar{\p})$ is braided commutative if and only if
\[
\m_{\He(\bar{\p}) \odo X}(\hat\iota_{B,X} \odo \hat\iota_{A,X})(\tilde\alpha^{\op}(x) \odo \tilde\beta^{\co}(y)) = \m_{\He(\bar{\p}) \odo X}(\hat\iota_{A,X} \odo \hat\iota_{B,X})(\tilde\beta^{\co}(y) \odo \tilde\alpha^{\op}(x))
\]
for all $x,y \in X$. The equivalence with the braided commutative of $(X,\alpha,\beta) \in \aYD^{ll}(\p)$ follows from Proposition~\ref{prop:dual-ll-bc-mod} and using the equations~\eqref{eq:yd_dual_1} and \eqref{eq:yd_dual_2}.
\fin

In the involutive case, we need to have a compatibility between the coactions and the involution on $X^{\op}$.

\begin{prop}\label{prop:duality_ydbc_inv}
Assume that $\alpha\circ\gamma = (S^{-2}_{A^{\op}} \odo \gamma)\circ\alpha$ and $\beta\circ\gamma = (S^{-2}_{B} \odo \gamma)\circ\beta$. Then, the following statements are equivalent
\begin{enumerate}[label=\textup{(\roman*)}]
\item $(X,\alp,\beta) \in \aYD^{ll}(\p)$;
\item $(X^{\op}_{\gamma},\beta^{\op},\alp^{\op}) \in \aYD^{ll}(\bar{\p})$.
\end{enumerate}
Moreover, $(X,\alp,\beta)$ is braided commutative if and only if $(X^{\op}_{\gamma},\beta^{\ops},\alp^{\ops})$ is braided commutative.
\end{prop}
\pr
Follows directly from Proposition~\ref{prop:op_co_coaction_involutive} and Proposition~\ref{prop:duality_ydbc}.
\fin

In particular, in the finite-dimensional Hopf algebra setting, we have

\begin{coro}
Let $A= H$, $B = \du{H}$ and $\p: A \times B \to \ku$ be the canonical pairing. The following statements are equivalent
\begin{enumerate}[label=\textup{(\roman*)}]
\item $(X,\alp,\beta) \in \aYD^{ll}(\p)$;
\item $(X,\lhd_{\beta},\alp) \in {}^{H}\aYD_{H}$;
\item $(X^{\op},\beta^{\op},\alp^{\op}) \in \aYD^{ll}(\bar{\p})$;
\item $(X^{\op},\lhd_{\alp^{\op}},\beta^{\op}) \in {}^{\du{H}}\aYD_{\du{H}}$.
\end{enumerate}
Moreover, if one of the Yetter--Drinfeld algebras in the above equivalence is braided commutative, then all the others are also braided commutative Yetter--Drinfeld algebras.
\end{coro}

\subsection{Equivalences between Yetter--Drinfeld algebras}

In this section, we study equivalences of Yetter--Drinfeld structures under isomorphism of the regular multiplier Hopf algebras.

\begin{nota}
Given an isomorphism of regular multiplier Hopf algebras $f: A \to A'$. We will denote by $f^{\ops}:A^{\op} \to A'^{\op}$ the canonical isomorphism given by $f^{\ops}(a^{\op}) := f(a)^{\op}$ for all $a \in A$.
\end{nota}

\begin{prop}\label{prop:equivalence_YD_pairing}
Let $\alp: X \to \M(A^{\op} \odo X)$ and $\beta: X \to \M(B \odo X)$ be two left coactions of the multiplier Hopf algebras $(A^{\op},\com_{A})$ and $(B,\com_{B})$ on $X$, respectively. For any pair of multiplier Hopf algebras isomorphisms $f: A \to A'$ and $g: B \to B'$, the following statements are equivalent
\begin{enumerate}[label=\textup{(\roman*)}]
\item The tuple $(X,\alp,\beta)$ is a left-left Yetter--Drinfeld algebra over $\p$.
\item The tuple $(X,{}^{f^{\ops}}\alp,\beta)$ is a left-left Yetter--Drinfeld algebra over ${}^{f^{-1}}\p$.
\item The tuple $(X,\alp,{}^{g}\beta)$ is a left-left Yetter--Drinfeld algebra over $\p^{g^{-1}}$.
\end{enumerate}
Here, we are using the notations ${}^{f^{\ops}}\alp := (f^{\ops} \odo \id)\alp$ and ${}^{g}\beta := (g \odo \id)\beta$ for the coactions; and ${}^{f^{-1}}\p: (a',b) \in A' \times B \mapsto \p(f^{-1}(a'),b) \in \ku$ and $\p^{g^{-1}}: (a,b') \in A \times B' \mapsto \p(a,g^{-1}(b')) \in \ku$ for the admissible pairings induced by the isomorphisms $f$ and $g$, respectively.
\end{prop}
\pr
Recall that $U^{\ops} := ({}^{\op} \odo \id)(U) \in \M(A^{\op} \odo B)$. The uniqueness of the canonical multiplier associated to a pairing implies $(f^{\ops} \odo \id)(U^{\ops}) = U^{\ops}({}^{f^{-1}}\p)$ and $(\id \odo g)(U^{\ops}) = U^{\ops}(\p^{g^{-1}})$. Then,
\begin{equation}\label{eq:e1}
\ad(U^{\ops}({}^{f^{-1}}\p)) = \ad((f^{\ops} \odo \id)(U^{\ops})) = (f^{\ops} \odo \id)\ad(U^{\ops})({f^{\ops}}^{-1} \odo \id)
\end{equation}
and
\begin{equation}\label{eq:e2}
\ad(U^{\ops}(\p^{g^{-1}})) = \ad((\id \odo g)(U^{\ops})) = (\id \odo g)\ad(U^{\ops})(\id \odo g^{-1}).
\end{equation}

\noindent The equivalence of the following equalities
\begin{equation}\label{eq:i}\tag{i}
(\id_{\M(B)} \odo \alp)\beta = \Sigma_{12}\ad(U^{\ops}_{12})(\id_{\M(A^{\op})} \odo \beta)\alp,
\end{equation}
\begin{equation}\label{eq:ii}\tag{ii}
(\id_{\M(B)} \odo {}^{f^{\ops}}\alp)\beta = \Sigma_{12}\ad(U^{\ops}({}^{f^{-1}}\p)_{12})(\id_{\M(A^{\op})} \odo \beta){}^{f^{\ops}}\alp,
\end{equation}
\begin{equation}\label{eq:iii}\tag{iii}
(\id_{\M(B)} \odo \alp){}^{g}\beta = \Sigma_{12}\ad(U^{\ops}(\p^{g^{-1}})_{12})(\id_{\M(A^{\op})} \odo {}^{g}\beta)\alp
\end{equation}
follows directly from the equalities~\eqref{eq:e1} and \eqref{eq:e2}. This prove the proposition.
\fin

\subsection{Yetter--Drinfeld algebras and coactions of the Drinfeld codouble}

It is known by the works of D. Radford \cite{R93} and S. Majid \cite{M95} that Yetter--Drinfeld modules over a Hopf algebra $H$ are close related with modules over the Drinfeld double of the same Hopf algebra $H$. In \cite{D05}, L. Delvaux  has shown that under certain conditions on the pairing between two regular multiplier Hopf algebras the relation is still true, see Theorem~\ref{th:delvaux_duality_right}. Here, we prove an equivalent relation using the ``only coaction'' characterization of Yetter--Drinfeld algebras over a pairing.

Let $\p: A \times B \to \ku$ be an admissible pairing between two regular multiplier Hopf algebras. Denote by $\hat{\p} : B^{\co} \times A^{\op} \to \ku $, $(b,a^{\op}) \mapsto \p(a,b)$, the admissible pairing between $(B,\com^{\co}_{B})$ and $(A^{\op},\com^{\op}_{A})$.

\begin{theo}\label{prop:left_equivalence_taipe} 
The functor
\[
\begin{array}{rccc}
\mathfrak{E} : & \aYD^{\,ll}(\p) & \leadsto & {}^{\T(\p)}\alg \\
& (X,\alp,\beta) & \mapsto & (X,\gamma = (\id \odo \beta)\alp) \\
& f: (X,\alp,\beta) \to (X',\alp',\beta') & \mapsto & f: (X,\gamma) \to (X,\gamma')
\end{array}
\]
yields an equivalence of categories. Moreover, we have the following equivalence between categories:
\begin{equation*}\label{eq:1}
\xymatrix@C=5pc@R=1.5pc{%
\shadowbox*{%
\begin{Bcenter}
$\;\aYD^{\,ll}(\p)\;$ \\[0.1cm]
\begin{Bitemize}
\item $(X,\alp,\beta)$
\item $f: (X,\alp,\beta) \to (X',\alp',\beta')$
\end{Bitemize}
\end{Bcenter}}\ar@{<~>}[r]\ar@{<~>}[d]%
&%
\shadowbox*{%
\begin{Bcenter}
$\;{}^{\T(\p)}\alg\;$ \\[0.1cm]
\begin{Bitemize}
\item $(X,\gamma=(\id \odo \beta)\alp)$
\item $f: (X,\gamma) \to (X',\gamma')$
\end{Bitemize}
\end{Bcenter}}\ar@{<~>}[d]%
\\
{\shadowbox*{%
\begin{Bcenter}
$\;{}^{A}\aYD_{A}\;$ \\[0.1cm]
\begin{Bitemize}
\item $(X,\lhd_{\beta,\p},\alp)$
\item $f: (X,\lhd_{\beta,\p},\alp) \to (X,\lhd_{\beta',\p},\alp')$
\end{Bitemize}
\end{Bcenter}}}\ar@{<~>}[r]%
&%
\shadowbox*{%
\begin{Bcenter}
$\;\alg_{\D(\hat{\p})}\;$ \\[0.1cm]
\begin{Bitemize}
\item $(X,\lhd_{\gamma})$
\item $f: (X,\lhd_{\gamma}) \to (X',\lhd_{\gamma'})$
\end{Bitemize}
\end{Bcenter}}}
\end{equation*}
\end{theo}
\pr
Take $(X,\alp,\beta) \in \aYD^{ll}(\p)$ and define $\gamma = (\id_{A^{\op}} \odo \beta)\alp: X \to \M(A^{\op} \odo B \odo X)$. We have
\begin{align*}
(\id_{\T(\p)} \odo \gamma)\gamma & = (\id_{\T(\p)} \odo (\id_{A^{\op}} \odo \beta)\alp)(\id_{A^{\op}} \odo \beta)\alp \\
& = (\id_{\T(\p)} \odo \id \odo \beta)(\id_{\T(\p)} \odo \alp)(\id \odo \beta)\alp \\
& = (\id_{\T(\p)} \odo \id \odo \beta)(\id_{A^{\op}} \odo (\id \odo \alp)\beta)\alp \\
& = (\id_{\T(\p)} \odo \id \odo \beta)(\id_{A^{\op}} \odo (\Sigma \odo \id)(\ad(U^{\ops}) \odo \id)(\id \odo \beta)\alp)\alp \\
& = ((\Sigma\circ \ad(U^{\ops}))_{23} \odo \beta)(\id \odo \id \odo \beta)(\id \odo \alp)\alp \\
& = ((\Sigma \ad(U^{\ops}))_{23} \odo \beta)(\id \odo \id \odo \beta)(\com_{A} \odo \id)\alp \\
& = ((\Sigma \ad(U^{\ops}))_{23} \odo \id)(\com_{A} \odo \id \odo \beta)(\id \odo \beta)\alp \\
& = ((\Sigma \ad(U^{\ops}))_{23} \odo \id)(\com_{A} \odo \com_{B} \odo \id)(\id \odo \beta)\alp \\
& = ((\id \odo \Sigma \odo \id)(\id \odo \ad(U^{\ops}) \odo \id)(\com_{A} \odo \com_{B}) \odo \id)(\id \odo \beta)\alp \\
& = (\com_{\T(\p)} \odo \id)\gamma
\end{align*}
and follows from the injectivity of $\alp$ and $\beta$ that $\gamma$ is injective. Then $(X,\gamma)$ is a left $\T(\p)$-comodule algebra. Moreover, if we compute the right action of $\D(\hat{\p})$ on $X$ induced by $\gamma$ and associated to the pairing $\mathds{P} = \hat{\p} \odo \p$, we have
\begin{align*}
x \lhd_{\gamma,\mathds{P}} (b^{\co} \bicroi a) & = ((\hat{\p} \odo \p)(b^{\co} \bicroi a, \cdot) \odo \id )(\gamma(x)) \\
& = ((\hat{\p} \odo \p)(b^{\co} \bicroi a, \cdot) \odo \id)(\id \odo \beta)(\alp(x)) \\
& = (\hat{\p}(b^{\co},\cdot) \odo (\p(a,\cdot) \odo \id)\beta)(\alp(x)) \\
& = (\p(a,\cdot) \odo \id)\beta)(x \lhd_{\alp,\hat{\p}} b^{\co}) \\
& = (x \lhd_{\alp,\hat{\p}} b^{\co}) \lhd_{\beta,\p} a
\end{align*}
for every $b\in B$, $a\in A$ and $x \in X$, i.e. $(X,\lhd_{\gamma,\mathds{P}})=\mathfrak{E}_{T}(X,\lhd_{\beta,\p},\alp)$ by Corollary~\ref{co:taipe_duality_right}. Now, take $(X,\gamma) \in {}^{\T(\p)}\alg$. By Proposition~\ref{prop:taipe_duality_right_drinfeld}, we have $(X,\lhd_{\gamma,\mathds{P}}) \in \alg_{\D(\hat{\p})}$, hence by Corollary~\ref{co:taipe_duality_right}, there is a right action $\lhd: X \odo A \to X$ and a left coaction $\alp: X \to \M(A^{\op} \odo X)$ such that $\mathfrak{E}_{T}(X,\lhd,\alp) = (X,\lhd_{\gamma})$. For all $x \in X$, we have
\begin{align*}
(\id \odo \beta_{\lhd})(\alp(x)) & = (\id \odo \beta_{\lhd})(U^{\op}_{1} \odo x \lhd_{\gamma} (U_{2} \bicroi 1)) = U^{\op}_{1} \odo \beta_{\lhd}(x \lhd_{\gamma} (U_{2} \bicroi 1))\\
& = U^{\op}_{1} \odo U_{2} \odo (x \lhd_{\gamma} (U_{2} \bicroi 1)) \lhd_{\gamma} (1 \bicroi U_{1}) \\
& = U^{\op}_{1} \odo U_{2} \odo (x \lhd_{\gamma} (U_{2} \bicroi U_{1})) \\
& = \gamma(x),
\end{align*}
because $\W = (i \odo {}^{\op})(\Sigma(U))_{12}(i \odo \id)(U)_{13} = ``U_{2} \bicroi U_{1} \odo U^{\op}_{1} \odo U_{2}"$ is the canonical multiplier for the pairing $\mathds{P} = \hat{\p} \odo \p: \D(\hat{\p}) \times \T(\p) \to \ku$, where $U = ``U_{1} \odo U_{2}" \in \M(A \odo B)$ is the canonical multiplier associated to the pairing $\p: A \times B \to \ku$. Finally, the commutative diagram between categories follows from Proposition \ref{prop:equivalence_st_oc_ll} and Theorem \ref{th:delvaux_duality_right}.
\fin

\section{Examples}\label{sec:exa}

\subsection{Transformation groups}

Let $\cdot: G \times S \to S$ be a left action of a (possibly infinite) group $G$ on a (possibly infinite) set $S$. Consider the two canonical multiplier Hopf algebras associated to $G$, see Example~\ref{example_group}, $(K(G),\com)$ and $(\C[G],\com')$. They admit a canonical pairing $\p_{G}:K(G) \times \ku[G] \to \ku$, $\p(f,\lambda_{g}) = f(g)$ which is admissible with canonical multiplier associated to $\p_{G}$ given by
\[
U_{G}=\sum_{g \in G}\delta_{g} \odo \lambda_{g} \in \M(K(G) \odo \ku[G]).
\]
Consider the left coaction of $(K(G),\com)$ on the algebra $K(S)$ induced by the action of $G$ on $S$
\[
\begin{array}{lccc}
\alp : & K(S) & \to & \M(K(G) \odo K(S)) \\ 
& p & \mapsto & \displaystyle\sum_{g \in G,\, s \in S}p(g \cdot s)\delta_{g} \odo \delta_{s} 
\end{array}
\]
and the trivial left coaction of $(\ku[G],\dcom)$ on $K(S)$, $\beta : K(S) \to \M(\ku[G] \odo K(S))$, $p \mapsto \lambda_{e}  \odo p$.

\begin{lemm}
The tuple $(K(S),\alp,\beta)$ is a braided commutative Yetter--Drinfeld algebra over $\p_{G}$.
\end{lemm}
\pr
First by direct computation, we have
\begin{align*}
(\ad(U_{G}) \odo \id)(\id \odo \beta)\alp(\delta_{s}) & = (\ad(U_{G}) \odo \id)\left(\sum_{g \in G} \delta_{g} \odo \lambda_{e} \odo \delta_{g^{-1} \cdot s}\right) \\
& = \sum_{g \in G} \ad(U_{G})(\delta_{g} \odo \lambda_{e}) \odo \delta_{g^{-1} \cdot s}  \\
& =
\sum_{g,h \in G} (\delta_{h} \odo \lambda_{h})(\delta_{g} \odo \lambda_{e})(\delta_{h} \odo \lambda_{h^{-1}}) \odo \delta_{g^{-1} \cdot s} \\
& = (\Sigma \odo \id)\left(\sum_{g \in G} \lambda_{e} \odo \delta_{g} \odo \delta_{g^{-1} \cdot s}\right) \\
& = (\Sigma \odo \id)(\id \odo \alp)(\lambda_{e} \odo \delta_{s}) \\
& = (\Sigma \odo \id)(\id \odo \alp)\beta(\delta_{s}),
\end{align*}
for all $s \in S$. Now observe that $\He(\p_{G})$ is the $\ku$-module $K(G) \# \ku[G]$ with algebra structure given by $\lambda_{g}\delta_{g'} = \delta_{g'g^{-1}}\lambda_{g}$ for all $g,g' \in G$. Thus inside the algebra $\M(\He(\p_{G}) \odo K(S))$, it holds
\begin{align*}
[\alp(\delta_{s}),(S' \odo \id)\beta(\delta_{s'})] & = \displaystyle\sum_{g\in G}(\delta_{g} \odo \delta_{g^{-1}\cdot s})(\lambda_{e} \odo \delta_{s'}) - (\lambda_{e} \odo \delta_{s'})(\delta_{g} \odo \delta_{g^{-1}\cdot s}) \\
& = \displaystyle\sum_{g\in G}(\delta_{g}\lambda_{e} - \lambda_{e}\delta_{g}) \odo \delta_{g^{-1}\cdot s}\delta_{s'} = 0,
\end{align*}
for all $s,s' \in S$. The last equality implies the braided commutative of the Yetter--Drinfeld algebra $(K(S),\alpha,\beta)$.
\fin

In the following lines, using the present example, we want to illustrate the motivation for the use of Yetter--Drinfeld algebras in order to construct quantum transformation groupoids. Denote by $G \ltimes S$ the transformation groupoid $G \times S \rightrightarrows S$ associated to the action $\cdot: G \times S \to S$. Recall that the unit space of this groupoid is $(G \ltimes S)^{(0)} = \{e\}\times S$, its domain and its range maps are $d: (g,x) \in G \ltimes S \mapsto s \in S$, $r: (g,s) \in G \times S \mapsto g\cdot s \in S$, respectively. The composable pairs space of this groupoid is
\[
(G \ltimes S)^{(2)} = \{((g,s),(g',s')) \in (G \ltimes S)^{2} : g'\cdot s' = s\} \subseteq (G\ltimes S) \times (G\ltimes S)
\]
and its groupoid operations are given by
\[
(g,g'\cdot t)(g',s') = (gg',s'), \qquad (g,s)^{-1} = (g^{-1},g\cdot s)
\]
for all $g,g' \in G$ and $s,s' \in S$.

\begin{lemm}
Consider the braided commutative Yetter--Drinfeld algebra $(K(S),\alpha,\beta)$.
\begin{enumerate}[label=\textup{(\arabic*)}]
\item We have
\[
(\ku[G],\com') \ltimes_{\beta,\p_{G}} K(S) \iso K(G) \#_{\lhd_{\beta,\p_{G}}} K(S) \iso K(G \ltimes S),
\]
where $K(G \ltimes S)$ denotes the (\as-)algebra of finitely supported $\ku$-valued functions on the transformation groupoid $G \ltimes S$.

\item We have
\[
(K(G),\com) \ltimes_{\alpha,\p_{G}} K(S) \iso \ku[G]\,\#_{\lhd_{\alpha,\p_{G}}}\,K(S) \iso \ku[G \ltimes S],
\]
where $\ku[G \ltimes S]$ denotes the groupoid (\as-)algebra of the transformation groupoid $G \ltimes S$.
\end{enumerate}
\end{lemm}
\pr
\begin{enumerate}[label=\textup{(\arabic*)}]
\item Consider the coaction $\beta : K(S) \to \M(\ku[G] \odo K(S))$, $p \mapsto \lambda_{e} \odo p$. The dual right action of $\beta$ is given by the trivial action
\[
\begin{array}{lccc}
\lhd_{\beta} : & K(S) \odo K(G) & \to & K(S) \\
& \delta_{s} \odo \delta_{g} & \mapsto & (\p_{G}(\delta_{g},\cdot) \odo \id)(\beta(\delta_{s})) = \delta_{g,e}\,\delta_{s}
\end{array}
\]
then, the algebraic crossed product associated with the action $\beta$ is given by 
\[
\begin{array}{ccccc}
(\ku[G],\com') \ltimes_{\beta,\p_{G}} K(S) & \cong & K(G) \odo K(S) & \cong  & K(G \ltimes S) \\
(\delta_{g} \odo 1_{\M(K(S))})\beta(\delta_{s}) & \mapsto & \delta_{g} \odo \delta_{s} & \mapsto & \displaystyle \delta_{(g,s)}
\end{array}
.
\]

\item Consider the coaction  $\alpha : K(S) \to \M(K(G) \odo K(S))$. By a direct computation, the dual multiplier Hopf (\as-)algebra action of $\alp$ is given by 
\[
\begin{array}{lccc}
\lhd_{\alp} : & K(S) \odo \ku[G] & \to & K(S) \\
& p \odo \lambda_{g} & \mapsto & (\p_{G}(\cdot,\lambda_{g}) \odo \id)(\alp(p)) = \displaystyle\sum_{s \in S}p(g \cdot s) \delta_{s}
\end{array}.
\]
then the algebraic crossed product associated with the action $\alp$ is
\[
\begin{array}{ccccc}
(K(G),\com) \ltimes_{\alpha} K(S) & \cong & \ku[G]\,\#_{\lhd_{\alpha,\p_{G}}}\,K(S) & \cong  & \ku[G \ltimes S] \\
(\lambda_{g} \odo 1_{\M(K(S))})\alp(p) & \mapsto & \lambda_{g}\,\#\,p& \mapsto & \displaystyle\sum_{s \in S} p(s) \lambda_{(g,s)}
\end{array}
,
\]
where $\ku[G \ltimes S]$ denotes the groupoid (\as-)algebra of $G \ltimes S$, i.e. the algebra generated by formal elements $\{\lambda_{(g,s)}\}_{(g,s) \in G \ltimes S}$, and relations
\[
\lambda_{(g,s)}\lambda_{(g',t)} = \delta_{s,g' \cdot s'}\lambda_{(gg',s')}
\]
for all $(g,s),(g',s') \in G \ltimes S$. For the involutive case, in $\ku[G \ltimes S]$ we also have the relation
\[
\lambda^{*}_{(g,s)} = \lambda_{(g,s)^{-1}} = \lambda_{(g^{-1},g \cdot s)}
\]
for all $(g,s) \in G\ltimes S$. 
\end{enumerate}
\fin
%

\begin{rema}
The lemma above only shows how using the coactions of the Yetter--Drinfeld algebra $(K(S),\alpha,\beta)$, we arrive to obtain two algebras associated with the transformation groupoid $G \ltimes S$. In order to extend the complete structure of the transformation groupoid, we need to use again the coactions but in a different way. This not will be done here, it will be treated in more detail in a separated work~\cite{Ta22_3}. The interested reader can also consult \cite{ET16} for the treatment of this construction in the operator algebra setting.
\end{rema}

\subsection{Trivial Yetter--Drinfeld algebras}

Let $\p: A \times B \to \ku$ be an admissible pairing of two regular multiplier Hopf algebras. Consider a non-degenerate algebra $X$ endowed with a trivial left coaction of $(A^{\op},\com^{\op}_{A})$ denoted by $\alpha_{\tr}$ and a trivial left coaction of $(B,\com_{B})$ denoted by $\beta_{\tr}$. Because
\[
(\id \odo \alpha_{\tr})\beta_{\tr}(x) = 1_{\M(B)} \otimes 1^{\op}_{\M(A)} \otimes x, \qquad (\id \odo \beta_{\tr})\alpha_{\tr}(x) = 1^{\op}_{\M(A)} \otimes 1_{\M(B)} \otimes x
\]
for all $x \in X$, and
\[
\ad(U^{\ops})(1^{\op}_{\M(A)} \otimes 1_{\M(B)}) = 1^{\op}_{\M(A)} \otimes 1_{\M(B)},
\]
then $(X,\alpha_{\tr},\beta_{\tr})$ is a Yetter--Drinfeld algebra over $\p$. Moreover, it is not hard to see that $(X,\alpha_{\tr},\beta_{\tr})$ is braided commutative if and only if $X$ is a commutative algebra.

\subsection{Actions of multiplier Hopf algebras}

Let $\p: A \times B \to \ku$ be an admissible pairing of two regular multiplier Hopf algebras, $(X,\alpha)$ be a left $(A^{\op},\com^{\op}_{A})$-comodule algebra and $(X,\beta_{\tr})$ be a left $(B,\com_{B})$-comodule algebra using the trivial coaction. If $A$ is commutative, then $(X,\alpha,\beta_{\tr})$ is a Yetter--Drinfeld algebra over $\p$. Moreover, $(X,\alpha,\beta_{\tr})$ is braided commutative if and only if $X$ is a commutative algebra.



\subsection{The Heisenberg algebra over the Drinfeld double}
\label{sec:heis-algebra-over-mu}

Let $\p: A \times B \to \ku$ be an admissible pairing of two regular multiplier Hopf algebras. Consider the Heisenberg algebra associated to the pairing $\p$, $\He(\p) = A \# B$, and the Drinfeld double associated to the pairing $\hb{\p}: A^{\co} \times B^{\op} \to \ku$, $\hb{\p}(a^{\co},b^{\op}) = \p(a,b)$, $\D(\hb{\p}) = A^{\co} \bicroi B \iso A \bicroi B^{\op} = \D(\p)$.

By Theorem 3.5 in~\cite{YZC17}, if we consider the left unital action
\[
\begin{array}{rccc}
\rhd : & \D(\hb{\p}) \odo \He(\p) & \to & \He(\p) \\
& (a^{\co} \bicroi b) \odo (a' \# b') & \mapsto & a_{(3)}(b_{(1)} \brhd a')S^{-1}_{A}(a_{(2)}) \# (b_{(2)}b'S_{B}(b_{(3)})) \blhd S^{-1}_{A}(a_{(1)})
\end{array}
\]
and the left coaction $\Gamma': \He(\p) \to \M(\D(\hb{\p}) \odo \He(\p))$ defined by
\[
((a'^{\co} \bicroi b') \odo 1)\Gamma'(a \# b) = (a'^{\co} \bicroi b')(a^{\co}_{(2)} \bicroi b_{(1)}) \odo a_{(1)} \# b_{(2)}
\]
for all $a,a'\in A$ and $b,b'\in B$, then the tuple $(\He(\p),\rhd,\Gamma')$ is a left-left Yetter--Drinfeld $(\D(\hb{\p}),\com_{\D(\hb{\p})})$-algebra.

In what follows, we give a new way to construct this same Yetter--Drinfeld algebra structure using our ``only coaction'' characterization. In other words, in this section we will prove the next theorem.

\begin{theo}
The tuple $(\He(\p),\Gamma',\Gamma)$ is a left-right Yetter--Drinfeld algebra over the natural pairing $\bar{\mathds{P}} = \hb{\p} \odo \bar{\p}: \D(\hb{\p}) \times \T(\bar{\p}) \to \ku$ between $\D(\hb{\p})$ and $\T(\bar{\p})$. Moreover, because
\[
(\He(\p),\Gamma',\Gamma) \in \aYD^{\,lr}(\bar{\mathds{P}}),
\]
then
\[
(\He(\p),\rhd_{\Gamma,\bar{\mathds{P}}},\Gamma') \in {}^{\D(\hb{\p})}_{\D(\hb{\p})}\aYD.
\]
\end{theo}

\subsubsection*{The left coaction of $\D(\hb{\p})$ on $\He(\p)$}


Consider the Lu's isomorphism $\bar{L}: \D(\hb{\p}) \to \He(\p)$, $a^{\co} \bicroi b \mapsto a \# b$. Hence, the map
\[
\begin{array}{lccc}
\Gamma' := (\id_{\D(\hb{\p})} \odo \bar{L})\com_{\D(\hb{\p})} \bar{L}^{-1}: & \He(\p) & \to & \M(\D(\hb{\p}) \odo \He(\p))
\end{array}
\]
defines a left coaction of the Drinfeld double $\D(\hb{\p})$ on $\He(\p)$. Indeed,
\begin{align*}
(\id_{\D(\hb{\p})} \odo \Gamma')\Gamma' & = (\id_{\D(\hb{\p})} \odo (\id_{\D(\hb{\p})} \odo \bar{L})\com_{\D(\hb{\p})} \bar{L}^{-1})(\id_{\D(\hb{\p})} \odo \bar{L})\com_{\D(\hb{\p})} \bar{L}^{-1} \\
& = (\id_{\D(\hb{\p})} \odo \id_{\D(\hb{\p})} \odo \bar{L})(\id_{\D(\hb{\p})} \odo \com_{\D(\hb{\p})})\com_{\D(\hb{\p})} \bar{L}^{-1} \\
& = (\com_{\D(\hb{\p})} \odo \id_{\He(\p)})(\id_{\D(\hb{\p})} \odo \bar{L})\com_{\D(\hb{\p})} \bar{L}^{-1} \\
& = (\com_{\D(\hb{\p})} \odo \id_{\He(\p)})\Gamma',
\end{align*}
and
\begin{align*}
(\cou_{\D(\hb{\p})} \odo \id_{\He(\p)})\Gamma' & = (\cou_{\D(\hb{\p})} \odo \id_{\He(\p)})(\id_{\D(\hb{\p})} \odo \bar{L})\com_{\D(\hb{\p})} \bar{L}^{-1} \\
& = \bar{L} (\cou_{\D(\hb{\p})} \odo \id_{\He(\p)})\com_{\D(\hb{\p})} \bar{L}^{-1} = \id_{\He(\p)}.
\end{align*}

\begin{rema}
Given $a,a'\in A$ and $b,b'\in B$, we have
\begin{align*}
((a'^{\co} \bicroi b') \odo 1)\Gamma'(a \# b) & = ((a'^{\co} \bicroi b') \odo 1)(\id_{\D(\hb{\p})} \odo \bar{L})\com_{\D(\hb{\p})}(a^{\co} \bicroi b) \\
& = (\id_{\D(\hb{\p})} \odo \bar{L})(((a'^{\co} \bicroi b') \odo 1)\com_{\D(\hb{\p})}(a^{\co} \bicroi b)) \\
& = (\id_{\D(\hb{\p})} \odo \bar{L})((a'^{\co} \bicroi b')(a^{\co}_{(2)} \bicroi b_{(1)}) \odo a^{\co}_{(1)} \bicroi b_{(2)}) \\
& = (a'^{\co} \bicroi b')(a^{\co}_{(
2)} \bicroi b_{(1)}) \odo a_{(1)} \# b_{(2)}, 
\end{align*}
then we can use the Sweedler type leg notation
\[
\begin{array}{lccc}
\Gamma' = (\id_{\D(\hb{\p})} \odo \bar{L})\com_{\D(\hb{\p})} \bar{L}^{-1}: & \He(\p) & \to & \M(\D(\hb{\p}) \odo \He(\p)) \\
& a \# b & \mapsto & a^{\co}_{(2)} \bicroi b_{(1)} \odo a_{(1)} \# b_{(2)}
\end{array}
\]
for the left coaction defined above.
\end{rema}


\subsubsection*{The right coaction of $\T(\bar{\p})$ on $\He(\p)$}

Here, we will consider two canonical right coactions of $B^{\op}$ and $A$ on $\He(\p)$ which allows to construct a right coaction of $\T(\bar{\p})$ on $\He(\p)$. Take the canonical multiplier $U \in \M(A \odo B)$ associated to the pairing $\p$ and the canonical inclusion homomorphisms $\iota_{A}: A \to \M(\He(\p))$, $a \mapsto a \# 1_{\M(B)}$, and $\iota_{B}: B \to \M(\He(\p))$, $b \mapsto 1_{\M(A)} \# b$. Denotes by $U^{\ops}$ the multiplier $(\id \odo {}^{\op})(U) \in \M(A \odo B^{\op})$ which is also the canonical multiplier associated to the pairing $\hb{\p}: A^{\co} \times B^{\op} \to \ku$.

\begin{rema}
In order to simplify notations in the following, we will consider the map $\iota_{A}\iota_{B}: A \odo B \to \He(\p)$ defined by $\iota_{A}\iota_{B}(a \odo b)=\iota_{A}(a)\iota_{B}(b) = a \# b$ for $a\in A$, $b\in B$, and the map $\iota_{\D,A}\iota_{\D,B}: A \odo B \to \D(\hb{\p})$ defined by $\iota_{\D,A}\iota_{\D,B}(a \odo b) = \iota_{\D,A}(a)\iota_{\D,B}(b) = a^{\co} \bicroi b$ for $a\in A$, $b\in B$. Observe that using the Lu's isomorphism, we have $\bar{L}\iota_{\D,A}\iota_{\D,B} = \iota_{A}\iota_{B}$.
\end{rema}

\begin{enumerate}[label=\textup{(\roman*)}]
\item Consider the map
\[
\begin{array}{lccc}
\alp : & \He(\p) & \to & \M(\He(\p) \odo B^{\op}) \\
& a \# b & \mapsto & (\iota_{A} \odo \id)(U^{\ops})(a \# b \odo 1^{\op}_{\M(B)})(\iota_{A} \odo \id)(U^{\ops})^{-1}
\end{array}
\]

\begin{itemize}
\item[($*$)] Denote $\bar{\V} := (\iota_{A} \odo \id_{B^{\op}})(U^{\ops}) \in \M(\He(\p) \odo B^{\op})$. Because
\small
\begin{align*}
(\id_{\He(\p)} \odo \com_{B^{\op}})(\bar{\V}) & = (\id_{\He(\p)} \odo \com_{B^{\op}})(\iota_{A} \odo \id)(U^{\ops}) \\
& = (\iota_{A} \odo \id \odo \id)(\id_{A} \odo \com_{B^{\op}})(U^{\ops}) \\
& = (\iota_{A} \odo \id \odo \id)((U^{\ops})_{12}(U^{\ops})_{13}) = \bar{\V}_{12}\bar{\V}_{13},
\end{align*}
\normalsize
then $\alp$ is a right coaction of $B^{\op}$  on the Heisenberg algebra $\He(\p)$.

\item[($*$)] For all $a,a'\in A$ and $b \in B$, follows from the equality
\small
\begin{align*}
a^{\co} \rhd_{\alp,\hb{\p}} (a' \# b) & = (\id \odo \hb{\p}(a^{\co},\cdot))(\alp(a' \# b)) = (\id \odo \hb{\p}(a^{\co},\cdot))(\bar{\V}(a \# b \odo 1^{\op})\bar{\V}^{-1}) \\
& = (\id \odo \hb{\p}(a^{\co}_{(3)},\cdot))(\bar{\V})(\id \odo \hb{\p}(a^{\co}_{(2)},\cdot))(a' \# b \odo 1^{\op})(\id \odo \hb{\p}(a^{\co}_{(1)},\cdot))(\bar{\V}^{-1}) \\
& = \iota_{A}((\id \odo \hb{\p}(a^{\co}_{(2)},\cdot))(U^{\ops}))(a \# b)\iota_{A}((\id \odo \hb{\p}(a^{\co}_{(1)},\cdot))((U^{\ops})^{-1})) \\
& = \iota_{A}(a_{(2)})(a \# b)\iota_{A}(S^{-1}_{A}(a_{(1)})).
\end{align*}
\normalsize
that, the left action induced by $\alp$ is given by
\[
\begin{array}{lccc}
\rhd_{\alp,\hb{\p}} : & A^{\co} \odo \He(\p) & \to & \He(\p)  \\
& a^{\co} \odo a' \# b & \mapsto & \iota_{A}(a_{(2)})(a' \# b)\iota_{A}(S^{-1}_{A}(a_{(1)}))
\end{array}.
\]
Moreover, we can see that
\small
\begin{equation}\label{eq:alpha}
a^{\co} \rhd_{\alp,\hb{\p}} (a' \# b) = (a_{(2)} \# 1)(a' \# b)(S^{-1}_{A}(a_{(1)}) \# 1) = a_{(3)}a'S^{-1}_{A}(a_{(2)}) \# (b \blhd S^{-1}_{A}(a_{(1)}))
\end{equation}
\normalsize
for all $a,a'\in A$, $b\in B$.
\end{itemize}

\item Consider the linear map
\[
\begin{array}{lccc}
\beta : & \He(\p) & \to & \M(\He(\p) \odo A) \\
& a \# b & \mapsto & (\iota_{B} \odo \id)(\Sigma(U))(a \# b \odo 1_{\M(A)})(\iota_{B} \odo \id)(\Sigma(U))^{-1}
\end{array}
\]

\begin{itemize}
\item[($*$)] Denote $\bar{\U} := (\iota_{B} \odo \id_{A})(\Sigma(U)) \in \M(\He(\p) \odo A)$. Because
\small
\begin{align*}
(\id \odo \com_{A})(\bar{\U}) & = (\id \odo \com_{A})(\iota_{B} \odo \id)(\Sigma(U)) = (\iota_{B} \odo \id \odo \id)(\id \odo \com_{A})(\Sigma(U)) \\
& = (\iota_{2} \odo \id \odo \id)\Sigma_{12}\Sigma_{23}(\com_{A} \odo \id)(U) = (\iota_{2} \odo \id \odo \id)\Sigma_{12}\Sigma_{23}(U_{13}U_{23}) \\
& = (\iota_{2} \odo \id \odo \id)(\Sigma(U)_{12}\Sigma(U)_{13}) = \bar{\U}_{12}\bar{\U}_{13},
\end{align*}
\normalsize
then $\beta$ is a right coaction of $A$ on the Heisenberg algebra $\He(\p)$.

\item[($*$)] Because
\small
\begin{align*}
b \rhd_{\beta,\p} (a \# b') & = (\id \odo \p(\cdot,b))(\beta(a \# b')) = (\id \odo \p(\cdot,b))(\bar{\U}(a \# b' \odo 1)\bar{\U}^{-1}) \\
& = (\id \odo \p(\cdot,b_{(1)}))(\bar{\U})(\id \odo \p(\cdot,b_{(2)}))(a \# b' \odo 1)(\id \odo \p(\cdot,b_{(3)}))(\bar{\U}^{-1}) \\
& = \iota_{B}((\p(\cdot,b_{(1)}) \odo \id)(U))(a \# b')\iota_{B}((\p(\cdot,b_{(2)}) \odo \id)(U^{-1})) \\
& = \iota_{B}(b_{(1)})(a \# b')\iota_{B}(S_{B}(b_{(2)}))
\end{align*}
\normalsize
for any $a\in A$ and $b,b' \in B$. The induced left action by $\beta$ is given by
\[
\begin{array}{lccc}
\rhd_{\beta,\p} : & B \odo \He(\p) & \to & \He(\p)  \\
& b \odo a \# b' & \mapsto & \iota_{B}(b_{(1)})(a \# b')\iota_{B}(S_{B}(b_{(2)}))
\end{array}.
\]
Moreover, we can see that
\begin{equation}\label{eq:beta}
b \rhd_{\beta,\p} (a \# b') = (1 \# b_{(1)})(a \# b')(1 \# S_{B}(b_{(2)})) = (b_{(1)} \brhd a) \# b_{(2)}b'S(b_{(3)}).
\end{equation}
for all $a,\in A$, $b,b'\in B$.
\end{itemize}
\end{enumerate}

\begin{lemm}\label{lemm:hola}
Let $U \in \M(A \odo B)$ be the canonical multiplier of the pairing $\p$. Consider
\[
\bar{\U} = (\iota_{B} \odo \id_{A})(\Sigma(U)) \in \M(\He(\p) \odo A), \qquad \bar{\V} = (\iota_{A} \odo \id_{B^{\op}})(U^{\ops}) \in \M(\He(\p) \odo B^{\op}).
\]
Then
\[
\bar{\V}_{12}\bar{\U}_{13}= (\bar{L} \odo \id_{\T(\bar{\p})})(\bar{\W})
\]
where $\bar{\W} \in \M(\D(\hb{\p}) \odo \T(\bar{\p}))$ is the canonical multiplier of the pairing $\bar{\mathds{P}}$ between the Drinfeld double $\D(\hb{\p})$ and the Drinfeld codouble $\T(\bar{\p})$.
\end{lemm}
\pr
Consider the canonical multiplier $\bar{U} = \Sigma(U) \in \M(B \odo A)$ associated to the pairing $\bar{\p}$. By Proposition~\ref{prop:pairing_drinfeld}, the canonical multiplier of the pairing $\bar{\mathds{P}} = \hb{\p} \odo \bar{\p}: A^{\co} \bicroi B \times B^{\op} \odo A \to \ku$ is given by the element
\begin{align*}
\bar{\W} & = (\iota_{\D,A} \odo \id_{B^{\op}})((\id \odo {}^{\op})(\Sigma(\bar{U})))_{12}(\iota_{\D,B} \odo \id_{A})(\bar{U})_{13} \\
& = (\iota_{\D,A} \odo \id_{B^{\op}})(U^{\ops})_{12}(\iota_{\D,B} \odo \id_{A})(\Sigma(U))_{13}
\end{align*}
in $\M(\D(\hb{\p}) \odo \T(\bar{\p}))$. Thus, we have
\small
\begin{align*}
\bar{\V}_{12}\bar{\U}_{13} & = (\iota_{A}\iota_{B} \odo \id_{B^{\op}} \odo \id_{A})(U^{\ops}_{13}\Sigma(U)_{24}) = (\bar{L}\iota_{\D,A}\iota_{\D,B} \odo \id_{B^{\op}} \odo \id_{A})(U^{\ops}_{13}\Sigma(U)_{24}) \\
& = (\bar{L} \odo \id_{B^{\op}} \odo \id_{A})(\iota_{\D,A}\iota_{\D,B} \odo \id_{B^{\op}} \odo \id_{A})(U^{\ops}_{13}\Sigma(U)_{24}) \\
& = (\bar{L} \odo \id_{B^{\op}} \odo \id_{A})((\iota_{\D,A} \odo \id_{B^{\op}})(U^{\ops})_{12}(\iota_{\D,B} \odo \id_{A})(\Sigma(U))_{13}) \\
& = (\bar{L} \odo \id_{B^{\op}} \odo \id_{A})(\bar{\W}) = (\bar{L} \odo \id_{\T(\bar{\p})})(\bar{\W}).
\end{align*}
\normalsize
\fin

\begin{prop}
Consider the coactions $\alp$ and $\beta$ defined as above. The tuple $(\He(\p),\alp,\beta)$ is a right-right Yetter--Drinfeld algebra over the pairing $\bar{\p}: B\times A \to \ku$. 
\end{prop}
\pr
Take the map $\Gamma := (\alp \odo \id)\beta: \He(\p) \to \M(\He(\p) \odo A \odo B^{\op})$. Consider the invertible elements
\[
\bar{\U} := (\iota_{B} \odo \id)(\Sigma(U)) \in \M(\He(\p) \odo A) \quad \text{ and } \quad \bar{\V} := (\iota_{A} \odo \id)(U^{\ops}) \in \M(\He(\p) \odo B^{\op}).
\]
Because $\alp(x) = \ad(\bar{\V})(x \odo 1^{\op}_{\M(B)})$ and $\beta(x) = \ad(\bar{\U})(x \odo 1_{\M(A)})$ for all $x \in \He(\p)$, then
\[
\Gamma(x) = (\alp \odo \id)\beta(x) = \ad(\bar{\V}_{12}\bar{\U}_{13})(x \odo 1^{\op}_{\M(B)} \odo 1_{\M(A)})
\]
for all $x \in \He(\p)$. Follows from Lemma~\ref{lemm:hola}, that
\begin{align*}
(\id_{\He(\p)} \odo \com_{\T})(\bar{\V}_{12}\bar{\U}_{13}) & = (\id_{\He(\p)} \odo \com_{\T})(\bar{L} \odo \id_{\T(\bar{\p})})(\bar{\W}) \\
& = (\bar{L} \odo \id_{\T(\bar{\p})} \odo \id_{\T(\bar{\p})})(\id_{\D(\hb{p})} \odo \com_{\T})(\bar{\W}) \\
& = (\bar{L} \odo \id_{\T(\bar{\p})} \odo \id_{\T(\bar{\p})})(\bar{\W}_{12}\bar{\W}_{13}) \\
& = (\bar{L} \odo \id_{\T(\bar{\p})})(\bar{\W})_{12}(\bar{L} \odo \id_{\T(\bar{\p})})(\bar{\W})_{13} \\
& = (\bar{\V}_{12}\bar{\U}_{13})_{12}(\bar{\V}_{12}\bar{\U}_{13})_{13}
\end{align*}
and because $\bar{\U}$ and $\bar{\V}$ are invertible multipliers, then $\Gamma$ is a right coaction of $\T(\bar{\p})$ on $\He(\p)$. Follows from the right version of Proposition~\ref{prop:left_equivalence_taipe}, that $(\He(\p),\alp,\beta)$ is a right-right Yetter--Drinfeld algebra over the flip pairing $\bar{\p}: B\times A \to \ku$.
\fin

\begin{coro}
$(\He(\p),\rhd_{\beta,\p},\alp)$ is a left-right Yetter--Drinfeld algebra over the multiplier Hopf algebra $(B,\com_{B})$.
\end{coro}
\pr
Follows directly from the right version of Proposition~\ref{prop:equivalence_st_oc_ll}.
\fin

\begin{coro}\label{coro:hola}
There is a right coaction of $\T(\bar{\p})$ on $\He(\p)$ given by
\[
\begin{array}{lccc}
\Gamma : & \He(\p) & \to & \M(\He(\p) \odo \T(\bar{\p})) \\
& a \# b & \mapsto & (\alp \odo \id)(\beta(a \# b))
\end{array}
\]
which using the pairing $\bar{\mathds{P}}$ induced a left action of $\D(\hb{\p})$ on $\He(\p)$ given by
\[
\begin{array}{lccc}
\rhd_{\Gamma,\bar{\mathds{P}}}: & \D(\hb{\p}) \odo \He(\p) & \to & \He(\p)  \\
& a' \bicroi b' \odo a \# b & \mapsto & a'_{(3)}(b'_{(1)} \brhd a)S^{-1}_{A}(a'_{(2)}) \# (b'_{(2)}bS_{B}(b'_{(3)}) \blhd S^{-1}_{A}(a'_{(1)}))
\end{array}.
\]
\end{coro}
\pr
Follows directly by proposition~\ref{prop:taipe_duality_right_drinfeld}. For the induced action, using the equalities~(\ref{eq:alpha}) and~(\ref{eq:beta}), we have
\begin{align*}
(a' \bicroi b') \rhd_{\Gamma,\bar{\mathds{P}}} (a \# b) & = a' \rhd_{\alp,\hb{\p}} (b' \rhd_{\beta,\p} (a \# b)) = a' \rhd_{\alp,\hb{\p}} ((b'_{(1)} \brhd a) \# b'_{(2)}bS(b'_{(3)})) \\
& = a'_{(3)}(b'_{(1)} \brhd a)S^{-1}_{A}(a'_{(2)}) \# (b'_{(2)}bS_{B}(b'_{(3)}) \blhd S^{-1}_{A}(a'_{(1)})),
\end{align*}
for all $a,a'\in A$ and $b,b'\in B$.
\fin

We finish given the proof that the Heisenberg algebra $\He(\p)$ is a Yetter--Drinfeld $(\D(\hb{\p}),\com_{\D(\hb{\p})})$-algebra.

\begin{theo}
Consider the left coaction of $\D(\hb{\p})$ on $\He(\p)$,
\[
\begin{array}{lccc}
\Gamma' = (\id_{\D(\hb{\p})} \odo \bar{L})\com_{\D(\hb{\p})} \bar{L}^{-1} : & \He(\p) & \to & \M(\D(\hb{\p}) \odo \He(\p)) \\
& a \# b & \mapsto & a^{\co}_{(2)} \bicroi b_{(1)} \odo a_{(1)} \# b_{(2)}
\end{array}
\]
and the right coaction of $\T(\bar{\p})$ on $\He(\p)$ construct in Corollary~\ref{coro:hola},
\footnotesize
\[
\begin{array}{lccc}
\Gamma : & \He(\p) & \to & \M(\He(\p) \odo \T(\bar{\p})) \\
& a \# b & \mapsto & (\iota_{A}\iota_{B} \odo \id \odo \id)(U^{\ops}_{13}\Sigma(U)_{24})(a \# b \odo 1 \odo 1)(\iota_{A}\iota_{B} \odo \id \odo \id)(U^{\ops}_{13}\Sigma(U)_{24})^{-1}
\end{array}
\]
\normalsize
where $U \in \M(A \odo B)$ is the canonical multiplier associated to the admissible pairing $\p$. If $\bar{\mathds{P}} = \hb{\p} \odo \bar{\p}: \D(\hb{\p}) \times \T(\bar{\p}) \to \ku$ denotes the natural pairing between $\D(\hb{\p})$ and $\T(\bar{\p})$, then
\[
(\He(\p),\Gamma',\Gamma) \in \aYD^{lr}(\bar{\mathds{P}})
\]
or equivalently
\[
(\He(\p),\rhd_{\Gamma},\Gamma') \in {}^{\D(\hb{\p})}_{\D(\hb{\p})}\aYD.
\]
\end{theo}
\pr
Consider the map $\Theta = (\Sigma\Gamma' \odo \id)\Gamma: \He(\p) \to \M(\He(\p) \odo \D(\hb{\p})^{\co} \odo \T(\bar{\p}))$. The twisted tensor coproduct $(\D(\hb{\p})^{\co} \odo \T(\bar{\p}),\com_{\bar{\W}})$ associated to the skew-copairing $\Sigma(\bar{\W})^{-1} \in \M(\T(\bar{\p}) \odo \D(\hb{\p})^{\co})$ yields a regular multiplier Hopf algebra called the Majid's version of the Drinfeld codouble associate to the pairing $\bar{\mathds{P}}$, it will be denoted by $\T_{M}(\bar{\mathds{P}})$. Its comultiplication is given by the map
\[
\com_{\bar{\W}} = (\id_{\D(\hb{\p})} \odo \Sigma\ad(\bar{\W})^{-1} \odo \id_{\T(\p)})(\com^{\co}_{\D} \odo \com_{\T}).
\]
Using the left-right version of Theorem~\ref{prop:left_equivalence_taipe}, we have the equivalence of the following statements:
\begin{enumerate}[label=\textup{(\roman*)}]
\item $\Theta = (\Sigma\Gamma' \odo \id)\Gamma$ is a right coaction of $\T_{M}(\bar{\mathds{P}})$ on $\He(\p)$;
\item $(\He(\p),\Gamma',\Gamma)$ is a left-right Yetter--Drinfeld algebra over the pairing $\bar{\mathds{P}}$;
\item $(\He(\p),\rhd_{\Gamma},\Gamma') \in {}^{\D(\hb{\p})}_{\D(\hb{\p})}\aYD$.
\end{enumerate}

The last equivalence show that in order to obtain the result, we need to proof the $\Theta$ is a right coaction. Denote $X = (\bar{L} \odo \id)(\bar{\W}) = \bar{\V}_{12}\bar{\U}_{13} \in \M(\He(\p) \odo \T(\bar{\p}))$, thus by definition $\Gamma(x) = \ad(X)(x \odo 1)$ for all $x \in \He(\p)$. By direct computation it hold
\small
\begin{align*}
(\Sigma\Gamma' \odo \id)(X) & = (\Sigma \odo \id)((\id \odo \bar{L})\com_{\D}\bar{L}^{-1} \odo \id)(\bar{\V}_{12}\bar{\U}_{13}) \\
& = (\Sigma \odo \id)((\id \odo \bar{L})\com_{\D}\bar{L}^{-1} \odo \id)(\bar{L} \odo \id)(\bar{\W}) \\
& = (\Sigma \odo \id)(\id \odo \bar{L} \odo \id)(\com_{\D} \odo \id)(\bar{\W}) \\
& = (\Sigma \odo \id)(\id \odo \bar{L} \odo \id)(\bar{\W}_{13}\bar{\W}_{23}) = \bar{\W}_{23}X_{13}, \\
& \\
(\id \odo \com_{\bar{\W}})(\bar{\W}_{23}) & = 1 \odo \com_{W}(\bar{\W}) = 1 \odo (\id \odo \Sigma\ad(\bar{\W}^{-1}) \odo \id)(\com^{\co}_{\D} \odo \com)(\bar{\W}) \\
& = 1 \odo (\id \odo \Sigma\ad(\bar{\W}^{-1}) \odo \id)(\com^{\co}_{\D} \odo \id \odo \id)(\bar{\W}_{12}\bar{\W}_{13}) \\
& = 1 \odo (\id \odo \Sigma\ad(\bar{\W}^{-1}) \odo \id)((\com^{\co}_{\D} \odo \id)(\bar{\W})_{123}(\com^{\co}_{\D} \odo \id)(\bar{\W})_{124}) \\
& = 1 \odo (\id \odo \Sigma\ad(\bar{\W}^{-1}) \odo \id)((\bar{\W}_{23}\bar{\W}_{13})_{123}(\bar{\W}_{23}\bar{\W}_{13})_{124}) \\
& = 1 \odo \Sigma_{23}(\bar{\W}^{-1}_{23}\bar{\W}_{23}\bar{\W}_{13}\bar{\W}_{24}\bar{\W}_{14}\bar{\W}_{23}) = 1 \odo \bar{\W}_{12}\bar{\W}_{34}\bar{\W}_{14}\bar{\W}_{32} \\
& = \bar{\W}_{23}\bar{\W}_{45}\bar{\W}_{25}\bar{\W}_{43}
\end{align*}
and
\begin{align*}
(\id \odo \com_{\bar{\W}})(X_{13}) & = (\id \odo \id \odo \Sigma\ad(\bar{\W}^{-1}) \odo \id)(\id \odo \com^{\co}_{\D} \odo \com)(X_{13}) \\
& = (\id \odo \id \odo \Sigma\ad(\bar{\W}^{-1}) \odo \id)((\id \odo \com)(X)_{145}) \\ 
& = (\id \odo \id \odo \Sigma\ad(\bar{\W}^{-1}) \odo \id)((\id \odo \com)((\bar{L} \odo \id)(\bar{\W}))_{145}) \\
& = (\id \odo \id \odo \Sigma\ad(\bar{\W}^{-1}) \odo \id)(((\bar{L} \odo \id)(\bar{\W})_{12}(\bar{L} \odo \id)(\bar{\W})_{13})_{145}) \\
& = (\id \odo \id \odo \Sigma\ad(\bar{\W}^{-1}) \odo \id)((\bar{L} \odo \id)(\bar{\W})_{14}(\bar{L} \odo \id)(\bar{\W})_{15}) \\
& = (\bar{L} \odo \id \odo \id \odo \id \odo \id)\Sigma_{34}(\bar{\W}^{-1}_{34}\bar{\W}_{14}\bar{\W}_{15}\bar{\W}_{34}) \\
& = (\bar{L} \odo \id \odo \id \odo \id \odo \id)(\bar{\W}^{-1}_{43}\bar{\W}_{13}\bar{\W}_{15}\bar{\W}_{43}) = \bar{\W}^{-1}_{43}X_{13}X_{15}\bar{\W}_{43}.
\end{align*}

Now, we are ready to finish the proof. Fix $x \in \He(\p)$, on one hand we have
\small
\begin{align*}
(\id \odo \com_{\bar{\W}})\Theta(x) & = (\id \odo \com_{\bar{\W}})(\Sigma\Gamma' \odo \id)\Gamma(x) \\
& = (\id \odo \com_{\bar{\W}})(\Sigma \odo \id)((\id\odo \bar{L})\com_{\D}\bar{L}^{-1} \odo \id)\ad(\bar{\V}_{12}\bar{\U}_{13})(x \odo 1_{\T(\bar{\p})}) \\
& = (\id \odo \com_{\bar{\W}})(\bar{\W}_{23}X_{13}(\Sigma\Gamma'(x) \odo 1_{\T(\bar{\p})})(\bar{\W}_{23}X_{13})^{-1}) \\
& = (\id \odo \com_{\bar{\W}})(\bar{\W}_{23}X_{13})(\id \odo \com_{\bar{\W}})(\Sigma\Gamma'(x) \odo 1)(\id \odo \com_{\bar{\W}})(\bar{\W}_{23}X_{13})^{-1} \\ 
& = \ad(\bar{\W}_{23}\bar{\W}_{45}\bar{\W}_{25}X_{13}X_{15}\bar{\W}_{43})((\id \odo \com_{\bar{\W}})(\Sigma\Gamma'(x) \odo 1))
\end{align*}
\normalsize
and on the other
\small
\begin{align*}
(\Theta \odo \id)\Theta(x) & = (\Theta \odo \id)(\bar{\W}_{23}X_{13}(\Sigma\Gamma'(x) \odo 1)(\bar{\W}_{23}X_{13})^{-1}) \\
& = (\bar{\W}_{23}X_{13})_{123}((\Sigma\Gamma' \odo \id \odo \id)(\bar{\W}_{23}X_{13}(\Sigma\Gamma'(x) \odo 1)(\bar{\W}_{23}X_{13})^{-1}))_{1245}(\bar{\W}_{23}X_{13})^{-1}_{123} \\
& = \bar{\W}_{23}X_{13}\bar{\W}_{45}((\Sigma\Gamma' \odo \id \odo \id)(X_{13}(\Sigma\Gamma'(x) \odo 1)X^{-1}_{13}))_{1245}\bar{\W}^{-1}_{45}X^{-1}_{13}\bar{\W}^{-1}_{23} \\
& = \bar{\W}_{23}X_{13}\bar{\W}_{45}(\Sigma\Gamma' \odo \id)(X)_{125}((\Sigma\Gamma' \odo \id)\Sigma\Gamma'(x) \odo 1)_{1245}(\Sigma\Gamma' \odo \id)(X)^{-1}_{125}\bar{\W}^{-1}_{45}X^{-1}_{13}\bar{\W}^{-1}_{23} \\
& = \bar{\W}_{23}X_{13}\bar{\W}_{45}(\Sigma\Gamma' \odo \id)(X)_{125}((\Sigma\Gamma' \odo \id)\Sigma\Gamma'(x))_{124}(\Sigma\Gamma' \odo \id)(X)^{-1}_{125}\bar{\W}^{-1}_{45}X^{-1}_{13}\bar{\W}^{-1}_{23} \\
& = \bar{\W}_{23}X_{13}\bar{\W}_{45}(\Sigma\Gamma' \odo \id)(X)_{125}((\id \odo \com^{\co}_{\D})\Sigma\Gamma'(x))_{124}(\Sigma\Gamma' \odo \id)(X)^{-1}_{125}\bar{\W}^{-1}_{45}X^{-1}_{13}\bar{\W}^{-1}_{23} \\
& = \bar{\W}_{23}X_{13}\bar{\W}_{45}\bar{\W}_{25}X_{15}((\id \odo \com^{\co}_{\D})\Sigma\Gamma'(x))_{124}X^{-1}_{15}\bar{\W}^{-1}_{25}\bar{\W}^{-1}_{45}X^{-1}_{13}\bar{\W}^{-1}_{23} \\
& = \ad(\bar{\W}_{23}X_{13}\bar{\W}_{45}\bar{\W}_{25}X_{15})(((\id \odo \com^{\co}_{\D})\Sigma\Gamma'(x))_{124}) \\
& = \ad(\bar{\W}_{23}X_{13}\bar{\W}_{45}\bar{\W}_{25}X_{15})(((\id \odo \com^{\co}_{\D})\Sigma\Gamma'(x))_{124}).
\end{align*}
\normalsize
Follows from the equality
\begin{align*}
(\id \odo \com^{\co}_{\D})\Sigma\Gamma' & = (\id \odo \com^{\co}_{\D})\Sigma(\id \odo \bar{L})\com_{\D}\bar{L}^{-1} \\
& = (\id \odo \com^{\co}_{\D})(\bar{L} \odo \id)\com^{\co}_{\D}\bar{L}^{-1} \\
& = (\bar{L} \odo \id \odo \id)(\id \odo \com^{\co}_{\D})\com^{\co}_{\D}\bar{L}^{-1} \\
& = (\bar{L} \odo \id \odo \id)(\com^{\co}_{\D} \odo \id)\com^{\co}_{\D}\bar{L}^{-1},
\end{align*}
that
\small
\begin{align*} (\id \odo \com_{\bar{\W}})(\Sigma\Gamma'(x) \odo 1) & = (\id \odo \id \odo \Sigma\ad(\bar{\W}^{-1}) \odo \id)(\id \odo \com^{\co}_{\D} \odo \com)(\Sigma\Gamma'(x) \odo 1) \\
& = (\id \odo \id \odo \Sigma\ad(\bar{\W}^{-1}) \odo \id)(\id \odo \com^{\co}_{\D} \odo \com)((\bar{L} \odo \id)\com^{\co}_{\D}\bar{L}^{-1} \odo \id)(x \odo 1) \\
& = (\bar{L} \odo \id \odo \Sigma\ad(\bar{\W}^{-1}) \odo \id)(\com^{\co}_{\D} \odo \id \odo \com)(\com^{\co}_{\D}\bar{L}^{-1} \odo \id)(x \odo 1) \\
& = (\bar{L} \odo \id \odo \id \odo \id^{\odo 2})(\com^{\co}_{\D} \odo \id \odo \id^{\odo 2})\Sigma_{23}\ad(\bar{\W}^{-1})_{23}(\com^{\co}_{\D} \odo \id^{\odo 2})(\bar{L}^{-1}(x) \odo 1^{\odo 2}) \\
& = \Sigma_{34}(\bar{L} \odo \id \odo \id \odo \id^{\odo 2})(\com^{\co}_{\D} \odo \id \odo \id^{\odo 2})(\com^{\co}_{\D} \odo \id^{\odo 2})(\bar{L}^{-1}(x) \odo 1^{\odo 2}) \\
& = \Sigma_{34}((\id \odo \com^{\co}_{\D})\Sigma\Gamma'(x) \odo 1 \odo 1) = ((\id \odo \com^{\co}_{\D})\Sigma\Gamma'(x))_{124},
\end{align*}
\normalsize
and this implies that $(\Theta \odo \id)\Theta(x) = (\Theta \odo \id)\Theta(x)$. The last equality, the fact that $\Gamma$ and $\Gamma'$ are coactions of multiplier Hopf algebras and the injectivity of $\Theta$ implies that $\Theta$ is a right coaction of $\T_{M}(\bar{\mathds{P}})$ on $\He(\p)$. 

\fin

\section{Yetter--Drinfeld \as-algebras over algebraic quantum groups}\label{sec:yd-quantum}

An algebraic quantum group is defined as a multiplier Hopf \as-algebra endowed with an invariant integral. Additionally to the algebraic objects that arise from the theory of multiplier Hopf algebras, namely the antipode and the counit, algebraic quantum groups are endowed with other enriched objects of an analytical nature. Moreover, because any algebraic quantum group admits a canonical admissible pairing with its dual algebraic quantum group, in this section we will treat Yetter--Drinfeld structures over algebraic quantum groups as a special case of Yetter--Drinfeld structures given in the setting of multiplier Hopf algebras. This section intends to introduce definitions and notations that will be used in a further work \cite{Ta22_3}.

First, we briefly recall the main objects related to Van Daele's algebraic quantum groups.

\subsection{Algebraic quantum groups}

An {\em algebraic quantum group} $\G = (\pol(\G),\com,\var)$ consist of a multiplier Hopf \as-algebra $(\pol(\G),\com)$ endowed with a non-zero functional $\var: A \to \C$ satisfying the conditions
\begin{enumerate}[label=\textup{(\roman*)}]
\item $\var$ is {\em positive}, i.e. $\var(a^{*}a) \geq 0$ for all $a \in \pol(\G)$;
\item $\var$ is {\em faithful}, i.e. $\var(a^{*}a) = 0$ implies $a=0$;
\item $\var$ is {\em (left) invariant}, i.e. $(\id \odo \var)\com(a) = \var(a)1_{\pol(\G)}$ for all $a \in \pol(\G)$.
\end{enumerate}
Denote by $S$ and $\cou$ the antipode and counit of the multiplier Hopf \as-algebra $(\pol(\G),\com)$. In \cite{DCVD10}, it was proved that the {\em scaling constant of $\G$} is $1$, i.e. $\var\circ S^{2} = \var$. There is an analytic one-parameter group $\tau$ of automorphism of $A$, called the {\em scaling group of $\G$}, such that $S^{2} = \tau_{-i}$ and $\var\circ\tau_{z} = \tau_{z}$ for all $z \in \C$. There is an anti-\as-automorphism $R$ of $\pol(\G)$, called {\em the unitary antipode of $\G$}, such that $R\circ\tau_{z} = \tau_{z}\circ R$, $S=R\circ\tau_{-i/2}$, $\Sigma\circ(R \odo R)\circ\com = \com\circ R$ and $R^{2} =\id$. Sometimes, we use the notation $R_{\G}$ instead of $R$.

The {\em dual algebraic quantum group of $\G$} is given by $\dG = (\du{\pol(\G)},\dcom^{\co},\du{\psi}:=\dvar\circ\du{R})$, where $(\du{\pol(\G)},\dcom,\dvar)$ denotes the dual multiplier Hopf \as-algebra of $(\pol(\G),\com,\var)$ \cite{VD98} and $\du{R}$ is the unitary antipode of the algebraic quantum group $(\du{\pol(\G)},\dcom,\dvar)$. Recall that
\[
\du{\pol(\G)} = \{\varphi_{a} := \varphi(\cdot a): \pol(\G) \to \C : a \in \pol(\G) \} \subset \pol(\G)^{\vee}.
\]

The {\em opposite algebraic quantum group} of $\G$ is given by $\G^{\ops} := (\pol(\G),\com^{\co},\var\circ R)$ and the {\em conjugate algebraic quantum group} of $\G$ is given by $\G^{\con} := (\pol(\G)^{\op},\com^{\op},\var^{\op})$. Here $\pol(\G)^{\op}$ denote the opposite algebra of $\pol(\G)$, $\com^{\op}(a^{\op})(b^{\op} \odo c^{\op}) := ((b \odo c)\com(a))^{\op}$ and $\var^{\op}(a^{\op}):=\var(a)$, for every $a,b,c \in \pol(\G)$. We have the relations $\du{(\G^{\ops})} = (\dG)^{\con}$, $\du{(\G^{\con})} = (\dG)^{\ops}$ and $(\G^{\con})^{\ops} = (\G^{\ops})^{\con} \cong \G$.

A {\em left action of $\G$ on a \as-algebra $X$} is a right coaction of the multiplier Hopf \as-algebra $(\pol(\G),\com)$ on $X$. Similarly, a {\em right action of $\G$ on a \as-algebra $X$} is a left coaction of the multiplier Hopf \as-algebra $(\pol(\G),\com)$ on $X$.

The {\em Heisenberg \as-algebra of $\G$} is defined as the smash product \as-algebra $\du{\pol(\G)}\#\pol(\G)$ constructed using the canonical left action of the multiplier Hopf \as-algebra $(\pol(\G),\com_{\G})$ on its dual multiplier Hopf \as-algebra $(\du{\pol(\G)},\dcom_{\G})$, i.e. the left unital action $\brhd : \pol(\G) \odo \du{\pol(\G)} \to \du{\pol(\G)}$, $a \odo \var_{b} \mapsto \var_{ab}$. We will use the notation $\He(\G)$ for this \as-algebra.

The {\em algebraic multiplicative unitary of $\G$}, denoted by $U$, is the unique unitary element in the \as-algebra $\M(\pol(\G) \odo \pol(\dG))$ satisfying the conditions
\begin{enumerate}
\item $(\com_{\G} \odo \id)(U) = U_{13}U_{23}$, $(\id \odo \com_{\dG})(U) = U_{13}U_{12}$;
\item $(R_{\G} \odo \id)(U) = (\id \odo R_{\dG})(U) = U^{*}$;
\item $\com_{\G}(a) = U^{*}(a \odo 1)U$ in $\M(\pol(\G) \odo \He(\G))$ and $\com_{\dG}(\ome) = \Sigma(U)(1 \odo \ome)\Sigma(U^{*})$ in $\M(\pol(\dG) \odo \He(\G))$.
\end{enumerate}
Moreover, $U$ is the canonical multiplicative for the canonical pairing between the multiplier Hopf \as-algebras $(\pol(\G),\com)$ and $(\pol(\dG),\dcom)$. Observe that $\hat{U}:= \Sigma(U^{*}) \in \M(\pol(\dG) \odo \pol(\G))$ is the algebraic multiplicative unitary of $\dG$.

\subsection{Yetter--Drinfeld $\G$-\as-algebras}

Let $\G = (\pol(\G),\com,\varphi)$ be an algebraic quantum group and $U$ be the algebraic multiplicative unitary of $\G$.

\begin{defi}\label{df:yd-quantum}
Let $X$ be a \as-algebra, $\theta : X \to \M(X \odo \pol(\G))$ be a left action of $\G$ on $X$ and $\du{\theta}: X \to \M(X \odo \pol(\dG))$ be a left action of $\dG$ on $X$. We say that $(X,\theta,\du{\theta})$ is a {\em (left) Yetter--Drinfeld $\G$-\as-algebra} if
\begin{equation*}
(\theta \odo \id_{\pol(\G))})\circ\du{\theta} = (\id_{\M(X)} \odo \Sigma)\circ(\id_{\M(X)} \odo \ad(\Sigma U))\circ(\du{\theta} \odo \id_{\pol(\G)})\circ\theta.
\end{equation*}
A morphism of (left) Yetter--Drinfeld $\G$-\as-algebras $f: (X,\theta,\du{\theta}) \to (X',\theta',\du{\theta}')$ is a non-degenerate \as-homomorphism $f: X \to \M(X')$ such that $(f \odo \id)\circ\theta = \theta'\circ f$ and $(f \odo \id)\circ\du{\theta} = \du{\theta}'\circ f$. The category of Yetter--Drinfeld $\G$-\as-algebras will be denoted by ${}_{\G}\aYD$.

Similarly, if $\theta : X \to \M(\pol(\G) \odo X)$ is a right action of $\G$ on a \as-algebra $X$ and $\du{\theta}: X \to \M(\pol(\dG) \odo X)$ is a right action of $\dG$ on $X$, we say that $(X,\theta,\du{\theta})$ is a {\em (right) Yetter--Drinfeld $\G$-\as-algebra} if
\begin{equation*}
(\id_{\pol(\dG)} \odo \theta)\circ\du{\theta} = (\Sigma \odo \id_{\M(X)})\circ(\ad(U) \odo \id_{\M(X)})\circ(\id_{\pol(\G)} \odo \du{\theta})\circ\theta.
\end{equation*}
A morphism of (right) Yetter--Drinfeld $\G$-\as-algebras $f: (X,\theta,\du{\theta}) \to (X',\theta',\du{\theta}')$ is a non-degenerate \as-homomorphism $f: X \to \M(X')$ such that $(\id \odo f)\circ\theta = \theta'\circ f$ and $(\id \odo f)\circ\du{\theta} = \du{\theta}'\circ f$. The category of right Yetter--Drinfeld $\G$-\as-algebras will be denoted by $\aYD_{\G}$.
\end{defi}

\begin{nota}
Let $\gamma$ be and algebra automorphism  on a \as-algebra $X$ such that $\gamma\circ * \circ \gamma \circ * = \id_{X}$. Given a left action $\theta: X \to X \odo \pol(\G)$ of the algebraic quantum group $\G$ on $X$ such that $\theta\circ\gamma = (\gamma \odo S^{-2}_{\G})\circ\theta$.
\begin{enumerate}[label=\textup{(\roman*)}]
\item We will use the notation $\theta^{\ops}$ for the co-opposite coaction $\theta^{\co}: X^{\op} \to X^{\op} \odo \pol(\G)^{\co}$. This is a left action of the algebraic quantum group $\G^{\ops}$ on $X$.

\item We will use the notation $\theta^{\con}$ for the opposite coaction $\theta^{\op}: X^{\op} \to X^{\op} \odo \pol(\G)^{\op}$. This is a left action of the algebraic quantum group $\G^{\con}$ on $X$.
\end{enumerate}
Similarly notations will be used for right actions of algebraic quantum groups.
\end{nota}

\begin{defi}\label{df:bc-quantum}
A (left) Yetter--Drinfeld $\G$-\as-algebra $(X,\theta,\du{\theta})$ is called {\em braided commutative} if for each $x,y \in X$, we have
\small
\begin{equation*}\label{eq:l-bc-quantum}
\theta^{\con}(x^{\op})\du{\theta^{\ops}}(y^{\op}) = \du{\theta}^{\ops}(y^{\op})\theta^{\con}(x^{\op})
\end{equation*}
\normalsize
inside the \as-algebra $\M(X^{\op} \odo \He(\dG^{\ops}))$. We will denote by ${}_{\G}\aYD^{\textrm{bc}}$, the subcategory of braided commutative left Yetter--Drinfeld $\G$-\as-algebras. Similarly, a (right) Yetter--Drinfeld $\G$-\as-algebra $(X,\theta,\du{\theta})$ is called {\em braided commutative} if for each $x,y \in X$, we have
\small
\begin{equation*}\label{eq:r-bc-quantum}
\theta^{\con}(x^{\op})\du{\theta^{\ops}}(y^{\op}) = \du{\theta}^{\ops}(y^{\op})\theta^{\con}(x^{\op})
\end{equation*}
\normalsize
inside the \as-algebra $\M(\He(\G) \odo X^{\op})$. We will denote by $\aYD^{\textrm{bc}}_{\G}$, the subcategory of braided commutative right $\G$-Yetter--Drinfeld \as-algebras.
\end{defi}

\begin{rema}
Let $\p_{\G}:\pol(\G)^{\op} \times \pol(\dG)^{\co} \to \C$ be the admissible pairing $\p_{\G}(a^{\op},\varphi_{b}) = \varphi(ab)$ for all $a,b \in \pol(\G)$. Our definitions of (braided commutative) Yetter--Drinfeld \as-algebras in the setting of algebraic quantum groups are given in order to have the following equivalence:
\begin{enumerate}[label=\textup{(\roman*)}]
\item $(X,\theta,\du{\theta})$ is in (resp. ${}_{\G}\aYD^{\textrm{bc}}$) ${}_{\G}\aYD$ if and only if $(X,\theta,\du{\theta})$ is in (resp. $\aYD^{\,rr}_{\textrm{bc}}(\p_{\G})$) $\aYD^{rr}(\p_{\G})$.
\item $(X,\theta,\du{\theta})$ is in (resp. $\aYD^{\textrm{bc}}_{\G}$) $\aYD_{\G}$ if and only if $(X,\theta,\du{\theta})$ is in (resp. $\aYD^{\,ll}_{\textrm{bc}}(\p_{\G})$) $\aYD^{ll}(\p_{\G})$;
\end{enumerate}
Moreover, our definitions are compatible with those ones given in the framework of von Neumann locally compact quantum groups \cite{ET16}.
\end{rema}



\bibliographystyle{abbrv}
\bibliography{biblio}

\addresseshere


\newpage
\appendix

\section{Standard Yetter--Drinfeld algebras}

For sake of completion, in this appendix section we collect the results on the other versions of Yetter--Drinfeld algebras over multiplier Hopf algebras. The details are similar to the results in the Hopf algebra setting, see for example \cite{Ta22_0}.

\subsection{Standard left-left Yetter--Drinfeld algebras}

If $X$ is a non-degenerate $(A,\com)$-algebra for the left action $\rhd: A \odo X \to X$ and a $(A,\com)$-algebra for the left coaction $\Gamma: X \to \M(A \odo X)$. Consider the Sweedler type leg notation $x_{[-1]} \odo x_{[0]} := \Gamma(x)$ for all $x \in X$.

We say that $(X,\rhd,\Gamma)$ is a {\em (left-left) Yetter--Drinfeld $(A,\com)$-algebra} if we have the following condition
\begin{equation}\label{eq:ll-YD-s-1}\tag{s-ll-YD}
a_{(1)}x_{[-1]}a' \odo (a_{(2)} \rhd x_{[0]}) = (a_{(1)} \rhd x)_{[-1]}a_{(2)}a' \odo (a_{(1)} \rhd x)_{[0]}
\end{equation}
for all $x \in X$ and $a,a' \in A$. A (left-left) Yetter--Drinfeld $(A,\com)$-algebra $(X,\rhd,\Gamma)$ is called {\em braided commutative} or {\em $A$-commutative}, if
\begin{equation}\label{eq:ll-bc-s-1}\tag{s-ll-BC}
xy = (x_{[-1]} \rhd y)x_{[0]} \quad \text{or equivalently} \quad xy = y_{[0]}(S^{-1}(y_{[-1]}) \rhd x)
\end{equation}
for all $x, y \in X$.


\begin{prop}
The tuple $(X,\rhd,\Gamma)$ is a left-left Yetter--Drinfeld $(A,\com)$-algebra if and only if
\begin{equation}\label{eq:ll-YD-s-2}\tag{s-ll-YD'}
\Gamma(a \rhd x)(a' \odo 1) = a_{(1)}x_{[-1]}S(a_{(3)})a' \odo (a_{(2)} \rhd x_{[0]})
\end{equation}
for all $x \in X$, $a, a' \in A$.
\end{prop}

\subsection{Standard left-right Yetter--Drinfeld algebras}

Let $X$ be a non-degenerate $(A,\com)$-module algebra for the left action $\rhd : A \odo X \to X$ and a $(A,\com)$-comodule algebra for the right coaction $\Gamma: X \to \M(X \odo A^{\op})$. Consider the Sweedler type leg notation $x_{[0]} \odo x^{\op}_{[1]} := \Gamma(x)$ for all $x \in X$.

We say that $(X,\rhd,\Gamma)$ is a {\em (left-right) Yetter--Drinfeld $(A,\com)$-algebra} if 
\begin{equation}\tag{s-lr-YD}
(a_{(1)} \rhd x_{[0]}) \odo a_{(2)}x_{[1]}a' = (a_{(2)} \rhd x)_{[0]} \odo (a_{(2)} \rhd x)_{[1]}a_{(1)}a'
\end{equation}
for all $x \in X$, $a,a' \in A$. A (left-right) Yetter--Drinfeld $(A,\com)$-algebra $(X,\rhd,\Gamma)$ is called {\em braided commutative} or {\em $A$-commutative}, if
\begin{equation}\label{eq:lr-bc-s-1}\tag{s-lr-BC}
xy = (S(x_{[1]}) \rhd y)x_{[0]} \quad \text{or equivalently} \quad xy = y_{[0]}(y_{[1]} \rhd x)
\end{equation}
for all $x, y \in X$.


\begin{prop}
The tuple $(X,\rhd,\Gamma)$ is a (left-right) Yetter--Drinfeld $(A,\com)$-algebra if and only if
\begin{equation}\tag{s-lr-YD'}
\Gamma(a \rhd x)(1 \odo a'^{\op}) = (a_{(2)} \rhd x_{[0]}) \odo (a'a_{(3)}x_{[1]}S^{-1}(a_{(1)}))^{\op}
\end{equation}
for all $x \in X$, $a,a' \in A$.
\end{prop}

\subsection{Standard right-right Yetter--Drinfeld algebras}

Let $X$ be a non-degenerate $(A,\com)$-module algebra for the left action $\lhd : X \odo A \to X$ and a $(A,\com)$-comodule algebra for the right coaction $\Gamma: X \to \M(X \odo A)$. Consider the Sweedler type leg notation $x_{[0]} \odo x_{[1]} := \Gamma(x)$ for all $x \in X$. 

We say that $(X,\lhd,\Gamma)$ is a {\em (right-right) Yetter--Drinfeld $(A,\com)$-algebra} if 
\begin{equation}\tag{s-rr-YD}
(x_{[0]} \lhd a_{(1)}) \odo a'x_{[1]}a_{(2)} = (x \lhd a_{(2)})_{[0]} \odo a'a_{(1)}(x \lhd a_{(2)})_{[1]}
\end{equation}
for all $x \in X$, $a,a' \in A$. A (right-right) Yetter--Drinfeld $(A,\com)$-algebra $(X,\lhd,\Gamma)$ is called {\em braided commutative} or {\em $A$-commutative}, if
\begin{equation}\label{eq:rr-bc-s-1}\tag{s-rr-BC}
xy = (y \lhd S^{-1}(x_{[1]}))x_{[0]} \quad \text{or equivalently} \quad xy = y_{[0]}(x \lhd y_{[1]})
\end{equation}
for all $x, y \in X$.

\begin{prop}
The tuple $(X,\lhd,\Gamma)$ is a (right-right) Yetter--Drinfeld $(A,\com)$-algebra if and only if
\begin{equation}\tag{s-rr-YD'}
(1 \odo a')\Gamma(x \lhd a) = (x_{[0]} \lhd a_{(2)}) \odo a'S(a_{(1)})x_{[1]}a_{(3)}
\end{equation}
for all $x \in X$, $a,a' \in A$.
\end{prop}

\section{Yetter--Drinfeld algebras over a pairing}

For sake of completion, in this appendix section we collect the results on the other versions of Yetter--Drinfeld algebras over pairings. The reader interested can consult \cite{Ta22_0} for the details in the Hopf algebra setting. 

\subsection{Right-right Yetter--Drinfeld algebras}

Let $\alp : X \to \M(X \odo A^{\op})$ be a right coaction of $(A^{\op},\com^{\op}_{A})$ on $X$ and $\beta : X \to \M(X \odo B)$ be a right coaction of $(B,\com_{B})$ on $X$. We say that $(X,\alp,\beta)$ is a {\em right-right Yetter--Drinfeld algebra over $\p$} if
\begin{equation}\tag{rr-YD}
(\alp \odo \id_{\M(B)})\beta = \ad(U^{\ops}_{23})^{-1}\Sigma_{23}(\beta \odo \id_{\M(A^{\op})})\alp.
\end{equation}
Sometimes, we will write $(X,\alp,\beta) \in \aYD^{\,rr}(\p)$. A right-right Yetter--Drinfeld algebra $(X,\alp,\beta)$ over $\p$ will be called {\em braided commutative} if
\small
\begin{equation}\label{eq:rr-bc-mod}\tag{rr-BC}
\m_{X^{\op} \odo \He(\bar{\p})}(\hat\iota_{X^{\op},A} \odo \hat\iota_{X^{\op},B})(\alp^{\op}(x^{\op}) \odo \beta^{\co}(y^{\op})) = \m_{X^{\op} \odo \He(\bar{\p})}(\hat\iota_{X^{\op},B} \odo \hat\iota_{X^{\op},A})(\beta^{\co}(y^{\op}) \odo \alp^{\op}(x^{\op}))
\end{equation}
\normalsize
for each $x,y \in X$. In this case, we will denote $(X,\alp,\beta) \in \aYD^{\,rr}_{bc}(\p)$.

\begin{prop}\label{prop:dual-rr-bc-mod}
A right-right Yetter--Drinfeld algebra $(X,\alp,\beta)$ over $\p$ is braided commutative if and only if
\[
\m_{X \odo \He(\p)}(\du{\iota}_{X,A} \odo \du{\iota}_{X,B})((\id \odo S_{A})\dot\alp(x) \odo \beta(y)) = \m_{X \odo \He(\p)}(\du{\iota}_{X,B} \odo \du{\iota}_{X,A})(\beta(y) \odo (\id \odo S_{A})\dot\alp(x))
\]
for every $x,y \in X$.
\end{prop}

\begin{prop}\label{prop:equivalence_st_oc_rr}
The following statements are equivalent:
\begin{enumerate}[label=\textup{(\arabic*)}]
\item $(X,\alp,\beta)$ is a (resp. braided commutative) right-right Yetter--Drinfeld algebra over $\p$.
\item $(X,\rhd_{\beta,\p},\alp)$ is a (resp. braided commutative) left-right Yetter--Drinfeld $(A,\com_{A})$-algebra.
\end{enumerate}
Here, the map $\rhd_{\beta,\p} : A \odo X \to X$, $x \odo a \mapsto (\id \odo {}_{a}\p)(\beta(x))$ denotes the left unital action of $(A,\com_{A})$ on $X$ induced by the right coaction $\beta$ and the pairing $\p$.
\end{prop}

\subsection{Left-right Yetter--Drinfeld algebras}

Let $\alp : X \to \M(A \odo X)$ be a left coaction of $(A,\com_{A})$ and $\beta : X \to \M(X \odo B)$ be a right coaction of $(B,\com_{B})$ on the same algebra $X$. We say that $(X,\alp,\beta)$ is a {\em left-right Yetter--Drinfeld algebra over $\p$} if
\begin{equation}\tag{lr-YD}
(\alp \odo \id_{\M(B)})\beta = \ad(U(\p)_{13})(\id_{\M(A)} \odo \beta)\alp.
\end{equation}
We will write $(X,\alpha,\beta) \in \aYD^{\,lr}(\p)$. A left-right Yetter--Drinfeld algebra $(X,\alpha,\beta)$ over $\p$ will be called {\em braided commutative} if
\begin{equation}\label{eq:lr-bc-mod}\tag{lr-BC}
\m_{\He(\p) \odo X}(\iota_{A,X} \odo \iota_{B,X})(\alp(x) \odo \Sigma\beta(y)) = \m_{\He(\p) \odo X}(\iota_{B,X} \odo \iota_{A,X})(\Sigma\beta(y) \odo \alp(x))
\end{equation}
for each $x,y \in X$. In this case, we will denote $(X,\alpha,\beta) \in \aYD^{\,lr}_{bc}(\p)$.

\begin{prop}\label{prop:equivalence_st_oc_lr}
The following statements are equivalent:
\begin{enumerate}[label=\textup{(\arabic*)}]
\item $(X,\alp,\beta)$ is a (resp. braided commutative) left-right Yetter--Drinfeld algebra over $\p$;
\item $(X,\rhd_{\beta,\p},\alp)$ is a (resp. braided commutative) left-left Yetter--Drinfeld $(A,\com_{A})$-algebra.
\end{enumerate}
Here, $\rhd_{\beta,\p}: A \odo X \to X$, $a \odo x \mapsto (\id \odo {}_{a}\p)(\beta(x))$ %
denotes the left unital action of $(A,\com_{A})$ on $X$ induced by the right coaction $\beta$ and the pairing $\p$.
\end{prop}

\subsection{Right-left Yetter--Drinfeld algebras}

Let $\alp : X \to \M(X \odo A)$ be a right coaction of $(A,\com_{A})$ and $\beta : X \to \M(B \odo X)$ be a left coaction of $(B,\com_{B})$ on the same algebra $X$. We say that $(X,\alp,\beta)$ is a {\em right-left Yetter--Drinfeld algebra over $\p$} if
\begin{equation}\tag{rl-YD}
(\id_{\M(B)} \odo \alp)\beta = \ad(U(\p)^{-1}_{13})(\beta \odo \id_{\M(A)})\alp.
\end{equation}
We will write $(X,\alpha,\beta) \in \aYD^{\,rl}(\p)$. A right-left Yetter--Drinfeld algebra $(X,\alpha,\beta)$ over $\p$ will be called {\em braided commutative} if
\begin{equation}\label{eq:rl-bc-mod}\tag{rl-BC}
\m_{X \odo \He(\bar{\p})}(\hat\iota_{X,A} \odo \hat\iota_{X,B})(\alp(x) \odo \Sigma\beta(y)) = \m_{X \odo \He(\bar{\p})}(\hat\iota_{X,B} \odo \hat\iota_{X,A})(\Sigma\beta(y) \odo \alp(x))
\end{equation}
for each $x,y \in X$. In this case, we will denote $(X,\alpha,\beta) \in \aYD^{\,rl}_{bc}(\p)$.

\begin{prop}\label{prop:equivalence_st_oc_lr}
The following statements are equivalent:
\begin{enumerate}[label=\textup{(\arabic*)}]
\item $(X,\alp,\beta)$ is a (resp. braided commutative) right-left Yetter--Drinfeld algebra over $\p$;
\item $(X,\lhd_{\beta,\p},\alp)$ is a (resp. braided commutative) right-right Yetter--Drinfeld $(A,\com_{A})$-algebra.
\end{enumerate}
Here, $\lhd_{\beta,\p}: X \odo A \to X$, $x \odo a \mapsto ({}_{a}\p \odo \id)(\beta(x))$ denotes the right unital action of $(A,\com_{A})$ on $X$ induced by the left coaction $\beta$ and the pairing $\p$.
\end{prop}

\subsection{Equivalence between Yetter--Drinfeld categories}

We end this appendix section with an equivalence between Yetter--Drinfeld categories.
\begin{prop}
Let $\p: A \times B \to \ku$ be an admissible pairing between regular multiplier Hopf algebras. We have the following equivalence of categories
\begin{equation*}\label{eq:left_equivalence_all_yd_m}
\xymatrix@C=5pc@R=1.5pc{%
\shadowbox*{%
\begin{Bcenter}
$\;\aYD^{\,rl}(\p)\;$ \\[0.1cm]
$(X,\alp,\beta)$
\end{Bcenter}}\ar@{<~>}[r]\ar@{<~>}[d]%
&%
\shadowbox*{%
\begin{Bcenter}
$\;\aYD^{\,lr}(\bar{\p})\;$ \\[0.1cm]
$(X,\beta,\alp)$
\end{Bcenter}}\ar@{<~>}[d]%
\\
{\shadowbox*{%
\begin{Bcenter}
$\;\aYD^{\,ll}(\p)\;$ \\[0.1cm]
$(X,(S^{-1}_{A} \odo \id)\Sigma\alpha,\beta)$
\end{Bcenter}}}\ar@{<~>}[r]%
&%
\shadowbox*{%
\begin{Bcenter}
$\;\aYD^{\,rr}(\bar{\p})\;$ \\[0.1cm]
$(X,(\id \odo S^{-1}_{B})\Sigma\beta,\alpha)$
\end{Bcenter}}}
\end{equation*}
Moreover, the above equivalence respect all the  braided commutativity conditions.
\end{prop}
\pr
It is straightforward.
\fin

\end{document}